\documentstyle{amltd2004}
\input amssym.def
\input amssym.tex
\input pic.mac
\def\hiha{\hangindent=27pt\hangafter=1}
\input boxedeps.tex 
\SetepsfEPSFSpecial 
\HideDisplacementBoxes
\def\figin#1#2{ 
$$
 {\BoxedEPSF{#1.eps scaled
#2}%
}%
$$
\noindent}
\def\joinrel{\mathrel{\mkern-4mu}}
\def\relbar{\mathrel{\smash-}}
\def\lrar{\relbar\joinrel\relbar\joinrel\rightarrow}

\def\llar{\leftarrow\joinrel\relbar\joinrel\relbar}
 
\begin{document}
\annalsline{155}{2002}
\received{May 21, 1999}
\startingpage{611}
\def\bye{\ \end{document}}
 \font\tenrm=cmr10
\def\ritem#1{\item[{\rm #1}]}
\catcode`\@=11
\font\twelvemsb=msbm10 scaled 1100
\font\tenmsb=msbm10
\font\ninemsb=msbm10 scaled 800
\newfam\msbfam
\textfont\msbfam=\twelvemsb  \scriptfont\msbfam=\ninemsb
  \scriptscriptfont\msbfam=\ninemsb
\def\msb@{\hexnumber@\msbfam}
\def\Bbb{\relax\ifmmode\let\next\Bbb@\else
 \def\next{\errmessage{Use \string\Bbb\space only in math
mode}}\fi\next}
\def\Bbb@#1{{\Bbb@@{#1}}}
\def\Bbb@@#1{\fam\msbfam#1}
\catcode`\@=12

 \catcode`\@=11
\font\twelveeuf=eufm10 scaled 1100
\font\teneuf=eufm10
\font\nineeuf=eufm7 scaled 1100
\newfam\euffam
\textfont\euffam=\twelveeuf  \scriptfont\euffam=\teneuf
  \scriptscriptfont\euffam=\nineeuf
\def\euf@{\hexnumber@\euffam}
\def\frak{\relax\ifmmode\let\next\frak@\else
 \def\next{\errmessage{Use \string\frak\space only in math
mode}}\fi\next}
\def\frak@#1{{\frak@@{#1}}}
\def\frak@@#1{\fam\euffam#1}
\catcode`\@=12


\newcommand{\inv}{^{-1}}
\def\myunderline#1{#1}

\newcommand{\bX}{{\Bbb X}}
\newcommand{\bL}{{\Bbb L}}
\newcommand{\bT}{{\Bbb T}}
\newcommand{\bP}{{\Bbb P}}
\newcommand{\bQ}{{\Bbb Q}}
\newcommand{\bZ}{{\Bbb Z}}
\newcommand{\bM}{{\Bbb M}}
\newcommand{\bN}{{\Bbb N}}
\newcommand{\bR}{{\Bbb R}}
\newcommand{\bC}{{\Bbb C}}
\newcommand{\bG}{{\Bbb G}}
\newcommand{\bA}{{\Bbb A}}
\newcommand{\bF}{{\Bbb F}}
\newcommand{\bV}{{\Bbb V}}
\newcommand{\bK}{{\Bbb K}}

\newcommand{\dX}{{\bf X}}
\newcommand{\dx}{{\bf x}}
\newcommand{\dy}{{\bf y}}
\newcommand{\dL}{{\bf L}}
\newcommand{\dT}{{\bf T}}
\newcommand{\dH}{{\bf H}}
\newcommand{\dP}{{\bf P}}
\newcommand{\dQ}{{\bf Q}}
\newcommand{\dZ}{{\bf Z}}
\newcommand{\dN}{{\bf N}}
\newcommand{\dR}{{\bf R}}
\newcommand{\dC}{{\bf C}}
\newcommand{\dG}{{\bf G}}
\newcommand{\dA}{{\bf A}}
\newcommand{\dV}{{\bf V}}
\newcommand{\dS}{{\bf S}}
\newcommand{\dK}{{\bf K}}

\newcommand{\dt}{{\bf t}}
\newcommand{\dtau}{{\bf \tau}}

\newcommand{\cC}{{\cal C}}
\newcommand{\cK}{{\cal K}}
\newcommand{\cR}{{\cal R}}
\newcommand{\cL}{{\cal L}}
\newcommand{\cG}{{\cal G}}
\newcommand{\cH}{{\cal H}}
\newcommand{\cM}{{\cal M}}
\newcommand{\cF}{{\cal F}}
\newcommand{\cP}{{\cal P}}
\newcommand{\cU}{{\cal U}}
\newcommand{\cE}{{\cal E}}
\newcommand{\cO}{{\cal O}}
\newcommand{\cA}{{\cal A}}
\newcommand{\cX}{{\cal X}}
\newcommand{\cN}{{\cal N}}
\newcommand{\cS}{{\cal S}}
\newcommand{\cV}{{\cal V}}
\newcommand{\cD}{{\cal D}}
\newcommand{\cT}{{\cal T}}

\newcommand{\fm}{\frak m}
\newcommand{\fB}{\frak B}
\newcommand{\fF}{\frak F}
\newcommand{\fR}{\frak R}

\newcommand{\hc}{\widehat c}
\newcommand{\hh}{\hat h}
\newcommand{\hs}{\hat s}
\newcommand{\hg}{\hat g}
\newcommand{\hf}{\hat f}
\newcommand{\hx}{\hat x}
\newcommand{\hl}{\hat l}
\newcommand{\hchi}{\hat \chi}

\newcommand{\hY}{\widehat Y}
\newcommand{\hT}{\widehat T}
\newcommand{\hF}{\widehat F}
\newcommand{\hX}{\widehat X}
\newcommand{\hU}{\widehat U}
\newcommand{\hM}{\widehat M}
\newcommand{\hP}{\widehat P}
\newcommand{\hG}{\widehat G}
\newcommand{\hE}{\widehat E}
\newcommand{\hH}{\widehat H}
\newcommand{\hR}{\widehat R}
\newcommand{\hL}{\widehat L}
\newcommand{\hK}{\widehat K}

\newcommand{\hone}{\widehat 1}
\newcommand{\hzero}{\widehat 0}
\newcommand{\hFun}{\widehat {{\rm Fun}}}

\newcommand{\hbK}{\widehat \bK}
\newcommand{\hbL}{\widehat \bL}
\newcommand{\hbX}{\widehat \bX}
\newcommand{\hbM}{\widehat \bM}

\newcommand{\tbM}{\widetilde \bM}
\newcommand{\tFun}{\widetilde {{\rm Fun}}}

\newcommand{\uA}{\underline{{\rm A}}}
\newcommand{\uP}{\underline{{\rm P}}}
\newcommand{\uR}{\underline{{\rm R}}}
\newcommand{\uS}{\underline{{\rm S}}}
\newcommand{\uG}{\underline{G}}
\newcommand{\uX}{\underline{X}}
\newcommand{\uY}{\underline{Y}}
\newcommand{\uT}{\underline{T}}

\newcommand{\uhX}{\underline{\widehat X}}
\newcommand{\uhK}{\underline{\widehat K}}

\newcommand{\ubG}{\underline{\bG}}
\newcommand{\ubX}{\underline{\bX}}
\newcommand{\ubK}{\underline{\bK}}
\newcommand{\ubF}{\underline{\bF}}
\newcommand{\ubT}{\underline{\bT}}
\newcommand{\ubL}{\underline{\bL}}

\newcommand{\uhbX}{\underline{\widehat \bX}}
\newcommand{\uhbL}{\underline{\widehat{\bL}}}
\newcommand{\uhbK}{\underline{\widehat{\bK}}}
\newcommand{\uhbF}{\underline{\widehat{\bF}}}

\newcommand{\wH}{\widetilde H}
\newcommand{\wG}{\widetilde G}
\newcommand{\wM}{\widetilde M}
\newcommand{\wU}{\widetilde U}
\newcommand{\wV}{\widetilde V}
\newcommand{\wP}{\widetilde P}
\newcommand{\wS}{\widetilde S}
\newcommand{\wT}{\widetilde T}
\newcommand{\wC}{\widetilde C}
\newcommand{\wR}{\widetilde R}
\newcommand{\wL}{\widetilde L}
\newcommand{\wK}{\widetilde K}
\newcommand{\wf}{\widetilde f}
\newcommand{\wGamma}{\widetilde\Gamma}
\newcommand{\wDelta}{\widetilde\Delta}
\newcommand{\Deltaw}{\widetilde\Delta}
\newcommand{\wdelta}{\widetilde\delta}
\newcommand{\deltaw}{\widetilde\delta}
\newcommand{\wphi}{\widetilde\phi}

\newcommand{\wdL}{\widetilde {\bf L}}

\newcommand{\ox}{\overline{x}}
\newcommand{\oy}{\overline{y}}
\newcommand{\oz}{\overline{z}}
\newcommand{\oC}{\overline{C}}
\newcommand{\oV}{\overline{V}}
\newcommand{\oW}{\overline{W}}
\newcommand{\oL}{\overline{L}}
\newcommand{\oX}{\overline{X}}
\newcommand{\oY}{\overline{Y}}
\newcommand{\oZ}{\overline{Z}}
\newcommand{\oP}{\overline{P}}
\newcommand{\oH}{\overline{H}}

\newcommand{\ocG}{\overline{\cG}}
\newcommand{\obM}{\overline{\Bbb M}}

\newcommand{\skipline}{\mbox{}\smallskip}

\newcommand{\isoto}{{\stackrel{\sim}{\rightarrow}}}
\newcommand{\ratmap}{- \kern -3pt \to}
\newcommand{\iiff}{\Leftrightarrow}
\newcommand{\acts}{\curvearrowright}
\newcommand{\onto}{\twoheadrightarrow}
\newcommand{\into}{\hookrightarrow}
\newcommand{\tendsto}{\rightsquigarrow}
\newcommand{\follows}{\Rightarrow}

\newcommand{\SP}{{{\rm Sp}}}
\newcommand{\diag}{{{\rm diag}}}
\newcommand{\Orb}{{{\rm Orb}}}
\newcommand{\Iso}{{{\rm Iso}}}
\newcommand{\LCM}{{{\rm LCM}}}
\newcommand{\GCD}{{{\rm GCD}}}
\newcommand{\emb}{{{\rm emb}}}
\newcommand{\ChC}{\hbox{\u{$C$}}}
\newcommand{\Star}{{{\rm Star}}}
\newcommand{\sn}{{{\rm sn}}}
\newcommand{\wn}{{{\rm wn}}}
\newcommand{\orb}{{{\rm orb}}}
\newcommand{\ord}{{{\rm ord}}}
\newcommand{\sat}{{{\rm sat}}}
\newcommand{\stab}{{{\rm stab}}}
\newcommand{\lin}{{{\rm lin}}}
\newcommand{\Fan}{{{\rm Fan}}}
\newcommand{\Cone}{{{\rm Cone}}}
\newcommand{\Conv}{{{\rm Conv}}}
\newcommand{\Fun}{{{\rm Fun}}}
\newcommand{\link}{{{\rm link}}}
\newcommand{\Sat}{{{\rm Sat}}}
\newcommand{\main}{{{\rm main}}}
\newcommand{\red}{{{\rm red}}}
\newcommand{\ch}{{{\rm ch}}}
\newcommand{\lev}{{{\rm lev}}}
\newcommand{\coker}{{{\rm coker}}}
\newcommand{\gps}{{{\rm gps}}}
\newcommand{\supp}{{{\rm supp}}}
\newcommand{\osupp}{{{\rm \overline{supp}}}}
\newcommand{\const}{{{\rm const}}}
\newcommand{\Vol}{{{\rm Vol}}}
\newcommand{\conv}{{{\rm conv}}}
\newcommand{\QC}{{{\rm QC}}}
\newcommand{\REP}{{{\rm REP}}}
\newcommand{\id}{{{\rm id}}}
\newcommand{\modulo}{\,{\rm mod}}
\newcommand{\Tor}{{{\rm Tor}}}
\newcommand{\Tors}{{{\rm Tors}}}
\newcommand{\Supp}{{{\rm Supp}}}
\newcommand{\GL}{{{\rm GL}}}
\newcommand{\Hilb}{{{\rm Hilb}}}
\newcommand{\Proj}{{{\rm Proj}}}
\newcommand{\bProj}{{{\rm \Bbb Proj}}}
\newcommand{\Stab}{{{\rm Stab}}}
\newcommand{\Spec}{{{\rm Spec}}}
\newcommand{\Diff}{{{\rm Diff}}}
\newcommand{\mult}{{{\rm mult} }}
\newcommand{\rank}{{{\rm rank}}}
\newcommand{\Hom}{{{\rm Hom}}}
\newcommand{\calHom}{\mathop{\cal Hom}}
\newcommand{\Coker}{{{\rm Coker}}}
\newcommand{\Ker}{{{\rm Ker}}}
\newcommand{\im}{{{\rm im}}}
\newcommand{\NS}{{{\rm NS}}}
\newcommand{\Pic}{{{\rm Pic}}}
\newcommand{\calPic}{{{\rm \cal Pic}}}
\newcommand{\Aut}{{{\rm Aut}}}
\newcommand{\calAut}{\mathop{\cal A\cal \cal t}}
\newcommand{\chr}{{{\rm char}}}
\newcommand{\Sym}{{{\rm Sym}}}
\newcommand{\Ext}{{{\rm Ext}}}
\newcommand{\Gal}{{{\rm Gal}}}
\newcommand{\Cor}{{{\rm Cor}}}
\newcommand{\Pol}{{{\rm Pol}}}

\newcommand{\bPol}{{{\rm \mathbf{Pol}}}}
\newcommand{\bCor}{{{\rm \mathbf{Cor}}}}
\newcommand{\bPic}{{{\rm \mathbf{Pic}}}}
\newcommand{\bPict}{{{\rm \mathbf{Pic}^{\tau}}}}
\newcommand{\bAut}{{{\rm \mathbf{Aut}}}}
\newcommand{\bExt}{{{\rm \mathbf{Ext}}}}
\newcommand{\bNS}{{{\rm \mathbf{NS}}}}
\newcommand{\bHom}{{{\rm Hom}}}

\newcommand{\uFun}{\underline{{\rm Fun}}}
\newcommand{\uPol}{\underline{{\rm Pol}}}
\newcommand{\uCor}{\underline{{\rm Cor}}}
\newcommand{\uPic}{\underline{{\rm Pic}}}
\newcommand{\uAut}{\underline{{\rm Aut}}}
\newcommand{\uExt}{\underline{{\rm Ext}}}
\newcommand{\uNS}{\underline{{\rm NS}}}
\newcommand{\uHom}{\underline{{\rm Hom}}}

\newcommand{\uhFun}{\underline{{\rm \widehat{Fun}}}}
\newcommand{\utFun}{\underline{{\rm \tilde{Fun}}}}


\renewcommand{\ChC}{C}

\newcommand{\bal}{^{{\rm bal}}}
\newcommand{\lwq}{{_{\bQ+}}}
\newcommand{\wtau}{\tilde\tau}
\newcommand{\wTheta}{\widetilde\Theta}

\newcommand{\bZge}{\bZ_{\ge0}}
\newcommand{\mom}{{{\rm Mom}}}
\newcommand{\MP}{{{\rm MP}}}
\newcommand{\Ver}{{{\rm Vert}}}
\newcommand{\Ustd}{\cU_{{\rm std}}}
\newcommand{\Umax}{\cU_{{\rm max}}}
\newcommand{\Umin}{\cU_{{\rm min}}}
\newcommand{\Vmin}{\cV_{{\rm min}}}
\newcommand{\simp}{{{\rm simp}}}
\newcommand{\lub}{{{\rm lub}}}
\newcommand{\Lbullet}{\stackrel{\bullet}{L}}
\newcommand{\cMbullet}{\stackrel{\kern 4pt\bullet}{\cM}}
\newcommand{\cPbullet}{\stackrel{\kern 3pt\bullet}{\cP}}

\newcommand{\Rdag}{R^{\dag}}
\newcommand{\Pdag}{P^{\dag}}
\newcommand{\Ldag}{L^{\dag}}
\newcommand{\thetadag}{\theta^{\dag}}
\newcommand{\Thetadag}{\Theta^{\dag}}
\newcommand{\Deltad}{\Delta^{\dagger}}
\newcommand{\Deltadd}{\Delta^{\ddagger}}
\newcommand{\deltadd}{\delta^{\ddagger}}
\newcommand{\deltad}{\delta^{\dagger}}
\newcommand{\Yd}{Y^{\dagger}}
\newcommand{\Ydd}{Y^{\ddagger}}
\newcommand{\Cd}{C^{\dagger}}
\newcommand{\Cdd}{C^{\ddagger}}
\newcommand{\ydd}{y^{\ddagger}}
\newcommand{\Pd}{P^{\dagger}}
\newcommand{\Pdd}{P^{\ddagger}}
\newcommand{\dPd}{\dP^{\dagger}}
\newcommand{\dPdd}{\dP^{\ddagger}}
\newcommand{\Ld}{L^{\dagger}}
\newcommand{\Ldd}{L^{\ddagger}}
\newcommand{\dLd}{\dL^{\dagger}}
\newcommand{\dLdd}{\dL^{\ddagger}}
\newcommand{\cUdd}{\cU^{\ddagger}}
\newcommand{\hchid}{\hchi^{\dagger}}
\newcommand{\hchidd}{\hchi^{\ddagger}}
\newcommand{\taud}{\tau^{\dagger}}
\newcommand{\tauw}{\tilde{\tau}}
\newcommand{\taudd}{\tau^{\ddagger}}
\newcommand{\phid}{\phi^{\dagger}}
\newcommand{\phidd}{\phi^{\ddagger}}
\newcommand{\varphidd}{\varphi^{\ddagger}}

\newcommand{\ConvHull}{{{\rm ConvHull}}}
\newcommand{\Tr}{{{\rm Tr}}}
\newcommand{\Sec}{\Sigma}
\newcommand{\Stratum}{{{\rm Strat}}}
\newcommand{\Strat}{{{\rm Strat}}}
\newcommand{\bSM}{{{\rm S\bM}}}
\newcommand{\TP}{{{\rm TP}}}
\newcommand{\cTP}{{{\rm \cT\cP}}}
\newcommand{\lge}{_{\ge0}}

\newcommand{\gp}{({\rm gp})}
\newcommand{\cov}{({\rm cov})}
\newcommand{\val}{{{\rm val}}}

\newcommand{\Funge}{\Fun_{ge0}}
\newcommand{\fixed}{{\rm fixed}} 
\newcommand{\framed}{{\rm fr}} 

\newcommand{\MDA}{{{\rm M_{(\Delta,A)}}}}
\newcommand{\MQ}{{{\rm M_Q}}}

\newcommand{\NT}{{{\rm NT}}}

\newcommand{\hubL}{\widehat\ubL}
\newcommand{\uPt}{\widetilde{\underline{P}}}

\newcommand{\cAg}{\cal A_{g}}
\newcommand{\cAgd}{\cal A_{g,d}}
\newcommand{\cAgdH}{\cal A_{g,d}^H}
\newcommand{\Agd}{{{\rm A}_{g,d}}}
\newcommand{\AgdH}{{{\rm A_{g,d}^H}}}

\newcommand{\Int}{{{\rm Int}}}
\newcommand{\SCT}{{{\rm SCT}}}
\newcommand{\SQA}{{\cal S\cal Q\cal A}}
\newcommand{\QA}{{\cal Q\cal A}}
\newcommand{\NR}{N\otimes\bR}
\newcommand{\NOR}{N_0\otimes\bR}
\newcommand{\XR}{X_{\bR}}
\newcommand{\bXR}{\bX_{\bR}}
\newcommand{\avor}{\overline{A}_g^{Vor}}
\newcommand{\Rp}{_1R}
\newcommand{\Pp}{_1P}
\newcommand{\Ru}{R_1}
\newcommand{\Rud}{R_d}
\newcommand{\oXR}{\oX_{\bR}}

\newcommand{\APgd}{{{\rm AP}_{g,d}}}
\newcommand{\APgdH}{{{\rm AP}^H_{g,d}}}
\newcommand{\APg}{{{\rm AP}_{g}}}
\newcommand{\AP}{{{\rm AP}}}
\newcommand{\cAP}{{\cal A \cal P}}
\newcommand{\cAPgd}{{\cal A \cal P}_{g,d}}

\newcommand{\Ag}{{{\rm A_g}}}
\newcommand{\oA}{\overline{{\rm A}}}

\newcommand{\phidel}{\phi_{\delta}}
\newcommand{\LDelR}{L_{\Delta,\bR}}
\newcommand{\Adel}{A_{\delta}}
\newcommand{\Adelp}{A_{\delta'}}
\newcommand{\alpydel}{\alpha_{\oy,\delta}}
\newcommand{\alpy}{\alpha_{\oy}}
\newcommand{\Kydelp}{\wK_{\oy,\delta}}
\newcommand{\Kdelp}{\wK_{\delta}}
\newcommand{\Jdelp}{\wJ_{\delta}}
\newcommand{\wKDel}{\wK_{\Delta}}
\newcommand{\Kydel}{K_{\oy,\delta}}
\newcommand{\KDel}{K_{\Delta}}
\newcommand{\JDel}{J_{\Delta}}
\newcommand{\JDelF}{J^F_{\Delta}}
\newcommand{\JdelF}{J^F_{\delta}}
\newcommand{\KDelp}{\wK_{\Delta}}
\newcommand{\SIQA}{\Sigma(Q,A)}
\newcommand{\MQA}{M_{(Q,A)}}
\newcommand{\MQAp}{\wM_{(Q,A)}}
\newcommand{\MQApred}{\wM_{(Q,A)\red}}
\newcommand{\MQAred}{M_{(Q,A)\red}}
\newcommand{\MDel}{M_{\Delta}}
\newcommand{\MDelp}{\wM_{\Delta}}
\newcommand{\UDel}{U_{\Delta}}
\newcommand{\UDelF}{U_{\Delta}^F}
\newcommand{\RDelF}{R_{\Delta}^F}
\newcommand{\UDelp}{\wU_{\Delta}}
\newcommand{\Kdel}{K_{\delta}}
\newcommand{\HDel}{H_{\Delta}}
\newcommand{\LQA}{L_{A}}
\newcommand{\Psq}{P_{Q}}
\newcommand{\Lsq}{L_{Q}}
\newcommand{\Tsq}{T_{Q}}
\newcommand{\Xsq}{X_{Q}}
\newcommand{\oXsq}{\oX_{Q}}
\newcommand{\dXsq}{\dX_{Q}}
\newcommand{\Doy}{D_{\oy}}
\newcommand{\Dsi}{D_{\upsilon}}

\newcommand{\Del}{{{\rm Del}}}
\newcommand{\SDel}{S_{\Delta}}
\newcommand{\GDel}{G_{\Delta}}
\newcommand{\Gdel}{G_{\delta}}
\newcommand{\hGDel}{\widehat{G}_{\Delta}}
\newcommand{\Xdel}{\dX_{\delta}}
\newcommand{\Rdel}{R_{\delta}}
\newcommand{\Ldel}{L_{\delta}}
\newcommand{\Pdel}{P_{\delta}}
\newcommand{\Thetadel}{\Theta_{\delta}}
\newcommand{\LDel}{L_{\Delta}}
\newcommand{\thetadel}{\theta_{\delta}}

\title{Complete moduli in the presence\\ of semiabelian group action}
\shorttitle{Complete moduli} 

\author{Valery Alexeev} 
 \institutions{University of Georgia, 
 Athens, GA\\
{\eightpoint {\it E-mail address\/}: valery@math.uga.edu}}

 \centerline{\bf Abstract}
\vglue12pt

  I prove the existence, and describe the structure, of moduli space
  of pairs $(P,\Theta)$ consisting of a projective variety $P$ with
  semiabelian group action and an ample Cartier divisor on it
  satisfying a few simple conditions. Every connected component of
  this moduli space is proper. A component containing a projective
  toric variety is described by a configuration of several polytopes,
  the main one of which is the secondary polytope. On the other hand,
  the component containing a principally polarized abelian variety
  provides a moduli compactification of $A_g$. The main irreducible
  component of this compactification is described by an ``infinite
  periodic" analog of the secondary polytope and coincides with the
  toroidal compactification of $A_g$ for the second Voronoi
  decomposition.
\vglue2pt
 \begin{small}

\def\sni#1{\vglue2pt\noindent\hglue12pt{#1}. }
\def\pni#1{\vglue2pt\noindent\hglue12pt{#1}.}
\def\ssni#1{\vglue-1pt\noindent\hskip38pt {#1}.}
\vglue2pt \centerline{\bf Contents} \vglue2pt
\noindent Introduction
\vglue3pt
\noindent Section 1. Overview
\vglue-1pt\noindent \hskip12pt 1.1. \ Fundamental definitions
\ssni{A} Stable toric and semiabelic varieties
\ssni{B} The relative case
\ssni{C} Complexes of cones and lattice  polytopes
\sni{1.2} Main results
\ssni{A} Moduli problem: abelic case
\ssni{B} Classification of stable toric varieties and pairs over closed fields
\ssni{C} Classification of SSAVs and pairs over closed fields
\ssni{D} Classification of SSAVs and pairs over $\Bbb C$
\ssni{E} Cohomologies of polarized stable semiablic varieties
\ssni{F} Types and singularities of SSAVs
\ssni{G} Moduli problem: the finite case
\ssni{H} Moduli problem: the infinite periodic case

\vglue3pt\noindent Section 2. Stable toric varieties and pairs
\sni{2.1} Complexes, posets and (co)sheaves on them
\sni{2.2} Basic (co)sheaves and relations between them
\sni{2.3} Affine stable toric varieties
\sni{2.4} Polarized stable toric varieties with a linearized line bundle
\pagebreak

\smallbreak\noindent \hglue12pt {2.5. } Cohomologies of polarized STVs
\sni{2.6} Stable toric pairs
\sni{2.7} The moment map for stable toric pairs over $\Bbb C$
\sni{2.8} One-parameter families of stable toric pairs
\sni{2.9} The modified complex $\Bbb M_\ast$
\pni{2.10} Moduli of stable toric pairs
\pni{2.11} Approximation of the moduli space
\pni{2.12} Generalized secondary polytopes for finite complexes
\pni{2.13} The moment map for the moduli space
\pni{2.14} Relation to families of Gelfand-Kapranov-Zelevinsky
\pni{2.15} Other families
\pni{2.16} First nontrivial example

\vglue2pt\noindent Section 3. Abelic pairs
\vglue2pt\noindent Section 4. Linearization of torus action
\sni{4.1} The polarization morphism and theorem of square
\ssni{A} The case of a closed field
\ssni{B} The case of families
\sni{4.2} Tori in Picard groups and infinite covers
\sni{4.3} Infinite covers arising from the torus action

\vglue2pt\noindent Section 5. Stable semiabelic varieties and pairs
\sni{5.1} (Co)sheaves on general cell complexes
\sni{5.2} Linearized varieties with a nontrivial abelian part
\sni{5.3} Arbitrary polarized SSAVs
\sni{5.4} Arbitrary stable semiabelic pairs
\sni{5.5} Stable semiabelic pairs over $\Bbb C$
\sni{5.6} Mumford-Faltings-Chai's uniformization of abelian varieties
\sni{5.7} One-parameter families of stable semiabelic pairs
\sni{5.8} The modified complex $\Bbb M_\ast$
\sni{5.9} Test families of stable semiabelic pairs
\ssni{A} The totally degenerate case: $\dim \Delta'=g$
\ssni{B} The general case
\pni{5.10} Moduli of stable semiabelic pairs
\pni{5.11} Approximation of the moduli space
\pni{5.12} Generalized secondary polytopes for periodic decompositions
\pni{5.13} Periodic non-Delaunay decompositions
\vglue2pt\noindent Section 6. Further questions
\vglue3pt\noindent References
 \end{small}

\intro

  In this paper I construct a proper moduli stack and its coarse
  moduli space, also proper, of pairs $(P,\Theta)$, where $P$ is a
  reduced projective scheme with a semiabelian group action and
  $\Theta$ is a Cartier divisor.  The pairs should satisfy a few
  simple conditions.  For me, the motivation for considering moduli of
  such pairs came from the {\it log minimal model program\/}; see
  [Ale2]. However, the minimal model program in dimensions
  $\ge4$ remains conjectural, and  so in this paper I do not use it.
  
  When the group acting on our variety is a torus and the sheaf
  $\cO_P(\Theta)$ is linearized, our theory is closely related to the
  Gelfand-Kapranov-Zelevinsky study of families of cycles and their
  toric degenerations; see [GKZ] and
  references there.  (The cycles appear as images of our pairs under a
  natural projection to $\bP^n$.)  The notions of regular
  decompositions and secondary polytopes defined there play an
  important role in our study. It turns out the closure of the stratum
  containing the standard projective toric variety $P[Q]$ for a
  lattice polytope $Q$ is itself the projective toric variety
  corresponding to the secondary polytope $\Sec (Q)$.  Examination of the
  other strata in our moduli space naturally leads to definition of the
  generalized secondary polytopes for each subdivision of $Q$. I
  describe an explicit stratification on the space and compute the
  dimensions of the strata. Their exact configuration depends on the
  resolution of the so-called {\it generalized Baues problem} in geometric
  combinatorics. (I thank Professor Sturmfels for pointing out this
  connection to me.)
  
  One of the most interesting strata in our moduli space is the one
  corresponding to abelian torsors with a divisor defining principal
  polarization. This stratum can be naturally identified with $A_g$,
  the moduli of principally polarized abelian varieties.  Since our
  moduli is proper, we obtain a geometric way of compactifying $A_g$,
  with every point on the boundary corresponding to a simple
  projective object. The ideas of the toric case turn out to apply
  here. The analog of a regular decomposition is the so-called Delaunay
  decomposition. I also construct the corresponding ``infinite
  $\GL$-invariant'' secondary polytope. The closure of $A_g$ is the
  toroidal compactification for the second Voronoi cone decomposition. As
  a consequence of our construction with the polytope, it is
  projective, which was apparently unknown until now. We note that
  for a higher degree of polarization $d$ there is a stratum in our
  moduli which is fibered over $A_{g,d}$ with   fibers of dimension
  $d-1$. There should be a way of looking at only some of the pairs to
  get a space with a finite morphism to $A_{g,d}$. There should also
  be a way to define appropriate level structures for our pairs. We
  will deal with   these questions elsewhere.

  According to the {\it stable reduction theorem}, when abelian varieties
  degenerate, the limit is a semi-abelian group variety which is an
  extension of an abelian variety by a torus. Therefore, the situation
  near the boundary of abelian moduli and the toric situation have a
  lot in common.  
  The main idea of this work is to understand everything first in the
  ``finite'' case where it is quite simple
  and then, when working in the ``infinite periodic'' case, exploit the
  analogy to the maximum.
  
  In the ``finite'' case, i.e.\ of pairs with a torus action and
  linearized $\cO_P(\Theta)$, the moduli are constructed 
  directly, bypassing any deformation theory. The moduli stack of such 
  toric pairs is glued directly from several smooth groupoids $[U/R]$
  with finite stabilizers. In this case $U$   is the spectrum of an
  explicit semigroup algebra, and $R$ is given by a torus action. Hence, we get an Artin stack and it is coarsely
represented
  by an algebraic space according to [KM]. Moreover, in
  our case this space turns out to be a projective scheme, again
  with an explicit combinatorial description.
  
  In the ``infinite periodic'' case, the one corresponding to abelian
  torsors and their degenerations, the situation is a little more
  complicated because in a semi-abelian family the toric rank may vary.
  However, for as long as we only look at an infinitesimal family over
  $S$, the toric rank stays constant and the semi-abelian group scheme
  $G/S$ is a global extension of an abelian scheme $A/S$ by a torus
  $T/S$. These families we can describe by the same toric methods and
  semigroup algebras as before. Hence, we get all the deformation
  theory essentially ``for free" here. To prove that our moduli stack is
  an algebraic Artin stack with finite stabilizer (again, coarsely
  representable by an algebraic space by [KM]) we use
  Artin's method in [Art2]. Our basic method in this
  case, therefore, is very similar to that of Faltings and Chai
  \cite{FC} but with the following difference. The
  considerations of \cite[IV.3, IV.4]{FC} are very delicate
  because they did not have an easily definable functor and could not
  quote\break M. Artin directly. Our considerations are not so delicate
  because we do and we can. We note that once one toroidal
  compactification is constructed, all others can be constructed from 
  it by blowing up and down on the boundary and, for moduli with level
  structures, also looking at finite covers and normalizing.

  The subject of compactifications of moduli of abelian varieties has
  a very rich history which we will not attempt to summarize here.  
  Some directly relevant papers, aside from those already  mentioned are: 
  [AMMT], [Mum1], [Nam1,2].

  About the base scheme, there was not much to gain by working over a
  (closed) field $k$; much of our work consists of working with
  semigroups and semigroup algebras, and a semigroup algebra $\bZ[H]$
  is just as simple as semigroup algebra $k[H]$. A toric scheme over
  $\bZ$ corresponding to a cone or a polytope is just as simple as the 
  similar toric variety over $k$. Moduli of abelian varieties also
  exist over $\bZ$. Hence, for   most of this paper we work over
  $\bZ$. One can make obvious changes to modify our arguments and
  notation  for the case of a field or some other   scheme.

\demo{Acknowledgements}
  I am indebted to many people with whom I interacted during the work
  on this paper. It started after a suggestion made to me by
  Y.~Kawamata concerning my article on the two-dimensional case.
  Collaboration with I.~Nakamura on \cite{AN} was a
  very important first step.  Discussions with J.~Koll\'ar have helped
  me to crystallize the definition of the functor.  I have also
  profited greatly from discussions with F.~Bogomolov about \hbox{theorem} of
  square and with R.~Erdahl about Delaunay decompositions. I had an
  opportunity to talk about parts of this work as well with R.~Varley,\break
  G.~Faltings, M.~Kapranov, H.~Lenstra, L. Moret-Bailly, S.~Mori, T.
  Oda, F.~Oort,  W. J. Whiteley and my colleagues at the University of
  Georgia. I am very grateful for their suggestions. I also thank
  \myunderline{T.~Kajiwara and} T.~Ueno for pointing out some
  misprints.
  
  I would also like to thank RIMS, Kyoto University, where part of
  this work was done, for hospitality and an excellent research
  environment, as well as NSF and the J. P. Sloan Foundation for providing
  partial support.
\enddemo

\section{Overview}
\label{sec:Overview}

1.1. {\it Fundamental definitions}.
 \vglue2pt

\hglue22pt A. {\it Stable toric and semiabelic varieties}.
  Fix an algebraically closed field $k$ of arbitrary characteristic.
  By a variety over $k$ we will mean a connected reduced but not
  necessarily irreducible separated scheme of finite type over~$k$.
  Recall the following:

\demo{Definition {\rm 1.1.1}} 
 1.  A {\it multiplicative torus\/} is a direct sum of several
    copies of the multiplicative group of the field $k$: $T=\bG_m^r$.
\begin{itemize}  
  \item[2.]  An {\it abelian variety\/} is a (connected) proper group
    variety $A$ over~$k$. It is then necessarily commutative and
    projective.
  \item[3.] A {\it semiabelian variety\/} is a group variety which is an
    extension
   $$ 
      1\to T \to G \to A \to 0
    $$ 
    of an abelian variety $A$ by a torus $T$. It is automatically
    commutative as well. Such extensions are classified by
    homomorphisms $c:X\to A^t$ from the group of characters of $T$ to
    the dual abelian variety of $A$.
  \end{itemize}
\enddemo
 
{\it Definition} 1.1.2.
  Let $L$ be an invertible sheaf on an abelian variety $A$. A~{\it polarization morphism\/} $\lambda(L):A\to A^t=\bPic^0 A$ is
  defined by $a\mapsto T^*_aL\otimes L^{-1}$. It is a homomorphism by
  theorem of square. A {\it polarization\/} on $A$ is a homomorphism
  $\lambda:A\to A^t$ which equals $\lambda(L)$ for some ample $L$. The
  sheaf $L$ with this property is not unique, and two such sheaves
  differ by $E\in \Pic^0 P$. The {\it degree\/} of polarization is defined
  as $d=h^0(L)=\sqrt{\deg\lambda}$.

\demo{Definition {\rm 1.1.3}}  1.
    By a {\it toric variety\/} we will mean a variety $P$ together
    with an action of a torus $T$ such that
    \begin{itemize} \item[] \begin{itemize}
    \item[(i)] $P$ is normal,
    \item[(ii)] there are only finitely many orbits.
    \end{itemize}
  \item[2.] By analogy, a {\it semi-abelic variety\/} is a variety $P$
    together with an action of a semi-abelian variety $G$ satisfying
    the same basic conditions.
    
    We will add the following mild condition on the action, both for
    toric and semiabelic cases:
    \begin{itemize}
    \item[] The stabilizer of the generic point is connected, reduced
      and lies in the toric part $T$ of $G$ (we view it as a group
      scheme over $k$).
    \end{itemize}
  \item[3.] If $G=A$ is abelian, we will call $P$ an {\it abelic
      variety.\/}
  \end{itemize}

\enddemo

{\it Remarks} 1.1.4.
  1) We will contrast toric varieties with {\it torus embeddings\/}
  $P\supset T$.  The difference is that for us a toric variety need
  not have an origin $1\in T\subset P$ fixed, and the torus action
  need not be faithful. However, if $T\acts P$ is a toric variety in
  our definition, $T'$ is the quotient of $T$ by the stabilizer of the
  generic point and $p\in P$ is an arbitrary point in the dense orbit,
  then by sending $1\in T'$ to $p$ we obtain the structure of a torus
  embedding $P\supset T'$.
  
  2) Since we assumed $P$ to be normal, it follows that the condition
  (iii) above holds for {\it every\/} point $p\in P$.
    
  3) In the abelian case there is already a perfectly good name
  ``abelian torsor''. We will use the term ``abelic variety'' in the
  interests of uniformity.  Unlike the toric case, allowing the action
  to be nonfaithful does not add any flexibility here, so we will
  forgo it.
\vglue12pt

  I thank S. Keel for suggesting the word ``semi-abelic''.

\demo{Definition {\rm 1.1.5}}
  A {\it stable semiabelic variety}, or SSAV for short, is a variety
  $P$ together with an action of a semiabelian variety $G$ {\it of
    the same dimension} such that:
  \begin{itemize}
  \item[(i)] $P$ is seminormal.
  \item[(ii)] There are only finitely many orbits.
  \item[(iii)] The stabilizer of every point is connected, reduced and lies
    in the toric part $T$ of $G$.
  \end{itemize}

\enddemo

{\it Definition} 1.1.6.
  A reduced scheme $P$ is called {\it seminormal\/} if every proper
  bijective morphism $f:P'\to P$ with reduced $P'$ inducing
  isomorphisms on the residue fields $\kappa(p')\supset\kappa(p)$ for
  each $p\in P$ is an isomorphism.  A reduced scheme $P$ is called
  {\it weakly normal\/} if every proper bijective morphism $f:P'\to P$
  inducing isomorphisms $\kappa(p')\supset\kappa(p)$ for each {\it generic point\/} $p\in P$ and a purely nonseparable extension for
  every other point is an isomorphism.
  \pagebreak

  Every reduced scheme with a finite normalization (i.e.\ $P$ is of
  finite type over a field or $\bZ$ or a localization of such) has
  unique seminormalization and weak normalization. They coincide if
  $\chr k=0$.  See \cite[I.7.2]{Kol} for an
  overview of their properties.
  \cite{GT} contains a wealth of
  information about seminormal schemes.
 \vglue12pt

  We will see that our varieties are weakly normal as well.

\demo{{R}emark {\rm 1.1.7}}
  Following the tradition that started with stable curves, we call our
  most general varieties stable semiabelic varieties. Note that a
  stable semiabelic variety need not be plainly semiabelic; in other
  words, {\it stable} acts as a widening term here, not a narrowing
  one.
\enddemo

{\it Definition} 1.1.8.
  A {\it polarized\/} stable semiabelic variety is a {\it projective\/} SSAV together with an ample invertible sheaf $L$.
  The {\it degree\/} of polarization is defined as $h^0(L)$.

\demo{Definition {\rm 1.1.9}}
  A {\it stable semiabelic pair\/} $(P,\Theta)$ consists of a
  \myunderline{stable} projective semiabelic variety $P$ and an
  effective ample Cartier divisor $\Theta$ on it which {\it does not
    contain any $G$-orbits.\/} $P$ is polarized by $L=\cO_P(\Theta)$.
  
  {\it Semiabelic, abelic, toric\/} and {\it stable toric\/} pairs are
  defined in the same way.
\enddemo

\hglue22pt  B. {\it The relative case}.  Instead of a closed field, we now work over an arbitrary base scheme
  $S$.  We need to recall the basics about semiabelian group schemes.
  The references for this are \cite[X]{SGA3},
  \cite[Ch.\ 6]{Mum2}, \cite[Ch.\ I]{FC}.

\vglue12pt {\it Definition} 1.1.10.  \hangindent=27pt\hangafter=1
 1.  A {\it split\/} (other names: constant or diagonalizable) {\it torus\/} of rank $r$ over a scheme $S$ is a fibered
product over
    $S$ of several copies of $\bG_{m,S}$. A {\it torus\/} $T$ of rank
    $r$ over $S$ is a smooth group scheme over $S$ which locally in
    fppf topology is a split torus of rank $r$.  In other words, there
    exists a covering family $\{U_i\to S\}$ of flat morphisms of
    finite presentation such that
    $T{\displaystyle\mathop{\times}\limits_{S}} U_i\simeq\bG^r_{m,U_i}$. In this definition,
    {\it fppf\/} can be replaced by {\it \'etale.\/}
    
\hglue32pt  \hangindent=27pt\hangafter=1   Equivalently, $T$ is a smooth separated group scheme of finite
    presentation over $S$ whose every geometric fiber is a torus of
    rank $r$ (\cite[XV.8.18]{SGA3}, \cite[I.2.11]{FC}).
\begin{itemize}
 
\item[2.] An {\it abelian scheme\/} over $S$ is a proper smooth group
    scheme $A$ with connected fibers. It is automatically commutative.

\item[3.] A {\it semiabelian scheme\/} over $S$ is a smooth separated
    group scheme over $S$ such that every geometric fiber $G_{\bar s}$
    is a semiabelian variety over $\kappa(\bar s)$. It is commutative
     by \cite[XVI.1.6]{SGA3}.
    
    Since the formation of maximal tori in algebraic groups is
    compatible with field extensions (\cite[XIV.1.1]{SGA3}), 
    a semiabelian group scheme over the spectrum of a nonclosed field
    $k'$ is an extension over $k'$ of a torus $T$ by an abelian
    variety $A$. For any base $S$, if the dimensions of abelian and
    toric parts of $G/S$ are constant then $G$ is a global extension
    of some abelian scheme $A$ by a torus $T$ over $S$
    (\cite[I.2.11]{FC}).
  \end{itemize}
\phantom{hi}
\vglue-12pt

 {\it Definition} 1.1.11.
  A {\it polarization\/} on an abelian scheme is a homomorphism
  $\lambda:A\to A^t$ to the dual abelian scheme $A^t=\bPic^0_{A/S}$
  whose restriction to every geometric fiber $A_{\bar s}$ is induced
  by an ample sheaf $L_{\bar s}$.
 
\vglue12pt {\it Remark} 1.1.12.
  The sheaf $L$ need not exist globally.
\vglue12pt

{\it Definition} 1.1.13.
  A {\it stable semiabelic scheme\/} over a scheme $S$ is a flat
  separated scheme $P$ of finite presentation over $S$ together with
  an action by a semiabelian group scheme $G/S$ such that every
  geometric fiber is a stable semiabelic variety over $\kappa(\bar
  s)$.
 
\vglue12pt {\it Definition} 1.1.14.
  A {\it polarized stable semiabelic scheme\/} over a scheme $S$ is a
  stable semiabelic scheme $G\acts P$ proper over $S$ and an
  invertible sheaf $L$ on $P$ such that every geometric fiber is a
  polarized stable semiabelic variety over $\kappa(\bar s)$.
\vglue12pt
{\it Definition} 1.1.15.
  A {\it stable semiabelic pair\/} over a scheme $S$ is a stable
  semiabelic scheme $G\acts P$ over $S$ and an effective relative
  Cartier divisor $\Theta$ on $P$ such that every geometric fiber is a
  stable semiabelic pair over $\kappa(\bar s)$.
\vglue12pt

  Once again, there are special cases of these definitions: for {\it abelic{\rm ,} toric{\rm ,} stable toric{\rm ,}
 semiabelic\/} schemes and pairs.

\vglue12pt 
\hglue22pt  C. {\it Complexes of cones and lattice polytopes}.
 \vglue6pt {\it Setup} 1.1.16.
  We fix 
 \vglue4pt
   1. a lattice $X\simeq\bZ^r$, 

  2.  a sublattice $i:Y\subset X$,

3.  a trivial $X$-torsor $\oX$.
  \vglue4pt\noindent 
  The difference between $X$ and $\oX$ is that in the latter we do not fix the origin.
  
  The symbol $\bX$ will denote the lattice $\bZ\oplus X\simeq
  \bZ^{r+1}$.  We will denote by $X_{\bR}$ (resp.\ $\oX_{\bR}$,
  $\bX_{\bR}$) the real vector space $X\otimes\bR$ (resp.\ $\oX\otimes\bR$, $\bX\otimes\bR$).
\vglue6pt

{\it Definition} 1.1.17.
  A {\it complex $\Omega$ of cones referenced by $X$} is a set
  $|\Omega|$ and a collection of distinct subsets $\{\omega\}$ called
  {\it cones\/} with $|\Omega|=\cup\omega$, together with a function
  $\rho:|\Omega|\to X_{\bR}$ which identifies each $\omega$ with a
  finitely generated rational polyhedral cone in $X_{\bR}$ with a
  vertex at the origin. The complex   is supposed to be {\it face-fitting\/}: for any
$\omega_1,\omega_2\in\Omega$ all faces of
  $\omega_i$ are in $\Omega$ and $\omega_1\cap\omega_2$ is a union of
  faces in~$\Omega$. The set $|\Omega|$ comes with a (weak) classical
  topology: a subset is open if its intersections with all $\omega$'s
  are open.
  
  The empty set is not in $\Omega$ by definition. The set $\Omega$ has
  a natural partial order -- {\it by the reverse\/} of inclusion of cones
  (so, if $\omega_1$ is a face of $\omega_2$, we say that $\omega_1$
  is greater) and we will use this partial order.  We will use the
  abbreviation {\it poset\/} for a {\it partially ordered set.\/}
  
  Complexes of cones referenced by $\bX$ are defined analogously.

\vglue12pt {\it Definition} 1.1.18.
  A {\it complex $\Delta$ of lattice polytopes referenced by $X$} is a
  set $|\Delta|$ and a collection of distinct subsets $\delta$ called
  {\it cells\/} with $|\Delta|=\cup\delta$, together with a function
  $\rho:|\Delta|\to X_{\bR}$ which identifies each $\delta$ with a
  polytope in $X_{\bR}$ with   vertices in $X$.  The complex is
  supposed to be {\it face fitting\/}: for any
  $\delta_1,\delta_2\in\Delta$ all faces of $\delta_i$ are in $\Delta$
  and $\delta_1\cap\delta_2$ is a union of polytopes in~$\Delta$.
  Once again, we will use the partial order on $\Delta$ {\it by
    reverse of inclusion.}  We will include $\emptyset$ in the poset
  of polytopes by definition.

  A complex of lattice polytopes referenced by $\oX$ is defined
  in the same way, replacing $X_{\bR}$ by $\oX_{\bR}$ and $X$ by
  $\oX$.
\vglue12pt

{\it Assumption} {\rm 1.1.19}.
  We will generally assume that $|\Delta|$ is connected and locally
  finite.

\vglue12pt {\it {R}emark} {\rm 1.1.20}.
  Note that the function $\rho$ need not be injective. The images of
  different cones or polytopes are allowed to intersect  
  arbitrarily.
\vglue12pt

{\it Definition} 1.1.21.
  Let $\wDelta$ be a not necessarily connected and necessarily not
  finite complex of polytopes referenced by $X$ with a fixed action
  $X\acts\wDelta$ compatible with the translation action of $X$ on
  itself. Then, for each cell $\wdelta$ and each $x\in X$ we have
  another cell $T_x\wdelta$ and $T_x \rho(\wdelta) = \rho
  (T_x\wdelta)$.
  \vglue1pt

  In this situation we will define a complex $\Delta=\wDelta/X$. Every
  cell $\delta$ in this complex is a polytope $\wdelta$ with some of
  its faces glued according to the equivalence relation given by
  $X$-translations. Obviously, $X$ acts on the set $|\wDelta|$
  properly discontinuously in the classical topology, and we have a
  set $|\wDelta|/X=|\Delta|=\cup\delta$.
  
  Assuming as before that $|\Delta|$ is connected, let $\Deltad$ be
  one of the connected components of $\wDelta$ and denote by $Y$ the
  subgroup of $X$ leaving $\Deltad$ invariant.  Then
  $\Delta=\Deltad/Y$. Note that the complex $\Deltad$ is referenced by
  $X$ and that it usually cannot be chosen canonically: the set of
  connected components of $\wDelta$ is in a 1-to-1 correspondence with
  $X/Y$.
  
  We will call the complex $\Delta$ defined in the above described way
  a {\it complex of lattice polytopes referenced by $\oX/Y$.\/} Note
  that the map $\rho:|\wDelta|\to X_{\bR}$ induces a map
  $\rho:|\Delta|\to \oX_{\bR}/Y$.
\pagebreak

  To deal with stable semiabelic varieties with nontrivial abelian
  part in a uniform way we will need   yet another generalization of
  the above concepts that allows cells to be noncompact.

\demo{Definition {\rm 1.1.22}}
  Assume:
  \begin{itemize}
  \item[1.] a surjective homomorphism $X\to X'$ to another lattice $X'$
    with kernel~$X_1$;
  \item[2.] a finite and connected complex of lattice polytopes $\Delta'=
    \Deltaw'/X'= \Deltad/Y'$ referenced by $\oXR'/Y'$.
  \end{itemize}
  Consider a complex $\wDelta$ which is the union of several copies of
  the pullback of the complex $\Deltaw'$ under $X\to X'$.  If $X_1$ is
  infinite then cells $\deltaw$, pullbacks of cells $\deltaw'$, are
  not compact anymore. Assume in addition that we are given an action
  of $X$ on $\Deltaw$ compatible with translation and extending the
  $X'$-action on $\Deltaw'$ such that the cell complex $\Deltaw/X$ is
  connected.  Then the stabilizer group $Y_1$ of every cell $\deltaw$
  is a subgroup of finite index in $X_1$.  This complex
  $\Delta=\Deltaw/X$ is our most general cell complex. It has compact
  cells which are quotients by $Y_1$ of pullbacks of polytopes. There
  exists a split exact sequence $0\to Y_1\to Y\to Y'\to 0$ with a
  sublattice $Y\subset X$, and $\Delta=\Deltad/Y$ is referenced by
  $\oX/Y$. Here, $\Deltad$ is one of the copies of the pullback of
  $\Delta'{}^{\dag}$.
\enddemo

{\it Example} 1.1.23.
  If $X'=0$ then $\Delta$ consists of one big cell $\oX_{\bR}/Y\simeq
  \bR^r/\bZ^r$.

\demo{{R}emark {\rm 1.1.24}}
  A cell complex referenced by $\oX/Y$ is a category which is {\it essentially a poset\/.} This means that for every
  $\delta_1\ge\delta_2$ there may be several arrows
  $\delta_1\to\delta_2$ but only the  arrows $\delta\to\delta$ are
  identities.
\enddemo

{\it Definition} 1.1.25.
  For a connected complex $\Deltad$ covering $\Delta$ there is a
  unique connected complex $\Deltadd$ covering $\Deltad$ such that
  $|\Deltadd|$ is {\it simply connected.\/} 
  Note that it is referenced by $X$ again and that $|\Delta|$ is the
  quotient of $|\Deltadd|$ by a free and discontinuous action of
  $\pi_1(\Delta)$.

\demo{Definition {\rm 1.1.26}}
  We   are in the
  \begin{itemize}
  \item[1.] {\it finite case\/} if $Y=0$,
  \item[2.] {\it infinite periodic case\/} if $Y$ has finite index in $X$,
  \item[3.] {\it intermediate case\/} otherwise.
  \end{itemize}
\enddemo

{\it Definition} 1.1.27.
  A {\it pointed cell complex $(\Delta,C)$\/} is defined by adding a
  set $C\subset|\Delta|\cap\rho^{-1}(\oX/Y)$ and subsets
  $C_{\delta}=C\cap\delta$ for each cell $\delta$ which {\it must
    contain all the vertices} of $\delta$. In the case of the
  pull-back from a lower-dimensional polytopal complex $\wDelta'$,
  the images of $C_{\delta}$ must contain the vertices of $\wdelta'$.

\demo{Definition {\rm 1.1.28 (order on complexes)}}
  Let $\cD_1,\cD_2$ be two cone or cell decompositions, or pointed
  cell decompositions referenced by $\XR$.  We say that $\cD_1$ is a
  {\it sub\/{\rm -}\/decomposition\/} of $\cD_2$ and denote it $\cD_1\le\cD_2$
  if the cones, resp.\ cells, resp.\ pointed cells of $\cD_1$ form a
  subdivision of those of $\cD_2$.  For pointed cells, $C_1$ has to be
  a subset of $C_2$.
\enddemo

{\it Definition} 1.1.29.
  We will call a complex $\cD$ {\it $1$-sheeted convex} if $\rho$ is
  injective and $\rho(|\cD|)$ is convex. Explicitly, this means that
  \begin{itemize}
  \item[1.] in the case of cones, $\Omega$ is a subdivision of a single
    finitely generated rational cone;
  \item[2.] in the case of cells, $|\Deltad|$ is a pullback of a single
    polytope of lower dimension.
  \end{itemize}

  The following definition is consistent with the one for simplicial
  complexes appearing in connection with Stanley-Reisner rings (see
   \cite{Sta1},
\cite{BH}), as well as with the more
  general definition of \cite{Bac} for arbitrary
  posets. In \cite{Yuz} CM posets were called
  $k$-spherical.

\demo{Definition {\rm 1.1.30}}\label{defn:CM_complex}
  A complex $\cD$ is called {\it almost Cohen-Macaulay\/} with respect to $k$,
  abbreviated ACM, or {\it locally Cohen-Macaulay\/} if it is locally
  finite and for every point $u\in|\cD|$ one has
  $\wH_i(|\cD|,|\cD|-u,k)=0$ for $0\le i<\dim\cD$.  Here $\wH_i$
  denotes the reduced singular homology. It is called {\it Cohen-Macaulay\/} with respect to $k$ if, in addition
$\wH_i(|\cD|,k)=0$
  for $0\le i< \dim\cD$.
\enddemo

{\it Example} 1.1.31.
  One has the following implications: $\cD$ is 1-sheeted convex
  $\follows$ $|\cD|$ is a topological manifold (possibly with
  boundary) $\follows$ $\cD$ is locally CM.

\demo{{\rm 1.2.} Main results}
\vglue4pt \hglue22pt A. {\it Moduli problem\/{\rm :} abelic case}.
 
\proclaimtitle{see \ref{lem:polarization_ca_to_a}}
\specialnumber{1.2.1}\proclaim{Lemma}
  Every polarized abelic scheme $(A\acts P,L)/S$ of degree $d$ defines
  in a natural way a polarization $\lambda:A\to A^t$ of degree
  $d=\sqrt{\deg\lambda}$ on~$A/S${\rm .}
\endproclaim

\proclaimtitle{\ref{thm:APgd}}
\specialnumber{1.2.2}\proclaim{Theorem}
  \hiha
{\rm 1.} The moduli stack $\cAPgd$ of abelic pairs $(P,\Theta)$ of
    degree $d$ is a separated Artin stack with finite stabilizers and it
    comes with a natural map of relative dimension $d-1$ to the stack
    ${\cal A}_{g,d}$ of polarized abelian varieties{\rm .}
\begin{itemize}  \ritem{2.} $\cAPgd$ has a coarse moduli space $\APgd$ which is a
    separated scheme and   comes with a natural projective map of
    relative dimension $d-1$ to~$\Agd$.
  \ritem{3.} Over $\bZ[1/d]$ the schemes $\APgd$ and $\Agd$ are disjoint
    unions of components naturally labeled by $g$\/{\rm -}\/dimensional $1$\/{\rm -}\/cell
    complexes $\oXR/Y$ with $|X/Y|=d${\rm .}
  \end{itemize}
\endproclaim

\proclaimtitle{\ref{cor:Ag=APg}, \ref{cor:APg_exists}}
\specialnumber{1.2.3}\proclaim{Theorem}
  In the principally polarized case{\rm ,} i.e.\ when $d=1$, $\cAP_g=\cA_g$
  and $\APg=\Ag${\rm .}
\endproclaim

\hglue22pt B. {\it Classification of stable toric varieties and pairs over
  closed fields}.
We start with the toric case because it is easier to formulate.
 
\demo{Setup {\rm 1.2.4}}
  We fix a torus $T$ of dimension $r$.  Varieties in this subsection
  all have $T$-action. We denote by $X\simeq\bZ^r$ the Cartier dual of
  $T$, so that $T=\bHom(X,\bG_m)$.
  
  Symbols $\uT$, $\ubT$, etc. denote certain constructible sheaves of
  commutative groups on the topological space $|\cD|$ to be
  defined in the main text. The reader can guess their meanings from
  the notation. The complex $\hbM^*$ is defined as the cylinder of the
  morphism of cochain complexes $C^*(\ubT)\to C^*(\uhFun_C)$.
\enddemo

\proclaimtitle{\ref{thm:every_affine_STV_is_such},
  \ref{thm:affine_STVs_isoclasses}} 
\specialnumber{1.2.5}\proclaim{Theorem}
  There is a $1$\/{\rm -}\/to\/{\rm -}\/$1$ correspondence between
  \begin{itemize}
  \ritem{1.} affine stable toric varieties $P${\rm ,}
  \ritem{2.} the following data\/{\rm :}\/
    \begin{itemize}
    \ritem{(a)} a finite complex $\Omega$ of cones referenced by $X$\/{\rm ,}\/
    \ritem{(b)} an element of the cohomology group
      $H^1(\Omega,\uT)/\Sym\, \Omega${\rm ,} where $\Sym\, \Omega$ is the {\rm (}\/finite\/{\rm )}
      group of symmetries of $\Omega$ compatible with automorphisms of
      $X$\/{\rm .} 
    \end{itemize}
  \end{itemize}
  Moreover{\rm ,} $\Aut P$ is a finite extension of $H^0(\Omega,\uT)${\rm .}
\endproclaim

\proclaimtitle{\ref{thm:polarized_STVs},
  \ref{thm:polarized_linearized_STVs},
  \ref{thm:classif_arbitrary_polarized_SSAV}}
\specialnumber{1.2.6}\proclaim{Theorem}
  There is a $1$\/{\rm -}\/to\/{\rm -}\/$1$ correspondence between
  \begin{itemize}
  \ritem{1.} polarized stable toric varieties $(P,L)${\rm ,}
  \ritem{2.} the following data\/{\rm :}
    \begin{itemize}
    \ritem{(a)} a sublattice $Y\subset X$ and an $X$-torsor $\oX${\rm ,}
    \ritem{(b)} a finite complex $\Delta$ of lattice polytopes referenced by
      $\oX/Y${\rm ,}
    \ritem{(c)} an element of $H^1(\Delta,\ubT)/\Sym\, \Delta${\rm .}
    \end{itemize}
  \end{itemize}
  Moreover{\rm ,} $\Aut (P,L)$ is a finite extension of $H^0(\Omega,\ubT)${\rm .}
\endproclaim

\proclaimtitle{\ref{thm:linearized_STpairs},
  \ref{thm:arbitrary_polarized_SSApairs}}
\specialnumber{1.2.7}\proclaim{Theorem}
  There is a $1$\/{\rm -}\/to\/{\rm -}\/$1$ correspondence between
  \begin{itemize}
  \ritem{1.} stable toric pairs $(P,\Theta)${\rm ,}
  \ritem{2.} the following data\/{\rm :}\/
    \begin{itemize}
    \ritem{(a)} a sublattice $Y\subset X$ and an $X$-torsor $\oX${\rm ,}
    \ritem{(b)} a finite complex $(\Delta,C)$ of pointed lattice polytopes
      referenced by $\oX/Y${\rm ,}
    \ritem{(c)} an element of the cohomology group $H^1(\hbM^*)/\Sym\, \Delta${\rm .}
    \end{itemize}
  \end{itemize}
  Moreover{\rm ,} $\Aut (P,\Theta)$ is a finite extension of $H^0(\hbM^*)${\rm .}
\endproclaim

\hglue22pt C. {\it Classification of {\rm SSAVs} and pairs over
  closed fields}. 
 
\demo{Setup {\rm 1.2.8}}
  We fix a semiabelian variety $G$ with the abelian part $A$ of
  dimension $a$ and the toric part $T$ of dimension $r$, $a+r=g$.
  Varieties in this subsection all have a $G$-action. As before, we
  denote by $X'\simeq\bZ^r$ the Cartier dual of $T$.  The extension $G$
  is defined by a homomorphism $c:X'\to A^t$.  
  
  We also fix a cell complex referenced by $\oX'/Y'$. We will denote
  by $i:\pi_1(|\Delta'|)\to X'$ the natural projection.

  Additionally, one has to have
  \vglue4pt
\hglue-6pt  1. a polarization $\lambda:A\to A^t$,
 \vglue6pt
\hglue-6pt 2.  a homomorphism $c^t:\pi_1(|\Delta'|)\to A$ such that
    $\lambda \circ c^t =c\circ i:\pi_1(|\Delta'|)\to A^t$. 
 \enddemo

  The first step is to classify the polarized varieties $(P,L)$ and
  pairs $(P,\Theta)$ with   additional data:
  \begin{itemize}
  \item[1.] a fixed point in one of the minimal-dimensional orbits on $P$, 
    which is isomorphic to $A$,
  \item[2.] an ample line bundle $\cM$ on $A$ such that
    $\lambda(\cM)=\lambda$.
  \end{itemize}
  With  these additional data, the theorems of the previous section
  remain   literally true if one understands $H^p$ as functions on
  the homology groups $H_p$ with the values in certain
  $\bG_m(k)$-torsors instead of $\bG_m(k)$; see
  \ref{thm:classif_arbitrary_polarized_SSAV},
  \ref{thm:arbitrary_polarized_SSApairs}.

  In general, the isomorphisms are obtained by dividing the family by
  the equivalence relation corresponding to forgetting the choice of
  the point in a minimal orbit.
  
  Hence, every polarized stable semiabelic variety (resp.\ pair) is
  labeled by a complex $\Delta'$ (resp.\ pointed complex) of lattice
  polytopes for the toric part. In addition, if $\chr k$ does not
  divide $d=\deg\lambda(\cM)$ then one naturally associates to $P$ or
  $(P,\Theta)$ a complex $\Delta$ of dimension $\dim \Delta'+a$ 
  with nonpolytopal cells if $a>0$.

\vglue6pt \hglue22pt D. {\it Classification of {\rm SSAVs} and pairs over $\bC$}.
Over $\bC$, a semiabelian variety $G$ can be written as $(\bC^*)^g/Y
=T/Y$ and a variety $P$ with $G$-action in this way can be considered
to be ``toric''. A different choice of such a representation
corresponds to a different choice of an isotropic sublattice in an
appropriate symplectic lattice.

\proclaimtitle{5.5.6}
\specialnumber{1.2.9}\proclaim{Theorem}
  Polarized stable semiabelic varieties {\rm (}\/resp.\ pairs\/{\rm )} with a fixed
  $T$\/{\rm -}\/action are parametrized by a cell complex $\Delta$ and an open
  subset in $H^1(\Delta,\uhbX)$ {\rm (}\/resp.\ $H^1((\Delta,C), \hbM^*)${\rm ),}
  consisting of classes satisfying the Riemann positivity condition{\rm ,}
  modulo the group of symmetries of $\Delta${\rm .}
\endproclaim

\proclaimtitle{2.7.2, 5.5.7}
\specialnumber{1.2.10}\proclaim{Theorem}
  A stable semiabelic pair with a fixed\break $T$\/{\rm -}\/action naturally defines a 
  moment map $\mom: P(\bC)\to |\Delta|${\rm .} A fiber over a point in the
  interior of a $k$\/{\rm -}\/dimensional cell is isomorphic to
  $(S^1)^k=U(1)^k${\rm .} 
\endproclaim

\hglue22pt E. {\it Cohomologies of polarized stable semiabelic varieties}.
 
\proclaimtitle{2.5.1, 5.4.1}
\specialnumber{1.2.11 }\proclaim{Theorem}
  For all $p,d>0${\rm ,}  $H^p(P,L^d)=0${\rm .}
\endproclaim

\hglue22pt F. {\it Types and singularities of} SSAVs.
We will call the complex $\Omega$, resp.\ $\Delta$, associated to an
  affine, resp.\ polarized, SSAV $P$ by the classification theorems
  above the {\it type\/} of $P$.

\demo{Example {\rm 1.2.12}}
  1) {\it Finite case.\/} The numerical type of a polarized toric
  variety is a lattice polytope $Q\subset \oX_{\bR}\simeq\bR^r$.
  
  2) {\it Infinite periodic case.\/} Let $(A\acts P,L)$ be a polarized
  abelic variety of degree $d$, $\chr k\nmid d$. Let $\lambda:A\to
  A^t$ be the induced polarization on $A$ with $K(\lambda)\simeq
  H\times \hH$, and say $H=X/Y$. In this case the numerical type of
  $(P,L)$ is $\oX_{\bR}/Y\simeq \bR^g/\bZ^g$.
\enddemo

\proclaimtitle{2.10.1, 5.7.4}
\specialnumber{1.2.13}\proclaim{Theorem}
  Let $(G\acts P,L)$ be a polarized stable semiabelic scheme over a
  connected scheme $S${\rm .} If for two geometric fibers the types are
  defined {\rm (}\/i.e.\ $\chr k$ does not divide the degree $d$ of
  polarization on the abelian part\/{\rm )} then $|\Delta_1|$ and $|\Delta_2|$
  can be noncanonically identified{\rm .}
\endproclaim

  In this paper we will only consider the moduli problem for stable
  semiabelic pairs that have the types 1) or 2) of the above example
  or their subdivisions.

\proclaimtitle{2.3.19, 2.4.12, 5.3.9}
\specialnumber{1.2.14}\proclaim{Theorem} \hglue-7pt
  In the affine case{\rm ,} $\Omega$ is {\rm CM} $\follows P$ is {\rm CM.}
  In the polarized case{\rm ,} $\Delta$ is locally {\rm CM} $\follows$ $P$ is {\rm CM.}
\endproclaim

\hglue22pt G. {\it Moduli problem{\rm :} the finite case}.  In the finite case we can make the situation more restricted by
  fixing $T,X$ 
  once and for all and requiring that all automorphisms of stable
  toric pairs preserve them. We then have:

\proclaimtitle{2.10.10, 2.11.11, 2.12.3, 2.12.13}
\specialnumber{1.2.15}\proclaim{Theorem}
  For each lattice polytope $Q${\rm ,}   over $\bZ$ or a field\/{\rm :}
  \begin{itemize}
  \ritem{1.} The moduli stack of stable toric pairs $(P,\Theta)$ of
    type~$\le Q$ is a proper Artin stack with finite stabilizers{\rm .}
  \ritem{2.} It has a coarse moduli space $\MQ${\rm .}
  \ritem{3.} $\MQ$ is projective{\rm .}
  \ritem{4.} $\MQ$ is naturally stratified{\rm ,} and every stratum corresponds
    in a $1$\/{\rm -}\/to\/{\rm -}\/$1$ way to a cell decomposition of $Q$ {\rm (}\/resp.\ a pointed
    cell decomposition of $(Q,C)${\rm ,} where $C=Q\cap X${\rm ).}
  \ritem{5.} The normalization of %
    every irreducible component of $\MQ$ is a projective toric scheme
    corresponding to a generalized secondary polytope{\rm .}
  \ritem{6.} The toric scheme $M_{\Sec(Q)}$ for the secondary polytope of
    $Q$ is isomorphic to the main irreducible component of
    $(M_Q)_{\red}${\rm .}
  \end{itemize}

\endproclaim

\hglue22pt H. {\it Moduli problem{\rm :} the infinite periodic case}.
 
\proclaimtitle{5.10.1}
\specialnumber{1.2.16}\proclaim{Theorem}
  \hiha
 {\rm 1.} The component $\overline{\cAP}_{g,d}$ of the moduli stack of
    semiabelic pairs containing $\cAPgd$ and pairs of the same
    numerical type is a proper Artin stack with finite stabilizers\/{\rm ;}\/
\begin{itemize}   \ritem{2.} It has a coarse moduli space $\overline{\AP}_{g,d}$ as a
    proper algebraic space{\rm ;}
  \ritem{3.} The space $\overline{\AP}_{g,d}$ is naturally stratified
    according to the toric part of $(P,\Theta)${\rm ,} and every stratum
    corresponds  $1$\/{\rm -}\/to\/{\rm -}\/$1$ to a cell decomposition of $\oXR'/Y'$
    with $\dim X'\le g$ and $|X'/Y'|$ dividing $d${\rm ,} modulo symmetries.

  \ritem{4.} Over $\bZ[1/d]$ the space $\overline{\AP}_{g,d}$ has a nicer
    stratification{\rm ,} with every stratum corresponding   $1$\/{\rm -}\/to\/{\rm -}\/$1$  
    to a cell decomposition of $\oXR/Y$ with $\dim X=g$ and $|X/Y|=d${\rm ,} modulo symmetries{\rm .}
  \end{itemize}  

\endproclaim

\proclaimtitle{5.11.6}
\specialnumber{1.2.17}\proclaim{Theorem}  
  In the principally polarized case {\rm (}\/when\break   $Y=X${\rm )} the
  toroidal compactification of $A_g$ for the second Voronoi decomposition
  is isomorphic to the main irreducible component of
  $\overline{\AP}_g${\rm ,} the one containing
  ${\rm A}_{g}=\AP_{g}${\rm .}
\endproclaim

\proclaimtitle{5.12.8}
\specialnumber{1.2.18}\proclaim{Theorem}
  The toroidal compactification of ${\rm A}_g$ for the second
  Voronoi decomposition is projective{\rm .}
\endproclaim

 \vglue-8pt
\section{Stable toric varieties and pairs}
\label{sec:Stable toric varieties and pairs}
\vglue-6pt

\demo{{\rm 2.1.} Complexes{\rm ,} posets and {\rm (}\/co\/{\rm )}\/sheaves on them}
Let $I$ be a {\it poset\/}, i.e.\ a partially ordered set. We can
  look at it as being a category associating an arrow $i_1\to i_2$ to
  each pair $i_1\ge i_2$. (Alternatively, we can associate to such a
  pair an arrow $i_2\to i_1$. The choice is entirely a matter of
  taste. Ours will be more convenient for working with cones and
  polytopes.)
 \enddemo

{\it Definition} 2.1.1.
  The {\it order topology} on $I$, which we will denote $I_{\ord}$, is
  the topology in which the open sets are the {\it increasing sets\/}:
  ${i}\in U$ and ${j}\ge{i} \follows {j}\in U$. The complements of the
  increasing sets are the decreasing sets, so these are closed.  The
  open cover $\Ustd$ by the sets $U_{\ge i}=\{j\ge i\}$ refines any
  other open cover.
\vglue12pt

  Although the major part of the results below requires very little of
  the poset $I$, all our posets will satisfy the following:

\demo{Assumption {\rm 2.1.2}} 
  There is a global bound on the length of all 
  chains in~$I$.
\enddemo

{\it Assumption} 2.1.3.
  $I$ is locally finite; i.e., for each $i$ there are only finitely many
  $j$ comparable to it. 
\vglue12pt

  Under Assumption 2.1.2, $I$ is covered by
  the sets $U_{\ge i_0}$ for all {\it minimal\/} $i_0$.  We will
  denote this locally finite open cover $\{U_{\ge i_0} \,|\, i_0\in
  I_{\min} \}$ by~$\Umin$.

\demo{Definition {\rm 2.1.4}}
  A {\it sheaf $F$ of sets} on $I_{\ord}$ is the same as an inverse
  system of sets $F_i$, called {\it stalks}, indexed by $I$.  Hence,
  for every pair $i_1\ge i_2$ we are given a map
  $\varphi_{i_2i_1}:F_{i_2}\to F_{i_1}$, $\varphi_{ii}=\id$ and
  $\varphi_{i_2i_1} \circ \varphi_{i_3i_2}= \varphi_{i_3i_1}$. For any
  subset $J\subset I$, one has the set of sections$$  
    F(J)={\displaystyle\lim_{{\llar}}}_{\raise2pt\hbox{$\scriptstyle J$}} \,F_j = \left\{(f_j)\in \prod_J
    F_j \,|\, \varphi_{j_2j_1}(f_{j_2})=f_{j_1} \right\} .
  $$  We
  have obvious restriction maps $F(J_2)\to F(J_1)$ for all $J_1\subset
  J_2$, and a presheaf defined this way is automatically a sheaf.
  
  Dually, a {\it cosheaf $G$ of sets} on $I_{\ord}$ is the same as a
  direct system $\{G_i,\psi_{i_1i_2}:G_{i_1}\to G_{i_2}\}$ of sets.
  For any subset $J\subset I$,
$$ 
    G(J)={\displaystyle\lim_{\lrar}}_{\raise2pt\hbox{$\scriptstyle J$}} \,G_j = \coprod_J G_j / \left( g_{j_1} \sim
      \psi_{j_1j_2}(g_{j_1}) \right) .
  $$  
  If the  $F_i$'s have an extra structure, i.e.\ they are groups, semigroups
  or algebras, and $\varphi$'s preserve it, the inverse limit has the
  same structure, so one has sheaves of groups etc. The situation with
  the cosheaves is a bit more complicated.  The direct limits in the
  categories of commutative groups and semigroups exist but as sets
  they do not coincide with the ones above.  Instead, they are given
  by the formula$$  
    G(J)={\displaystyle\lim_{\lrar}}_{\raise2pt\hbox{$\scriptstyle J$}} \,G_j = \mathbold{\oplus}_J G_j / 
    \left< g_{j_1} \sim \psi_{j_1j_2}(g_{j_1}) \right>.
  $$  So, we have to be a little careful about the
  category we work in. This will always be clear from the context.
\enddemo

  When working with the (co)sheaves of commutative groups -- and this
  will be the default case -- we can apply the standard machinery of
  sheaf cohomologies \cite{God} and cosheaf
  homologies (see  \cite{Bre}).  By the very
  definition, we have$$  H^0(I_{\ord},F)=
  {\displaystyle\lim_{{\llar}}}_{\raise2pt\hbox{$\scriptstyle I$}} \,F_i, \qquad H_0(I_{\ord},G)= 
{\displaystyle\lim_{\lrar}}_{\raise2pt\hbox{$\scriptstyle I$}} \,G_i.
  $$  
  Since the open cover $\Ustd$ refines any other, the Cech
  cohomologies \linebreak \u{\it H}\hskip.5pt$^q(I_{\ord},F)$ are isomorphic to 
  \u{\it H}\hskip.5pt$^q(\Ustd,F)$. The natural homomorphisms\break \u{\it H}\hskip.5pt$^q(I_{\ord},F)\to
  H^q(I_{\ord},F)$ are isomorphisms for $q=0,1$ and are injective for
  $q\ge2$ (since $I_{\ord}$ is not paracompact, they are not
  guaranteed to be isomorphisms in general).  Dually, the natural
  homomorphisms $H_q(I_{\ord},F)\to \hbox{{\rm \u{\it H}}}_q(I_{\ord},F)$ are
  isomorphisms for $q=0,1$ and are
  surjective for $q\ge2$.

  The following case will be of special interest:

\demo{Definition {\rm 2.1.5}}
  The poset $I$ is called a {\it join prelattice} or a {\it join
    semilattice} if any two elements ${i_1},{i_2}$ 
  have the least upper bound $\lub({i_1},{i_2})$.
  In this case for all ${i_1},{i_2}\in I$ the set $U_{\ge i_1}\cap
  U_{\ge i_2}$
  equals $U_{\ge \lub({i_1},{i_2})}$.
\enddemo

\specialnumber{2.1.6} \proclaim{Lemma}\label{lem:computing_cohs_for_special_sheaves}
  \hiha {\rm 1.} Any sheaf $F$ {\rm (}\/resp.\ cosheaf\/{\rm )} on the open set $U_{\ge i}$ is
    flabby and hence acyclic{\rm .} \begin{itemize}
  \ritem{2.} If $I$ is a join prelattice then $H^q(I_{\ord}, F)=
    H^q(\Umin,F)$ {\rm (}\/resp.\ for a cosheaf\/{\rm ).}
  \end{itemize}

\endproclaim

\demo{Proof}   
  We will prove these statements for a sheaf. For a cosheaf use a dual
  argument. 1) is obvious since any open subset $U$ of $U_{\ge i}$
  contains $i$ and so $\Gamma(U,F)=F_i$.  Part 2) follows because $F$
  is acyclic when restricted to any of $\mathbold{\cap}\, U_{\ge i}$  (see  \cite[5.9.2]{God}).
\enddemo

  The order topology is quite convenient for computations because
  there are just so few open sets. However, it is not very intuitive.
  Another drawback is that a free action by a group such as $\bZ/n\bZ$
  or $\bZ^r$ on $I_{\ord}$ is often not discontinuous.  There is a
  standard construction associating to each poset $I$ (and, in fact,
  to any category) an ordinary simplicial complex $BI$ which
  comes with the (weak) classical topology: a subset is open if and only if its
  intersection with each simplex is open.

\demo{Definition {\rm 2.1.7}}
  $BI$ is a simplicial complex which has a simplex $\sigma_{i_0\dots
    i_q}$ for each strictly increasing chain $i_0<\dots <i_q$ in $I${\rm ,}
  with obvious inclusions{\rm .}
\enddemo

  Clearly, $BI$ is a disjoint union of the interiors
  $\sigma^0_{i_0\dots i_q}$ of the simplices.

\demo{Definition {\rm 2.1.8}}
  Define two maps 
  \begin{eqnarray*}
    &\phi_l:BI \to I, \qquad \sigma^0_{i_0\dots i_q}\mapsto i_0 
    \qquad & \hbox{mapping to the {\it least} } i,
    \\
    &\phi_g:BI \to I, \qquad \sigma^0_{i_0\dots i_q}\mapsto i_q 
    \qquad & \hbox{mapping to the {\it greatest} } i.
  \end{eqnarray*} \pagebreak

\noindent Denote 
  \begin{eqnarray*}
    U_i=\phi_g^{-1}(U_{\ge i}),&& \quad W_i=\phi_g^{-1}(i)  ,\\
    V_i=\phi_l^{-1}(U_{\ge i}), &&\quad V^0_i=\phi_l^{-1}(i).
  \end{eqnarray*}
\enddemo
 
\specialnumber{2.1.9}\proclaim{Lemma} \label{lem:basics_for_UVWs}
  Let $U$ be an open subset of $I_{\ord}${\rm .} Then
  \begin{itemize}
  \ritem{1.} $\phi_g^{-1}(U)$ is open{\rm ,} i.e{\rm .,} $\phi_g$ is continuous{\rm .} In
    particular{\rm ,} $U_i$ are open and $\oW_i=\phi_g^{-1}(U_{\le i})$ are
    closed{\rm .}
  \ritem{2.} $\phi_l^{-1}(U)$ is closed{\rm ,} in particular $V_i$ are closed{\rm .}
  \ritem{3.} $W_i=U_i\cap \oW_i$ and $V_i^0$ are locally closed and
    contractible to the point $\sigma_i=W_i\cap V_i^0${\rm .}
  \end{itemize}

\endproclaim

\demo{Proof}
This requires only a routine check.
\enddemo

  Since $\phi_g$ is continuous, we can consider the pullback
  $\phi_g^{-1}F$. This is a constructible sheaf on $BI$ whose
  restriction to each $W_i$ is a constant sheaf $F_i$ and the stalk
  at a point $x\in W_i$ is $F_i$. Analogously, $\phi_g^{-1}G$ is a
  constructible cosheaf whose restriction on each $W_i$ is a constant
  cosheaf $G_i$.

\specialnumber{2.1.10}  \proclaim{Lemma}\label{lem:all_cohs_coinside}
  For all $q${\rm ,}
  $H^q(BI,\phi_g^{-1}F)=H^q(I_{\ord},F)$ and \\
  $H_q(BI,\phi_g^{-1}G)=H_q(I_{\ord},G)$.
\endproclaim

\demo{Proof}
  For a sheaf $F$ the statement is an application of the Vietoris
  mapping theorem. Indeed, each fiber $W_i$ is taut in $BI$
  since the latter is a CW complex, and its higher cohomologies with
  the coefficients in $F_i$ vanish because it is contractible. The
  statement for a cosheaf $G$ is obtained by a dual argument.
\enddemo

  From now we will write $H^q(I,F)$ or simply $H^q(F)$ if the poset
  $I$ is understood (resp.\ $H_q(I,G)$, $H_q(G)$).

\demo{Definition {\rm  2.1.11 (Dual sheaf to a cosheaf of semigroups)}}   Let $G$ be a co-
  sheaf of commutative semigroups on $I_{\ord}$, which we now can also
  consider as a cosheaf on its geometric realization $BI$.  The system
  $G_i$ defines a cosheaf of semigroup algebras $\bZ[G_i]$ and,
  dually, an inverse system of affine \myunderline{semi}group schemes
  $\Spec\bZ[G_i]$.  For any scheme $S$, we obtain a sheaf of
  commutative semigroups for the multiplication which we will denote
  by $\hG(S)$:
  \begin{eqnarray*}
    \hG(S)_i &=& \Spec\bZ[G_i](S) =
    \{ \hbox{morphisms } S \to \Spec\bZ[G_i] \} \\
    &=& \{ \hbox{semigroup homomorphisms } G_i \to 
    \Gamma(S,\cO_S) \}.
  \end{eqnarray*} \pagebreak
  If the semigroups $G_i$ are groups then$$ 
    \hG(S)_i =     \{ \hbox{group homomorphisms }
    G_i \to \bG_m(S)=\Gamma(S,\cO_S)^*  \}
  $$ 
  and $\hG(S)$ is a sheaf of commutative groups. The cosheaf $G$ here
  has an ``absolute'' character and the dual cosheaf is relative to
  the choice of a scheme~$S$.
\enddemo

\specialnumber{2.1.12}\proclaim{Lemma}\label{2.1.12}
  \hiha  {\rm 1.} If the  $G_i$ are free abelian groups then there is the following
    noncanonically split exact sequence\/{\rm :}\/
$$ 
      0 \to \Ext\left(\hbox{{\rm \u{\it H}}}_{p-1}(G),\bG_m(S) \right) 
      \to \hbox{{\rm \u{\it H}}\hskip1pt}^p(\hG(S)) \to 
      \Hom\left( \hbox{{\rm \u{\it H}}}_{p}(G) , \bG_m(S) \right) \to 0 .
    $$  \begin{itemize}
  \ritem{2.} Without any assumptions on $G_i${\rm ,} when $S=\Spec\, k$ is the
    spectrum of an algebraically closed field{\rm ,} 
    {\rm \u{\it H}}\hskip1pt$^p(\hG(S))= \Hom(\hbox{{\rm \u{\it H}}}_p(G),\bG_m(k)).$
  \end{itemize}

\endproclaim

\demo{Proof}
  Both parts are immediate applications of the universal coefficient
  theorem (e.g.\ \cite[Thm.5.5.3]{Spa}) and the fact
  that for an algebraically closed field, $\bG_m(k)=k^*$ is a divisible
  group, hence an injective $\bZ$-module.
\enddemo

2.1.13.
  Now consider a complex of lattice polytopes $\Delta=\{\delta_i \,;\,
  i\in I\}$ referenced by a lattice $X$ with the partial order on $I$
  which is the {\it opposite} to the order by inclusion. Sometimes,
  we will write $\Delta$ instead of $I$. Clearly, our assumptions on
  the posets are satisfied.  Pick a point $p_i$ in the interior of
  each polytope $\delta_i$ and identify the simplex $\sigma_{i_0\dots
    i_q}$ with a subset of $\delta_{i_q}$ with vertices $p_{i_0}\dots
  p_{i_q}$.  Then $BI$ is nothing but the barycentric decomposition of~$\Delta$. In particular, $BI=|\Delta|$.  The locally
closed sets
  $V^0_i$ are interiors $\delta_i^0 $ of the polytopes $\delta_i$, and
  the closed sets $V_i$ are the polytopes $\delta_i$ themselves.  For
  a minimal by reverse inclusion polytope $\delta_i$ the set $W_i$ is
  the barycenter of $\delta_i$.
\vglue12pt

  Next, consider a complex of cones $\Omega=\{\omega_i \,;\, i\in I\}$
  referenced by a lattice~$X$. Again, we order the index set $I$ by
  the order {\it opposite} to the inclusion of cones.  The
  relationship between $|\Omega|$ and $BI=B\Omega$ in this case is
  slightly more complicated. Assume that $|\Omega|$ is connected.
  Since $\Omega$ is face-fitting, the minimal faces of all $\omega$'s
  are equal to a linear subspace $\omega_{\min}\subset X_{\bR}$.  If
  we denote $X'_{\bR}=X_{\bR}/\omega_{\min}$ then every $\omega_i$ is
  the preimage of a strictly convex cone $\omega'_i$ in~$X'_{\bR}$.

\demo{Definition {\rm 2.1.14}}
  Fix an arbitrary norm on $X'_{\bR}$ and let $S^{r'-1}$ be the unit
  sphere centered at the origin. The {\it spherical complex\/}
  $S\Omega=S\Omega'$ induced by $\Omega$ consists of the sets
  $\omega'_i\cap S^{r'-1}$. It will include the empty set by
  definition. 
\enddemo

2.1.15.
  Pick a point $p_i$ in the interior of each spherical cell
  $\omega'_i\cap S^{r'-1}$ and, in addition, a point $p_{\min}$ for
  the minimal cone. It is now clear that $B(I\setminus \min)$ is the
  barycentric decomposition of $|S\Omega|$ and that $B\Omega$ is the
  cone over it with the apex $p_{\min}$.

\demo{{R}emark {\rm 2.1.16}}
  It is clear that $|\Omega|$ and $|B\Omega|$ are homotopy equivalent
  and that the sheaf $\phi_g\inv F$ on $B\Omega$ has the same
  cohomologies as the natural constructible sheaf on $|\Omega|$.
\enddemo

\demo{Definition {\rm 2.1.17}}
  For a complex $\Delta$ of polytopes in $X$ 
  we define the complex $\Cone\, \Delta$ as  consisting of cones over all
  $\delta_i\subset (1,X_{\bR})$ in  
$\bX_{\bR}$, plus $\{0\}$.
\enddemo

  Note that the posets $\Cone\, \Delta$ and $\Delta$ (which includes
  $\emptyset$)
  are isomorphic, as well as are the posets $\Omega$
  and $S\Omega$, and that the topological spaces $|S\Cone\, \Delta|$ and
  $|\Delta|$ are homeomorphic.

\demo{{\rm 2.2.} Basic {\rm (}\/co\/{\rm )}\/sheaves and relations between them}
 \enddemo

{\it Definition} 2.2.1.
  The cosheaf $\uX$ on $\Omega_{\ord}$ is defined by associating to
  each cone $\omega_i$ the lattice $X_i=X\cap\bR\,\omega_i$, the
  saturated sublattice of $X$ generated by $\omega_i$. If $i>j$, i.e.\ $\omega_i$ is a face of $\omega_j$, then  
$X_i\subset X_j$,
  so this defines a cosheaf.  The dual sheaf $\hX=T$ is defined by the
  groups $\hX_i=T_i=\Hom(X_i,\bG_m(S))$.

\demo{Definition {\rm 2.2.2}}
  Let $\Delta$ be a complex of lattice polytopes referenced by $X$.
  For every polytope $\delta_i\subset(1,X_{\bR})$ let $\bX_i$ be the
  lattice $\bX\cap\bR\delta_i$. These groups define the cosheaf
  $\ubX$. The groups $\uhbX_i=\ubT_i=\Hom(\bX_i,\bG_m(S))$ define the
  dual sheaf $\uhbX=\ubT$.  Now, 
  $\Delta_{\ord}\simeq \Cone\, \Delta_{\ord}$, and the cosheaf $\ubX$ is
  the same as the cosheaf $\uX$ on $\Cone\, \Delta$ defined above.
\enddemo

{\it Definition} 2.2.3.
  For a pointed cell complex $(\Delta,C)$ we define a cosheaf
  $\uFun=\uFun[\Delta,C]$ of commutative groups by setting
  $\uFun_{i} = \Fun(C_i,\bZ)$.
  and the dual sheaf $\widehat{\uFun}=\widehat{\uFun}[\Delta,C](S)$ by
  $\widehat{\uFun}_{i}=\Fun(C_i,\bG_m(S))$.  Further, we define a
  cosheaf and the dual sheaf of
  semigroups by setting 
  $\uFun_{\ge0,i}(S)=\Fun(X\cap \delta_i\setminus C_i,\bZge)\oplus
  \Fun(C_i,\bZ)$ and
  $\uhFun_{\ge0,i}(S)=\Hom(\Fun_{\ge0,i},\cO_S(S)(S))$.  Also, define
  the cosheaves $\ubL,\ubK$ and the dual sheaves $\uhbL(S),\uhbK(S)$
  by
  \begin{eqnarray*}
    && \ubL_{i} = \ker\phi_{C_i}, \quad 
    \ubK_{i} = \coker\,\phi_{C_i} \\
    &&\mbox{where } \phi_{C_i}:\Fun(C_i,\bZ) \to \bX_{i},
    \quad (n_x) \mapsto \sum n_x (1,x) .
  \end{eqnarray*}
  All these (co)sheaves are relative to our choice of $(\Delta,C)$.
\vglue12pt

  Note that since we assumed that $\Conv\, C_i=\delta_i$, the groups
  $\bK_i$ are finite.  We have the following exact sequence of
  cosheaves on $\Delta$:$$ 
    0 \to \ubL \to \uFun \stackrel{\phi}{\to} 
    \ubX \to \ubK \to 0 .
  $$ 
  When $S=\Spec\, k$, $k=\bar k$, the dual sequence of sheaves is also
  exact:$$ 
    1\to \uhbK \to \uhbX \stackrel{\hat\phi}{\to} 
    \uhFun \to \uhbL \to 1    .
  $$

\specialnumber{2.2.4}\proclaim{Lemma}\label{2.2.4}
  $H_p(\Delta,\uFun)=0$ for $p>0$ and
  $H_0(\Delta,\uFun)=\Fun(C,\bZ)$. 
  Dually{\rm ,}
  $H^p(\Delta,\uhFun)=0$ for $p>0$ and
  $H^0(\Delta,\uhFun)=\Fun(C,\bG_m(S))$. 
\endproclaim

\demo{Proof}
  Indeed, the sheaf $\uhFun$ is a direct sum of sheaves $\uhFun_c$,
  $c\in C$, each of which is a constant sheaf on a contractible set
  $\overline{W}_c$.
\enddemo

\demo{Definition {\rm 2.2.5}}
  We introduce the complex $\bM_*=\bM_*[\Delta,C]$ as the mapping
  cylinder of the induced homomorphism of Cech chain complexes
  $\phi:\ChC_*(\Ustd,\uFun)\to \ChC_*(\Ustd,\ubX)$, or, if the poset
  $\Delta$ is a join semilattice, of Cech chain complexes with
  respect to the cover $\Umin$.  Explicitly, one
  has $\bM_p=\ChC_{p-1}(\uFun)\oplus \ChC_p(\ubX)$ with the
  differential $(d,(-1)^p\phi)\oplus d:\bM_{p+1}\to\bM_p$.
  
  Dually, we have a complex $\hbM^*=\hbM^*[\Delta,C]$ which is the
  cylinder of the map of Cech cochain complexes $\hat\phi:
  \ChC^*(\Ustd,\uhbX(S)) \to \ChC^*(\Ustd,\uhFun(S))$. Now
  $\hbM^p(S)=\Hom( \bM_p,\bG_m(S))$.
\enddemo

\demo{Definition {\rm 2.2.6}}
  By adding to each of the complexes above one more group --
  $\Fun(C,\bZ)$, $\bX$ etc., we obtain the {\it augmented complexes.}
  We will call their cohomology groups the {\it reduced homology
    groups \/} and denote them by $\oH_*$   (resp.\ for the reduced
  cohomologies $\oH^*$).  By taking the cylinder of the augmented
  chain complexes we obtain a new complex $\obM_*$ (resp.\ $\overline{\hbM^*}$).
\enddemo

  We will mostly be interested in the groups $H_p(\bM(S))$ and
  $H^p(\hbM(S))$ for $p=0,1$. Here is the basic tool for computing
  them:

\specialnumber{2.2.7}\proclaim{Lemma}\label{2.2.7} \hglue-7pt
  For the dual sheaves let $S\!=\!\Spec\, k${\rm ,} $k=\bar k$ in this
  statement\/{\rm :} 
  \begin{itemize}
  \ritem{1.} \u{\it H}\hskip.5pt$_p(\ubX)\simeq H_p(\bM_*)$ for $p\ge2$ and there is the
    following exact sequence\/{\rm :}
$$ 
      0\to H_1(\ubX) \to H_1(\bM_*) \to H_0(\uFun) 
      \to H_0(\ubX) \to H_0(\bM_*)      \to 0.
    $$ 
  \ritem{2.} Dually{\rm ,} $H^p(\hbM^*)\simeq \hbox{{\rm \u{\it H}}\hskip.5pt}^p(\uhbX)$ for $p\ge2$ and there is the following
exact sequence\/{\rm :}
$$ 
      1 \to H^0(\hbM^*) \to H^0(\uhbX) \to
      H^0(\uhFun) \to H^1(\hbM^*) \to H^1(\uhbX) \to 1.
    $$ 
  \ritem{3.} $H_0(\bM_*)\simeq H_0(\ubK)$ and there is the following exact
    sequence\/{\rm :}
$$ 
      H_2(\ubK)\to H_0(\ubL) \to H_1(\bM_*) \to H_1(\ubK) \to 0.
    $$ 
  \ritem{4.} Dually{\rm ,} $H^0(\hbM^*)\simeq H^0(\uhbK)$ and there is the
    following exact sequence\/{\rm :} 
$$ 
      1\to H^1(\uhbK) \to H^1(\hbM^*) \to
      H^0(\uhbL) \to H^2(\uhbK).
    $$ 
  \end{itemize}

\endproclaim

\demo{Proof}
  The short exact sequence of complexes$$ 
    0\to \ChC_*(\ubX) \to \bM_* \to \ChC_*(\uFun)[-1] \to 0
  $$ 
  of Cech chains gives rise to a long exact sequence on
  homologies. Part 1) follows from this sequence immediately by
  application the vanishing Lemma~\ref{2.2.4} and the
  fact that $H_p=\hbox{{\rm \u{\it H}}}_p$ for $p=0,1$. 

  Similarly, we have the following exact sequences:
  \begin{eqnarray*}
    && 0\to \Cone(\ChC_*(\uFun) \to \im\, \phi) \to
    \bM_* \to \ChC_*(\ubK) \to 0, \\    
    && 0\to \ChC_*(\ubL) \to \ChC_*(\uFun) \to \im\, \phi \to 0.
  \end{eqnarray*}
  The latter sequence implies that the complexes $\Cone(\ChC_*(\uFun)
  \to \im\, \phi)$ and $\ChC_*(\ubL)[-1]$ are quasi-isomorphic.
  Therefore, in the derived category we have a triangle$$ 
    \ChC_*(\ubL)[-1] \to \bM_* \to \ChC_*(\ubK) .
  $$ 
  Consequently, it also gives a long exact sequence of homologies.
  This, together with the fact that $H_2\to\hbox{{\rm \u{\it H}}}_2$ is surjective,
  gives 3). Parts 2) and 4) follow by a dual argument.
\enddemo

{\it Remark} 2.2.8.
  There are also similar exact sequences for the reduced
  (co)homologies which we skip for brevity.

\demo{Definition {\rm 2.2.9}} 
  Instead of the sheaf of groups $\uhFun$ in the definition of the
  complex $\hbM^*$ we can try to use the sheaf of semigroups for
  multiplication $\uhFun_{\ge0}$. The first two terms of the complex
  $\hbM_{\ge0}$ are$$ 
    \hbM_{\ge0}^0=\ChC^0(\uhbX) = \hbM^0
    \qquad \hbox{and} \qquad
    \hbM_{\ge0}^1=\ChC^1(\uhbX) \oplus \ChC^0(\uhFun_{\ge0}).
  $$ 
  The group $\hbM_{\ge0}^0$ naturally acts on $\hbM_{\ge0}^1$ according
  to the formula$$ 
    \hx^0.(\hx^1,\hf^0)=
    \left(d\hx^0\cdot\hx^1,\hx^0\hf^0\right).
  $$ 
  We will denote by $Z^1(\hbM_{\ge0}^*)$ the subset of $\hbM_{\ge0}^1$
  of pairs $(\hx^1,\hf^0)$ such that$$ 
    d\hx^1=0 \quad \hbox{and} \quad
    \hf^0_{i_1} = \hx^1_{i_1i_2} \hf^0_{i_2} \hbox{ in }
    \Fun(X_{i_1},\cO_S(S)) \hbox{ for all } i_1 > i_2.
  $$ 
\enddemo

 {\it Example} 2.2.10. 
  \label{exmp:homologies_Esherian_broken_ladder}
  Let $\Delta$ be the complex of lattice polytopes as on the picture.
  Explicitly, the coordinates of the vertices are $(0,0)$, $(0,4)$,
  $(4,0)$, $(1,1)$, $(1,2)$ and $(2,1)$.  From the long exact sequence
  of homologies associated to
  \begin{equation}\label{eqn:X_bX_Z}
    0\to \uX \to \ubX \to \bZ \to 0,
  \end{equation}
    taking into account that $H_p(\Delta,\bZ)=0$ for $p>0$ and
  $H_0(\Delta,\bZ)=\bZ$ (since $|\Delta|$ is contractible), we obtain
  $H_p(\Delta,\ubX)=H_p(\Delta,\uX)$ for $p>0$ and \pagebreak
  $H_0(\Delta,\ubX)=H_0(\Delta,\uX)\oplus \bZ$. An explicit
  computation of Cech homologies shows that\break
  $H_0(\Delta,\uX)\simeq\bZ^3$, $H_1(\Delta,\uX) \simeq\bZ$ and higher
  homologies vanish.
\figin{fig1}{700}
\enddemo

{\it Example}  2.2.11.
  Consider a complex $\Delta$ consisting of four squares meeting at the
  origin, with the outside vertices $(\pm1,\pm1)$, $(0,\pm2)$,
  $(\pm2,0)$. The group $C_1(\uX)$ is the direct sum of four copies of
  $\bZ$ corresponding to four common sides of the squares. The image
  $B_0$ of this group in $C_0$ contains the element
  $\{(2,0),(-2,0),(0,2),(0,-2)\}$ but not the half of it.  Therefore,
  $H_0(\uX)=C_0/B_0$ has a 2-torsion. If $\Delta^{(p)}$, $p>0$,
  denotes the direct product of $\Delta$ with $p$ copies of the
  complex from the example 2.2.10 then by
  the K\"uneth formula we see that
  $H_p(\Delta^{(p)},\uX)=H_p(\Delta^{(p)},\ubX)$ has 2-torsions.
\vglue12pt

2.3. {\it Affine stable toric varieties}.
  The first example is an ordinary affine toric variety over an
  algebraically closed field $k$. For a rational polyhedral cone
  $\omega\subset X_{\bR}$, which need not be maximal-dimensional, let
  $S_{\omega}=\omega\cap X$ be the saturated semigroup of lattice
  elements in $\omega$, and denote the semigroup algebra
  $k[S_{\omega}]$ by $R_{\omega}$.  Then $P_{\omega}=\Spec\, R_{\omega}$
  is an affine torus embedding, and it is well-known that it is normal
  (hence seminormal) and that our condition on the stabilizers is
  satisfied for the standard torus action $T\acts P_{\omega}$. We are
  now going to generalize this.

\vglue12pt 2.3.1.
  Let $\Omega$ be a cone complex referenced by $X$. For each cone
  $\omega_i$ we have an associated semigroup and  a semigroup algebra
  \begin{eqnarray*}
    S_i= X\cap\omega_i=X_i\cap\omega_i, \qquad 
    R_i=k[S_i] = \oplus_{s\in S_i} k\zeta^s .
  \end{eqnarray*}
  To each semigroup $S_i$ we will add one more element $\infty$ with
  the property $s+\infty=\infty$ for all $s\in S_i$.  We will make a
  convenient but nonstandard choice: in the corresponding semigroup
  algebra we will set formally $\zeta^{\infty}=0$.  To avoid confusion
  with the ordinary semigroup algebra $k[S_i\cup\infty]$ we will call
  this algebra $k_{\infty}[S_i\cup\infty]$.  Note that the collection
  $\{S_i\}$ is naturally a cosheaf on $\Omega_{\ord}$: for two cones
  $\omega_{i_1}$ and $\omega_{i_2}$ with $i_1>i_2$ (i.e.\ $\omega_{i_1}$ is a face of $\omega_{i_2}$) one has an embedding
  $S_{i_1}\into S_{i_2}$.  Since $S_{i_2}\setminus S_{i_1}$ is a
  semigroup ideal, we can make the collection $\{S_i\cup\infty\}$ into
  a sheaf on $\Omega_{\ord}$: for $i_1>i_2$ one has the Reese
  homomorphism $S_{i_2}\cup\infty\to S_{i_1}\cup\infty$ defined by$$ 
    S_{i_2}\ni s\mapsto \left\{
        \begin{array}{l@{\quad}l}
          s & \hbox{if }
          s\in S_{i_1}\\  
          \infty & \hbox{otherwise.}
        \end{array}
        \right.
  $$ 
  The Reese semigroup homomorphism induces an epimorphism of semigroup
  algebras ${\rm pr}_{i_2i_1}:R_{i_2} \to R_{i_1}$. We will
  denote these induced (co)sheaves of semigroups and algebras by
  $\uS[\Omega]$ and $\uR[\Omega]$ respectively.

\specialnumber{2.3.2}\proclaim{Lemma}
 {\rm 1.} In the category of sets ${\displaystyle\lim_{{\llar}}} (S_i\cup\infty) =
    ({\displaystyle\lim_{\lrar}}\, S_i)\cup\infty${\rm .}
\begin{itemize}
  \ritem{2.} For any ring $k$ (including $k=\bZ$) in the category of
    $k$-algebras{\rm ,}
$$ 
      H^0(\Omega,\uR[\Omega]) = 
      k_{\infty}[ H^0(\Omega,\uS[\Omega]) ].
    $$ 
    This is a free $k$\/{\rm -}\/module with the basis $\{\zeta^s,\
    s\in{\displaystyle\lim_{\longrightarrow}}\,  S_i\}$ and the multiplication$$ 
      \zeta^s\cdot\zeta^t = \left\{ 
        \begin{array}{l@{\qquad}l}
         \sum \zeta^{s+t} & \hbox{for all 
            $S_i$ so that $s,t\in S_i$} \\
          0 & \hbox{otherwise}.
        \end{array}
        \right.
    $$ 
  \end{itemize}

\endproclaim

\demo{Proof}
  All   parts are immediate from the definitions.
\enddemo

{\it Definition} 2.3.3.
  Define an algebra $R[\Omega]=H^0(\Omega,\uR[\Omega])$ and the affine
  variety $P[\Omega]$ as $\Spec\, R[\Omega]$.
\vglue4pt

    There is a natural torus action on $R[\Omega]$ and $P[\Omega]$: for
  each $t\in T=\Hom(X,\bG_m(k))$ and $u\in S_{\omega}$ one has
  $t.\zeta^u=t(\rho u)\zeta^u$.

\demo{Example {\rm 2.3.4}}
  Divide the line $\bR\supset\bZ$ into three cones: two half-lines
  $\omega_1$ and $\omega_2$ and the origin $\omega_{12}$. One has
  $R_{1}= k[x]$, $R_{2}= k[y]$, $R_{\omega_{12}}= k$.  Then
  $H^0(\Omega,\uR[\Omega])=k[x,y]/(xy)$. The variety $P[\Omega]$ is
  the union of two transversally intersecting lines.
\enddemo

{\it Remark} 2.3.5.
  In the case when the function $\rho:|\Omega|\to X_{\bR}$ is
  injective, the algebra $R[\Omega]$ was introduced by R. Stanley in
  \cite[\S4]{Sta2}.
\vglue4pt

  Next, we define a twisted version of the above. Consider a
  collection $\{t_{i_1i_0}\in T_{i_0}\}$ where $(i_1,i_0)$ goes over
  all pairs of indices with $i_0> i_1$ so that for all $i_0>i_1>i_2$
  one has$$ 
    t_{i_2i_0}=\varphi_{i_1i_0}(t_{i_2i_1})\cdot t_{i_1i_0}
    \hbox{ in } T_{i_0}
  $$ 
  (this can be supplemented by the elements $t_{ii}=1$). A collection
  like this is equivalent to a cocycle $t\in Z^1(\Ustd,\uT)$.  This
  cocycle defines a sheaf $\uR[\Omega,t]$ with the same stalks $R_i$
  and twisted homomorphisms$$ 
    p_{i_1i_0}=t_{i_1i_0} \cdot {\rm pr}_{i_1i_0}.
  $$

\demo{Definition {\rm 2.3.6}}
  Similarly, we define an algebra and a variety$$ 
    R[\Omega,t]= H^0(\Omega,\uR[\Omega,t]), \qquad
    P[\Omega,t]=\Spec\, R[\Omega,t].
  $$ 
  As before, it comes with a natural $T$-action.
\enddemo

  The algebras $R[\Omega,t]$ are a particular case of ``glued
  algebras'' studied by Dayton and Weibel in
  \cite{DW} and, on the other hand, of
  sections of rings of flasque (flabby) sheaves of integral algebras
  studied by Yuzvinsky in \cite{Yuz}. Thus,
  we can describe the precise structure of $R$'s and $P$'s by
  specializing and combining the quite general results of these
  papers. We note that Dayton and Weibel work with ``almost posets''
  which are slightly more general than what we need, and Yuzvinsky
  works with the semilattices which are slightly less general than what
  we need. In most cases, the general case can be deduced from the
  special case of a semilattice by a simple trick.

\demo{{R}emark {\rm 2.3.7}}
  If the reference function $\rho:|\Omega|\to X_{\bR}$ is one-sheeted, 
  i.e.\ injective, then $\Omega$ is a semilattice.
\enddemo

\specialnumber{2.3.8}\proclaim{Lemma}
  The sheaf $\uR[\Omega,t]$ is flabby {\rm (}\/flasque\/{\rm ).}\/
\endproclaim

\demo{Proof}
  If $\Omega$ is a semilattice, this is an application of the
  criterion in \cite[Cor.1.11]{Yuz} -- one
  has to check two properties for a semigroup algebra $k[S_{\omega}]$
  corresponding to a cone $\omega$ and a collection of ideals
  corresponding to its faces, which is immediate.
  
  The general case is reduced to the case of a semilattice by
  subdividing some cones into several ``dummy'' cones and assigning to
  these interior cones the same semigroups and identical connecting
  homomorphisms.
\enddemo

\specialnumber{2.3.9}\proclaim{{C}orollary}\label{lem:POmega_points_of}
  As a topological space{\rm ,}
$$ 
    P[\Omega,t] = {\displaystyle\lim_{\longrightarrow}}\,  \Spec\, k[S_{i}]
    = {\displaystyle\lim_{\longrightarrow}}\,  \Spec\, P_{i}.
  $$ 
  In particular{\rm ,} the irreducible components of $P[\Omega,t]$ are the
  toric varieties $P_{i_0}$ for maximal cones $\omega_{i_0}${\rm ,} they
  intersect along unions of smaller toric varieties $P_i$~etc\/{\rm .}
\endproclaim

\demo{Proof}
  Since the sheaf $\uR[\Omega,t]$ is flabby, for each $i$ there is a
  prime ideal $P_i$ such that $R_i=R[\Omega,t]/P_i$. With this in
  mind, \cite[Thm.1.2]{DW} and the
  previous lemma imply that the ring $R[\Omega,t]$ satisfies the
  condition ${\rm (CRT)_l}$ of
  [DW,\break p.\ 37] for all $l\ge2$.  The
  statement then is proved in
  \cite[Prop.\ 1.6]{DW}.
\enddemo

\specialnumber{2.3.10}\proclaim{{C}orollary}\label{lem:ROmega_is_seminormal}
  $R[\Omega,t]$ is seminormal{\rm .}  
\endproclaim

\demo{Proof}
  By \cite[Thm.1.8]{DW}.
\enddemo

\specialnumber{2.3.11}\proclaim{Lemma}
  $P[\Omega,t]=\Spec\, R[\Omega,t]$ is an affine stable toric
  variety {\rm (STV).} It is weakly normal as well{\rm .}
\endproclaim

\demo{Proof}
  $R[\Omega,t]$ is obviously reduced.  The only property that we
  have not yet checked is that the stabilizer of every point is
  connected and reduced. But every point lies in one of the
  irreducible components, and for normal toric varieties this property
  is well-known. For the last statement, consider a bijection $f:P'\to
  P$ inducing an isomorphism on the residue fields for all generic
  points. Since every irreducible component of $P$ is normal, $f$ has
  to be an isomorphism on each of them, and so it has to induce an
  isomorphism on the residue fields of each scheme point $x\in P$.
  Hence, the seminormality implies the weak normality in this case.
\enddemo

2.3.12.  Now let $P$ be {\it any} affine STV.  First of all, note that
  for each $p\in P$, the orbit $\orb\, p= Tp$ is isomorphic to the
  quotient torus $T_p=T/\Stab_p$. If $X_p\subset X$ is its group of
  characters then the stabilizer $\Stab_p$ is the multiplicative group
  dual to $X/X_p$. Hence, it is connected and reduced $\iiff$ $X/X_p$
  is torsion free; i.e.\ the sublattice $X_p\subset X$ is saturated.
 
\specialnumber{2.3.13}\proclaim{Lemma}\label{lem:irred_affine_STVs}
  Let $P$ be an irreducible affine variety with $T$\/{\rm -}\/action and only
  finitely many $T$\/{\rm -}\/orbits{\rm .} Assume that for each $p\in P$ the
  stabilizer $\Stab_p$ is connected and reduced{\rm .} Then\/{\rm :}
  \begin{itemize}
  \ritem{1.} The normalization $P^{\nu}$ is isomorphic to an affine torus
    embedding containing a quotient torus $T'${\rm .}
  \ritem{2.} $\nu:P^{\nu}\to P$ is bijective and induces isomorphisms on
    the residue fields of all scheme points{\rm .}
  \ritem{3.} If $P$ is seminormal then $P=P^{\nu}${\rm .}
  \end{itemize}
\endproclaim

\demo{Proof}
  Say $P=\Spec\, R$ for a finitely generated $k$-algebra $R$.  Let $p$
  be an arbitrary point in the dense orbit that is isomorphic to a
  torus $T'$ with a group of characters $X'$. The $T$-action on $R$
  defines an $X$-grading, and the embedding $T'\subset P$ defines an
  embedding $R\subset k[X']$ of $X$-graded algebras since $T'$ is
  dense.  Therefore, $R$ is isomorphic to a semigroup algebra $k[S]$,
  where $S$ is a finitely generated sub-semigroup of $X'$. The
  normalization $P^{\nu}$ obviously is $\Spec\, k[S^{\rm sat}]$, a torus
  embedding corresponding to a cone $\omega\subset X'_{\bR}$.  The
  $T$-orbits of $P$ are in a 1-to-1 correspondence with faces of this
  cone.  Since the stabilizers are connected and reduced, for each
  face $F$ the semigroup $S\cap F$ generates $F\cap X'$ as a group.
  Hence, $P^{\nu}\to P$ is bijective and induces an isomorphism
  $k(x)\simeq k(x')$ for every scheme point $x'\in X^{\nu}$ with
  $\nu(x')=x$.  It has to be an isomorphism if $P$ is seminormal.
\enddemo

  Therefore, an irreducible affine STV is simply an affine toric
  variety corresponding to a cone $\omega\subset X_{\bR}$.

\specialnumber{2.3.14}\proclaim{Theorem}\label{thm:every_affine_STV_is_such}
  Every affine {\rm STV} $P$ is isomorphic to one of $P[\Omega,t]${\rm .}
\endproclaim

\demo{Proof}
  $P$ is a union of finitely many $T$-orbits $\Orb_i$. Let $P_i$ be
  the closure of $\Orb_i$. It is irreducible, our condition on the
  stabilizers is satisfied,  and so by the above it defines a cone
  $\omega_i$ in $X$. A natural order by inclusion on $P_i$'s makes the
  index set into a poset. Further, since $P_i^{\nu}\to P_i$ is
  bijective and $T$-orbits on an affine toric variety correspond to
  faces of the cone, if $P_{i_1}\subset P_{i_2}$ then the cone
  $\omega_1$ is a face of $\omega_2$. An intersection of two $P_i$'s
  is a union of (closures) of other orbits, i.e.\ $P_j$'s.  Hence, the
  variety $P$ defines uniquely a complex of cones $\Omega$.
  
  Further, for each $i$ choose a point $p_i\in \Orb_i$, i.e.\ an origin
  in the subvariety~$P_i$.  This choice identifies the normalization
  $P_i^{\nu}$ with a standard affine toric variety $\Spec\, S_i$.  It
  also gives an origin in each lower dimensional $T$-orbit on $\Spec\, S_i$: it is $\lim_{1-{{\rm PS}}} p_{i}$ for any 1-parameter
  subgroup for which the limit lands in that orbit.  Since
  $P_i^{\nu}\to P_i$ is an isomorphism when restricted to any
  $T$-orbit, for all $i_1>i_2$ in $\Orb_{i_1}$ we have two ``origins''
  which differ by the torus action: for a uniquely defined $t_{i_1i_2}
  \in T_{i_1}$ one has$$ 
    \lim_{1-{{\rm PS}}} p_{i_2} = 
    t_{i_1i_2}^{-1} . p_{i_1}.
  $$ 
  The collection $\{t_{i_1i_2}\}$ is obviously compatible
  with respect to  triples $i_1>i_2>i_3$ and so is equivalent to a cocycle
  $t\in Z^1(\Ustd,\uT)$. 
  
  Say $P=\Spec\, R$. The collection of morphisms $P_i^{\nu}\to P$ is
  equivalent to a collection of homomorphisms $R\to R_i$ and hence to
  a homomorphism from $R$ to ${\displaystyle\lim_{{\llar}}}\, R_i= R[\Omega,t]$. By the
  construction and Lemma \ref{lem:POmega_points_of} the dual morphism 
  $P[\Omega,t]\to P$ is a bijection and by
  \ref{lem:irred_affine_STVs}.2 it induces an isomorphism on residue
  fields. Since $P$ is seminormal, it has to be an isomorphism.
\enddemo

2.3.15.
  This lemma establishes a 1-to-1 correspondence between the
  \linebreak ``marked'' affine STVs and the data $[\Omega,t^{1}]$. The
  standard origin $p_i$ on $\Spec\, k[S_i]$ is the same as a
  homomorphism $p_i:k[S_i]\to k$ sending each $\zeta^s$ to 1. Via the
  connecting homomorphism
  $p_{i_2i_1}=t_{i_2i_1}\cdot{\rm pr}_{i_2i_1}$ we see that
  $\lim_{1-{{\rm PS}}} p_{i_2}:t_{i_2i_1}.\zeta^s\mapsto~1$;
  hence $\lim_{1-{{\rm PS}}} p_{i_2} = t_{i_1i_2}^{-1} .
  p_{i_1}$.
 
\demo{Definition {\rm 2.3.16}}
  Consider the ``framed'' category $M^{\framed}[\Omega](k)$ of all
  affine STVs over a fixed algebraically closed field $k$
  corresponding to a given complex $\Omega$ with the action of a fixed
  torus $T$, in which we let the arrows be the isomorphisms giving the
  identity on $T$.  As a version, we can also consider the category
  $M[\Omega](k)$ where we allow arbitrary automorphisms $T\isoto T$
  (these are in a 1-to-1 correspondence with automorphisms of $X$,
  i.e.\ elements of $\GL(X)$.
\enddemo

  We will postulate them to be small categories.  Such a category is
  called a groupoid and it is given by two sets: $U$ -- the set of
  objects and $R$ -- the set of arrows and a map $j=(b,e):R\to U\times
  U$ associating to an arrow its beginning and end, and another map
  $e:U\to R$ associating to an object the identity arrow. These maps
  should satisfy the well-known axioms. This groupoid will be denoted
  $[U/R]$.  In particular, if $G$ is an abstract \pagebreak group acting on a set
  $U$ then $[U/G]$ will denote the groupoid obtained by taking
  $R=G\times U$ and $R\to U\times U$ defined as $(g,u)\mapsto
  (u,g.u)$.  The symbol $1$ will denote a\break $1$-point set.

\specialnumber{2.3.17}\proclaim{Lemma}
  If $\Omega$ consists of a single cone $\omega$ and its faces then
  the groupoid $M^{\framed}[\Omega](k)$ is equivalent to
  $[1/T_{\omega}(k)]$ and $M[\Omega](k)$ is equivalent to
  $[1/T_{\omega}(k)\ltimes \Sym\, \omega]${\rm ,} where the latter means the
  extension of the group $T_{\omega}(k)$ by the subgroup of $\GL(X)$
  of symmetries of the cone $\omega${\rm .}
\endproclaim

\demo{Proof}
  Trivial.
\enddemo

\specialnumber{2.3.18}\proclaim{Theorem}\label{thm:affine_STVs_isoclasses}
  In general{\rm ,} $M^{\framed}[\Omega](k)$ is equivalent to$$ 
    [Z^1(\Ustd,\uT)/C^0(\Ustd,\uT)] \sim
    [H^1(\Omega,\uT)/H^0(\Omega,\uT)]
  $$ 
  with the last group acting trivially{\rm .} In other words{\rm ,} the set of
  isomorphism classes of framed affine {\rm STVs} corresponding to a complex
  $\Omega$ is $H^1(\Omega,\uT)$ and for each $t${\rm ,} $\Aut
  P[\Omega,t]=H^0(\Omega,\uT)${\rm .}  The groupoid $M[\Omega](k)$ is
  equivalent to
$$ 
    [Z^1(\Ustd,\uT)/C^0(\Ustd,\uT)\ltimes\Sym\, \Omega] \sim
    [H^1(\Omega,\uT)/H^0(\Omega,\uT)\ltimes\Sym\, \Omega].
  $$ 
\endproclaim

\demo{Proof}
  By 2.3.15, the marked STVs over $k$ are
  classified by the set\break $Z^1(\Ustd,\uT)$ and the automorphism groups
  of marked STVs are trivial.  The groupoid of unmarked STVs is
  obtained by dividing this set by the equivalence relation induced by
  choosing different origins. That is obviously given by the action of
  $C^0(\Ustd,\uT)$, and the statement follows.
\enddemo

  As   mentioned above, the rings $R[\Omega,t]$ generalize the
  rings $R[\Omega]$ introduced by Stanley for a complex of cones lying
  {\it inside} $X_{\bR}$, and on the other hand -- the
  Stanley-Reisner rings of simplicial complexes. As in these
  particular cases, the singularities of $P[\Omega,t]=\Spec\, R[\Omega,t]$ strongly depend on the topology of the geometric
  realization $|\Omega|$. See Definition 1.1.30 for 
  the notions of CM and locally CM complexes.

\specialnumber{2.3.19}\proclaim{Theorem}\label{thm:CM_affine_case}
  Assume that the complex $\Omega$ is a semilattice{\rm .}  If $\Omega$ is
  {\rm CM} then $P[\Omega,t]$ is {\rm CM.}
\endproclaim

\demo{Proof}
  This is an application of a theorem of Yuzvinsky
  \cite[5.1]{Yuz} and the proof proceeds in
  the same way as that of \cite[4.6]{Sta2} which is a special
  case.
  
  For each $q>0$ and each $q$-dimensional cone $\omega$, choose
  $s_{\omega}\in\omega^0$. This defines an element
  $\zeta^{s_{\omega}}$ in $k[S_{\omega}]$ and so also in
  $R[\Omega,t]$.  Let $r_q=\sum_{\dim\omega=d+1-q}
  \zeta^{s_{\omega}}$, where $d=\dim\Omega$. Then $r_1,\dots, r_d$ is
  a standard system of parameters as defined in
  \cite{Yuz}, and the statement follows from
  Theorem 5.1. there.
\enddemo
\pagebreak

{\it Remark} 2.3.20.
  Assume that $\Omega$ is not a semilattice but is a quotient of
  another complex $\Omega'$ that is a semilattice and CM, by a finite
  group $G$ whose order is coprime to $\chr k$. Then $P[\Omega,t]=
  P[\Omega', t']/G$ and so is CM as well. I do not know whether this
  can be done for any CM complex of cones.

  On the other hand, it seems that the results of Yuzvinsky can be
  extended to the general (nonsemilattice) case.


\demo{{\rm 2.4.} Polarized stable toric varieties with a
  linearized line bundle}

\vglue4pt 2.4.1.  \hskip.4pt Affine torus embeddings over a closed field correspond in a \hbox{1-to-1}
  way to finitely generated rational polyhedral cones in $X_{\bR}$.
  Analogously, projective torus embeddings $(P,L)$ together with an
  extended action $T\acts L$ correspond 1-to-1   to lattice
  polytopes in $X_{\bR}$. Let us recall this connection.
  
  Let $\delta\subset X_{\bR}$ be a polytope with vertices in $X$. Let
  $\Cone\,\delta\subset \bX_{\bR}$ be the cone over $\delta$ lying in
  the hyperplane $(1,X_{\bR})$, and denote the semigroup $\bX
  \cap\Cone\,\delta$ by $S_{\delta}$.  Then $P_{\delta}=\Proj\,
  k[S_{\delta}]$ is a projective torus embedding and
  $L_{\delta}=\cO(1)$ is a linearized ample invertible sheaf.  One
  has a canonical identification$$ 
    \oplus_{d\ge0} H^0(P_{\delta},L_{\delta}^d) =
    k[S_{\delta}]= R_{\Cone\,\delta} .
  $$ 
  The affine cone $\Spec\, R_{\Cone\,\delta}$ is an affine toric variety
  for the action of the torus $\bT=\bG_m\oplus T$.  The complement of
  the apex of this cone is the $\bG_m$-torsor $\dL_{\delta}^{-1}$
  which is obtained from the line bundle $L^{-1}_{\delta}$, considered
  as a variety, by removing the zero section.
 
\demo{Definition {\rm 2.4.2}}
  Given a complex of lattice polytopes $\Delta$ referenced by $X$ and
  a cohomology class $[\tau]\in H^1(\Delta,\ubT)$ we define a
  projective variety $P[\Delta,\tau]$ as the $\Proj$ of
  $R[\Delta,\tau]=R[\Cone\, \Delta,\tau]$ and an ample sheaf
  $L[\Delta,\tau]$ on it as $\cO(1)$. The corresponding
  $\bG_m$-torsor will be denoted $\dL[\Delta,\tau]$.
\enddemo

\demo{{R}emark {\rm 2.4.3}}  
  The algebras $k[S_{\delta}]$ and $R[\Delta,\tau]$ are not
  necessarily generated in degree 1.  However, it is easy to see that
  each of them is finite over the subalgebra generated by elements of
  degree 1. In this situation $\cO(1)$ is an invertible sheaf; see, for example, \cite{Mor}.
\enddemo

\specialnumber{2.4.4}\proclaim{Theorem}
  The pair $(P,L)[\Delta,\tau]$ is a polarized stable toric variety
  with a $T$\/{\rm -}\/linearized action on $L${\rm .} It is weakly normal as well{\rm .}
\endproclaim

\demo{Proof}
  Take a polytope $\delta\subset(1,X_{\bR})$ and an element $u\in
  S_{\delta}$. It is easy to see that the localization of the algebra
  $R_{\delta}$ at $\zeta^u$ is the semigroup algebra $k[\bX\cap\Star_u
  \Cone\,\delta]$ and the subalgebra in the latter of elements of degree
  0 -- the semigroup algebra $k[X\cap\Star_u \Cone\,\delta]=
  k[S_{\Star_u\delta}]$.
  
  Given two polytopes $\delta_{i_1}>\delta_{i_2}$, a homomorphism
  $\tau_{i_2i_1}.{\rm pr}_{i_2i_1}: k[S_{\Cone\,\delta_{i_2}}]
  \to k[S_{\Cone\,\delta_{i_2}}]$ induces the homomorphism
  $t_{i_2i_1}.{\rm pr}_{i_2i_1}: k[S_{\Star_u\delta_{i_2}}]
  \to k[S_{\Star_u\delta_{i_2}}]$, where $t_{i_2i_1}$ is the image of
  $\tau_{i_2i_1}$ under the natural projection $\bT_{i_1}\to T_{i_1}$.
  Since the inverse limit commutes with the localization and 
  the subalgebra is of\break degree 0, we see that$$ 
    (R[\Delta,\tau][(\zeta^u)^{-1}])_0 = R[\Star_u\Delta,t].
  $$ 
  The algebra $R[\Delta,\tau]$ is finite over the subalgebra
  generated by $\zeta^u$ for all elements $u$ of degree 1 which appear 
  as vertices of some polytopes $\delta_i$. Therefore, $\Proj\, R[\Delta,\tau]$ is covered by finitely many affine STVs, and so is 
  a projective STV. The $\bX$-grading on $R[\Delta,\tau]$ defines
  the $\bT$-action on $L=\cO(1)$. This proves the statement.
\enddemo

\specialnumber{2.4.5}\proclaim{{C}orollary}\label{lem:PDelta_points_of}
  As a topological space{\rm ,}
$$ 
    P[\Delta,\tau] = {\displaystyle\lim_{\longrightarrow}}\,  P_i.
  $$ 
  In particular{\rm ,} the irreducible components of $P[\Delta,\tau]$ are
  the projective toric varieties $P_{i_0}$ for maximal polytopes
  $\delta_{i_0}${\rm ;} they intersect along smaller $P_{i}$\/{\rm '}\/s etc{\rm .} The
  restriction of $L[\Delta,\tau]$ to each $P_i$ is isomorphic to
  $L_i$.
\endproclaim

\demo{Proof}
  This is immediate from \ref{lem:POmega_points_of}.
\enddemo

2.4.6.
  Let $P$ be an arbitrary projective variety with the $T$-action and
  $L$ be an ample line bundle with the extended $T$-action. On the
  other hand, we have a natural scalar $\bG_m(k)$-action on $L$.
  Because $H^0(P,L)=k$, the only automorphisms of $L$ covering the
  identity on $P$ are in $\bG_m(k)$. Therefore, the two actions
  generate the action by a bigger group $\bT$ which is the extension
  of $\bG_m(k)$ by $T$. Every such extension is commutative, so
  $\bT=\bG_m(k)\oplus T$.

\specialnumber{2.4.7}\proclaim{Theorem}\label{thm:polarized_STVs}
  Every polarized stable toric variety $(P,L)$ with a linearized ample
  sheaf $L$ is isomorphic to one of polarized varieties
  $P[\Delta,\tau]$ for a complex $\Delta$ of lattice polytopes
  referenced by $X${\rm .}
\endproclaim

\demo{Proof}
  Consider the algebra $R=\oplus_{d\ge0}H^0(P,L^d)$. Then $P=\Proj\, R$,
  and $\Spec\, R$ is an affine cone over $P$. Note that the maximal
  ideal $R_+$ of elements of positive degree is not an embedded prime.
  Since $P$ is reduced, $R$ is reduced.  The algebra $R$ has a natural
  $\bT$-action. We claim that the seminormalization of $R$ is
  isomorphic to one of the algebras $R[\Delta,t]$. Indeed, $\Spec\, R^{{\rm sn}}\to \Spec\, R$ is an isomorphism except possibly
  at the apex since the complement of the apex is the $\bG_m$-torsor
  $\dL^{-1}$ over a seminormal variety. Hence, our condition on the
  orbits is satisfied and $\Spec\, R^{{\rm sn}}$ is an affine
  STV. Since the two affine varieties differ possibly only at the
  apex, it follows that $P=\Proj\, R=\Proj\, R[\Delta,t]$ and $L=\cO(1)$.
\enddemo

{\it Remark} 2.4.8.
  The results of the next subsection will imply that indeed $R\simeq
  R[\Delta,\tau]$.
\vglue12pt

2.4.9.
  As in the affine case, there is a unique isomorphism $P\isoto
  P[\Delta,\tau]$ if we mark the origins $p_i$ in each $\bT$-orbit on
  $P$, i.e.\ in each $\bG_m$-torsor $\dL_i^{-1}$ over $P_i$ (the
  $0$-dimensional orbit, the apex, can be omitted).

\demo{Definition {\rm 2.4.10}}
  For a complex $\Delta$ of lattice polytopes referenced by a lattice
  $X$, $M[\Delta](k)$ will denote the groupoid of polarized
  STVs {\it with a linearized line bundle} $L$ over $k=\bar k$ in
  which the arrows are the isomorphisms identical on $T$. As before,
  $M^{\framed}[\Delta](k)$ is a similar category in which we allow the
  isomorphisms which are not necessarily identities on $T$.
\enddemo

  We skip the proof of the following theorem since it is almost
  exactly the same as in the affine case.

\specialnumber{2.4.11}\proclaim{Theorem}\label{thm:polarized_linearized_STVs}
  The groupoid $M^{\framed}[\Delta](k)$ is equivalent to$$ 
    [Z^1(\Ustd,\ubT)/C^0(\Ustd,\ubT)] \sim
    [H^1(\Delta,\ubT)/H^0(\Delta,\ubT)]
  $$ 
  with the last group acting trivially{\rm .} In other words{\rm ,} the set of
  isomorphism classes of framed polarized {\rm STVs} with the linearized
  action corresponding to a complex $\Delta$ is $H^1(\Delta,\ubT)$ and
  for each $\tau$ one has $\Aut P[\Delta,\tau]=H^0(\Delta,\ubT)${\rm .}
  The groupoid $M[\Delta](k)$ is equivalent to$$ 
    [Z^1(\Ustd,\ubT)/C^0(\Ustd,\ubT)\ltimes\Sym\, \Delta] \sim
    [H^1(\Delta,\ubT)/H^0(\Delta,\ubT)\ltimes\Sym\, \Delta] .
  $$ 
\endproclaim

  Concerning the singularities, we have a statement similar to the affine case.

\specialnumber{2.4.12}\proclaim{Lemma}\label{thm:CM_proj_polarized_case}
  Assume that the complex $\Delta$ is a semilattice{\rm .}  If $\Delta$ is
  locally {\rm CM} then $P[\Omega,t]$ is {\rm CM.}
\endproclaim

\demo{Proof}
  Indeed, each cone complex $\Star_v\Delta$, $v\ne0$, in this case is
  CM.
\enddemo


2.5. {\it Cohomologies of polarized} STV{\it s}.
 
\specialnumber{2.5.1}\proclaim{Theorem}\label{thm:cohomologies_polarized_STVs}
  If $(P,L)=(P,L)[\Delta,\tau]$ then
  \begin{itemize}
  \ritem{1.} $H^p(P,L)=0$ for $p>0${\rm ;}
  \ritem{2.} $H^0(P,L)=\pi^{-1}(\tau)${\rm ,} where
    $\pi:Z^1(\hbM_{\ge0}[\Delta,\rho^{-1}X]) \to Z^1(\uhbX)$ is the
    natural projection{\rm ,} and $h^0(P,L)=\#\rho^{-1}X${\rm ,} the number of
    integral points in $|\Delta|${\rm .}
  \end{itemize} 

\endproclaim

\specialnumber{2.5.2}\proclaim{{C}orollary}
  For $(P,L)=(P,L)[\Delta,\tau]$ there is a canonical isomorphism of
  graded algebras$$ 
    \oplus_{d\ge0}H^0(P,L^d) = R[\Delta,\tau].
  $$ 
\endproclaim

  A similar result in a special case can be found in
  \cite{AN}. 
  The proof will be based on the following two lemmas.

\specialnumber{2.5.3}\proclaim{Lemma}
  For any subcomplex $\Omega'\subset\Omega$ {\rm (}\/resp.\ $\Delta'\subset\Delta${\rm ),}
 in other words for any open in
the order
  topology subset there exists  a closed affine {\rm (}\/resp.\ polarized\/{\rm )} subvariety
  $P[\Omega',t']\subset P[\Omega,t]$ {\rm (}\/resp.\ $P[\Delta',\tau']\subset
  P[\Delta,\tau]${\rm ).}
\endproclaim

\demo{Proof}
  The projective case follows from the affine since $P[\Delta,\tau]$
  is covered by $P[\Omega,t]$'s. The cocycle $t'$ is defined by
  restricting $t$: for each cone $\omega_{i'}\in\Omega'$ choose and
  fix a cone in $\Omega$ such that $\omega'\subset\omega$. If
  $\omega_{i'_1}\subset\omega_{i_1}$ and
  $\omega_{i'_2}\subset\omega_{i_2}$ then $t'_{i'_1i'_2}$ is the image
  of $t_{i_1i_2}$.  Define the map $f:R[\Omega,t]\to [\Omega',t']$ by
  the formula$$ 
    \zeta^u\mapsto
    \left\{ 
      \begin{array}{l@{\qquad}l}
        \zeta^u & \hbox{if } 
        u\in S_{\omega'_{i'}}\subset S_{\omega_i}\\
        0 & \hbox{if } u\notin |\Omega'|.
      \end{array}
    \right.
  $$ 
  Clearly, $f$ is well defined, is an algebra homomorphism and is
  surjective. This gives a closed embedding.
\enddemo

2.5.4.
  Now choose an open cover $\cU=\{U_{i_0} = U_{\ge i_0} \}$ of
  $\Omega_{\ord}$ (resp.\ of $\Delta_{\ord}$). For each of the open
  subsets $U_{i_0}\cap\cdots \cap U_{i_q}$ denote the corresponding
  closed subscheme by $P_{i_0\dots i_q}$. (Warning: if $\cU\ne\Umin$,
  the subschemes $P_{i_0}$ are not necessarily irreducible.) We will
  identify the sheaf $\cO_{P_{i_0\dots i_q}}$ with its direct image
  under the embedding $P_{i_0\dots i_q}\subset P$. Order the elements
  $\{i_s\}$ somehow (this has nothing to do with the partial order on
  $I$).  For each collection $i_0< \dots <i_q$ define a homomorphism
  $d: \cO_{P_{i_0\dots \hat i_s \dots i_q}} \to \cO_{P_{i_0\dots
      i_q}}$ as the restriction $\times (-1)^s$. Putting these all
  together, we obtain a complex$$ 
    0\to \cO_P \to \oplus \cO_{P_{i_0}} 
    \to \oplus \cO_{P_{i_0i_1}} \to \cdots
  $$ 
  in which clearly $d\circ d=0$.

\specialnumber{2.5.5} \proclaim{Lemma} 
  For $P=P[\Omega,t]$ {\rm (}\/resp.\ $P[\Delta,\tau]${\rm )} and any closed cover
  $\cU$ the above complex is an exact resolution of $\cO_P$ by the
  sheaves of $\cO_P$\/{\rm -}\/modules{\rm .}
\endproclaim

\demo{Proof}
  As before, the projective case follows from the affine. In the
  affine case, each of the sheaves is represented by an algebra
  splitting into a direct sum of $0$- or $1$-dimensional vector spaces
  $V_x$, one for each  $x\in\rho^{-1}X$. The exactness is sufficient to
  check on each of the $V_x$, which reduces to the acyclicity of the
  complex when we  compute the reduced cohomologies
  $\bar H^p(D_{n-1},k)$ of an $(n-1)$-dimensional disk using the
  Leray covering by simplices $D_{n-1}^{(i)}$, $i=0,\dots,n$.
  Hence, we have an exact
  sequence of commutative groups.  These groups are $R$-modules by the
  previous lemma.
\enddemo

\demo{Proof  of Theorem {\rm 2.5.1}}
  We will prove this statement by induction on the dimension and the
  number of minimal polytopes in $\Delta$, the case of one polytope,
  i.e.\ an ordinary projective toric variety being standard (see for example
  \cite{Oda}, \cite{Ful}).  In the above
  resolution take $\cU=\Umin$. The cohomologies of $L$ can be computed
  as the hypercohomologies of the above resolution twisted by $L$. By
  the induction assumption, all the higher cohomologies of
  $\cO_{P_{i_0\dots i_q}}\otimes L$ vanish, so that  we are reduced to the
  complex of $H^0$'s which are also known. Looking at each space
  $V_x$, we see that this complex is exact except for the first term.
  The kernel of the first homomorphism consists of the collections
  $\{f_{i_0}\in H^0(P_{i_0},L)\}$ such that for all $i_0,i_1$ one has
  in $H^0(P_{i_0i_1},L)$:
  \begin{eqnarray*}
    f_{i_0}t\inv_{i_0,i_0i_1}=f_{i_1}t\inv_{i_1,i_0i_1}, \qquad
    \hbox{i.e.} 
    \qquad f_{i_0}=f_{i_1}t\inv_{i_1i_0},
  \end{eqnarray*}
  which is precisely $\pi^{-1}(\tau)$.
\enddemo

  The same argument and the fact that on a projective toric variety
  higher cohomologies of the structure sheaf vanish give:

\specialnumber{2.5.6}\proclaim{Theorem}
  $H^p(P[\Delta,\tau],\cO)= H^p(|\Delta|,k)${\rm .}
\endproclaim

2.6. {\it Stable toric pairs}.
  A category of pairs $(P,\Theta)$ of varieties and Cartier divisors
  on them is equivalent to the category $(P,L,\theta)$, where
  $L=\cO(\Theta)$ and $0\ne\theta\in H^0(P,L)$ is an equation of
  $\Theta$. In this section we deal with the case when $L$ is
  $T$-linearized. To classify the stable toric pairs in this case we
  have to add to the data of the previous subsection an equation
  $\theta$ so that $\Theta$ does not contain any $T$-orbits entirely.
  We will start with the case when $P$ is irreducible, i.e.\ $(P,L)=(P_{\delta},L_{\delta})$. One has a canonical eigenspace
  decomposition $H^0(P_{\delta},L_{\delta}) = \mathbold{\oplus}_{x\in X\cap
    \delta} k \xi^{x}$, and so $\theta=\sum e_{x}\xi^x$.

\specialnumber{2.6.1}\proclaim{Lemma}
  $\Theta$ does not contain any $T$\/{\rm -}\/orbits entirely $\iiff$ $e_x\ne0$
  for the vertices of $\delta${\rm .}
\endproclaim

\demo{Proof}
  It suffices to consider only $0$-dimensional $T$-orbits, i.e.\ the
  fixed points of the $T$-action. These fixed points correspond in a
  1-to-1 way to the vertices $x_i\in X\cap \delta$. Say $x_0$ is one
  of them.  Then the section $\xi^{x_0}$ is the only one among $\xi^x$
  which does not vanish at that point. Therefore, one must have
  $e_{x_0}\ne0$.
\enddemo

{\it Definition} 2.6.2.
  $C(\Theta) = \{x\in X \,|\, e_x\ne0 \}$.

\specialnumber{2.6.3}\proclaim{{C}orollary}
  $\Theta$ does not contain any $T$\/{\rm -}\/orbits entirely $\iiff$
  ${\rm Vert}\, \delta\subset C(\Theta) \iiff \delta=\Conv\, C(\Theta)${\rm .}
\endproclaim

{\it Definition} 2.6.4.
  More generally, say $(P,L)=(P,L)[\Delta,\tau]$. Then for each
  $\theta\in R[\Delta,\tau]_1$ and for each $i$ we can consider its
  restriction on $R_i=k[S_i]$, and the set $C(\Theta)= \{x\in
  \rho^{-1}X \,|\, e_x\ne0 \}$ is well-defined.
 
\specialnumber{2.6.5}\proclaim{{C}orollary}
  If $(P,L)=(P,L)[\Delta,\tau]$ then $\Theta$ does not contain any
  $T$\/{\rm -}\/orbits entirely $\iiff$ ${\rm Vert}\, \Delta \subset
  C(\Theta)${\rm .}
\endproclaim

{\it Definition} 2.6.6.
  For a complex $\Delta$ of lattice polytopes referenced by a lattice
  $X$, $\MP^{\framed}[\Delta,C](k)$ will denote the groupoid of stable
  toric pairs $(P,\Theta)$ {\it with a linearized line bundle}
  $L=\cO(\Theta)$ over $k=\bar k$ in which the arrows are the
  isomorphisms identical on $T$. As before, $\MP[\Delta,C](k)$ is a 
  similar category in which we allow the isomorphisms which are not
  necessarily identities on~$T$.  Further, for a pointed complex
  $(\Delta,C)$ of lattice polytopes referenced by a lattice $X$ and
  with $\Ver\, \Delta\subset C$, $\MP^{\framed}[\Delta,C](k)$ and
  $\MP[\Delta,C](k)$ will denote the groupoids of pairs $(P,\Theta)$
  with $C(\Theta)=C$. Clearly,$$ 
    \MP[\Delta](k) = \coprod_{C,\, \Ver\, \Delta\subset C} 
    \MP[\Delta,C](k).
  $$

\specialnumber{2.6.7}\proclaim{Lemma}
  If $\Delta$ consists of the lattice polytope $\delta_i$ and its
  faces then
  \begin{eqnarray*}
    \MP^{\framed}[\Delta](k) &\sim&
    [ \hFun_{\ge0,i} / \bT_i ] ,
    \\  
    \MP[\Delta](k) &\sim&
    [ \hFun_{\ge0,i} / \bT_i \ltimes\Sym\, \delta] 
  \end{eqnarray*}
  {\rm (}\/see subsection {\rm 2.2} for definitions of these groups{\rm ).}  If{\rm ,} in addition{\rm ,} a set
  $C\supset\Ver\, \delta$ is chosen then
  \begin{eqnarray*}
    \MP^{\framed}[\Delta,C](k) &\sim&
    [ \hFun_{i} / \bT_i ] \sim [\hbL_i / \hbK_i ],
    \\  
    \MP[\Delta,C](k) &\sim&
    [ \hFun_{i} / \bT_i \ltimes\Sym\, \delta] 
    \sim [\hbL_i / \hbK_i \ltimes\Sym\, \delta]
  \end{eqnarray*}
  with $\hbK_i$ acting trivially{\rm .}
\endproclaim

\demo{Proof}
  All we have to do is to divide the set of the admissible $\theta$'s
  by the equivalent relation induced by the choice of the origin in
  the torsor $\dL_i^{-1}$. The latter is given by the action of the
  group $\bT_i$.
\enddemo

\specialnumber{2.6.8}\proclaim{Theorem}\label{thm:linearized_STpairs}
  \begin{eqnarray*}
    \MP^{\framed}[\Delta](k) &\sim&
    [ Z^1(\hbM^*_{\ge0}[\Delta,\Ver\, \Delta])/
      \hbM^0[\Delta] ] 
    \quad\hbox{and} \\
    \MP[\Delta](k) &\sim&
    [ Z^1(\hbM^*_{\ge0}[\Delta,\Ver\, \Delta]) /
    \hbM^0[\Delta] \ltimes\Sym\, \Delta] ]  . 
  \end{eqnarray*}
  For a pointed complex $(\Delta,C)$ with $\Ver\, \Delta\subset C${\rm ,}
  \begin{eqnarray*}
    \MP^{\framed}[\Delta,C](k) &\sim&
    [ Z^1(\hbM^*[\Delta,C])/
      \hbM^0[\Delta] ]     \\
    &\sim&
    [ H^1(\hbM^*[\Delta,C])/
      H^0(\hbM^*[\Delta,C])    ] ,     \\
    \MP[\Delta,C](k) &\sim&
    [ Z^1(\hbM^*[\Delta,C]) /
      \hbM^0[\Delta] \ltimes\Sym(\Delta,C)] ]   \\
    &\sim&
    [ H^1(\hbM^*[\Delta,C]) /
      H^0(\hbM^*[\Delta,C]) \ltimes\Sym(\Delta,C)] ]
  \end{eqnarray*}
  with $H^0(\hbM^*[\Delta,C])=H^0(\hbK[\Delta,C])$ acting trivially{\rm .}
  For every element\break  $l\in~H^0(\Delta,\uhbL[\Delta,C])$
  {\rm (}\/which describes a collection of pairs $(P_i,\Theta_i)${\rm )} whose image
  under the obstruction map$$ 
    H^0(\Delta,\uhbL[\Delta,C]) \to
    H^2(\Delta,\uhbK[\Delta,C])
  $$ 
  is zero{\rm ,} there is an associated finite set of stable toric pairs and
  this set is a torsor over the group $H^1(\Delta,\uhbK[\Delta,C])${\rm .}
  For each of these pairs$$ 
    \Aut(P,\Theta)= H^0(\hbM^*[\Delta,C])
    = H^0(\Delta,\uhbK[\Delta,C])
  $$ 
  is a finite group{\rm .} 
\endproclaim

\demo{Proof}
  A marked framed polarized STV $(P,L)[\Delta,\tau]$ does not have any
  automorphisms. By Theorem \ref{thm:cohomologies_polarized_STVs} a
  choice of $\theta$ is equivalent to a choice of an element of
  $Z^1(\hbM^*[\Delta,C])$. To describe the above groupoids we have to
  divide this set by the equivalence relation induced by choosing a
  different marking.  The second part of the theorem is an application
  of Lemma \ref{2.2.7}.
\enddemo

2.7. {\it The moment map for stable toric pairs over $\bC$}.
  Starting with a stable toric pair $(P,\Theta)$ over $k=\bC$, we are
  going to define a map $P(\bC)\to |\Delta|$ generalizing the usual
  moment map for toric varieties; see \cite{Oda} for
  example. This construction further shows how well the cell decomposition $\Delta$
  describes the essential properties of a stable toric pair.

\demo{Definition {\rm 2.7.1}}
  Let $\theta=\sum_{x\in C(\Theta)} \xi_x\in H^0(L)$ be an equation of
  $\Theta$.  The eigenfunctions $\xi_x$ are sections of the line
  bundle $L$ and at every point at least one of these functions is not
  zero.
  
  The {\it moment map\/} $P\to |\Delta|$ is defined by the formula$$ 
    p \mapsto 
    \frac{\sum |\xi_{x}(p)|^2 x}{\sum |\xi_{x}(p)|^2} \in \XR .
  $$ 
  Let $P_i$ be an irreducible component of $P$ corresponding to a
  polytope $\delta_i$. The sections $\xi_x$ with $x\not\in \delta$ are
  identically zero on $P_i$. Therefore, restriction of the above map
  to $P_i$ is a usual moment map for a projective toric variety. Its
  image is precisely the polytope $\delta$.
\enddemo

  Since stable toric varieties are ``built'' of ordinary projective
  toric varieties, the next result easily follows from the description
  of the moment map for toric varieties.

\specialnumber{2.7.2}\proclaim{Theorem}\label{thm:moment_map_STPs}
  \hiha {\rm 1.} The moment map defines a continuous map between the
    topological spaces $P(\bC)$ and $|\Delta|${\rm .} The fiber over a point
    in the interior of a $q$-dimensional cell $\delta$ is
    $S^q=U(1)^q${\rm .} \begin{itemize}
  \ritem{2.} Choose a set of representatives for $\xi_x${\rm '}s {\rm (}\/they are
    defined up to a common complex constant\/{\rm ).} Then the restriction of
    the moment map $P \to |\Delta|$ to the subset $P(\bR_{\ge0})$ of
    points with nonnegative and real $\xi_x$ is a homeomorphism{\rm .}
  \end{itemize}

\endproclaim

{\it Definition} 2.7.3.
  By analogy with the usual terminology for toric varieties, we can
  call the $C^{\infty}$-varieties locally isomorphic to
  $P[\Delta,|\tau|](\bR_{\ge0})$ (i.e.\ locally isomorphic to
  $P[\Omega,|t|](\bR_{\ge0})$ {\it varieties with sharp corners.\/}
 
\vglue12pt 
2.8. {\it One-parameter families of stable toric pairs}.
  Let $(\cR,\fm)$ be a DVR with the fraction field $\cK$, $k$ be the
  residue field $\cR/\fm$, $\eta=\Spec\cK$ and $O=\Spec\, k$ be the
  generic and the special points of $\Spec\cR$ respectively.  We will
  denote $\cS=\Spec\cR$.  We will also denote by $s$ the uniformizing
  parameter, i.e., the generator of the ideal $\fm$.

\specialnumber{2.8.1}\proclaim{Theorem}\label{thm:STPs_1paramfamily}
  Let $T_{\eta}$ be a torus over $\cK$ and $(P_{\eta},\Theta_{\eta})$
  be a family over $\cK$ of stable toric pairs with the linearized
  $T_{\eta}$\/{\rm -}\/action on $L_{\eta}=\cO_{P_{\eta}}(\Theta_{\eta})${\rm .} Then
  after a finite{\rm ,} possibly ramified cover this family can be completed
  to a family of stable toric pairs $(P,\Theta)$ over $\cR${\rm ,} and such
  a completion is unique up to an isomorphism{\rm .}
\endproclaim

{\it Reduction Step} 1.
  After an \'etale base change we can assume that the torus $T_{\eta}$
  is split with the character group $X\isoto\bZ^r$. We will continue
  to use the same letters $\cR,\cK$ etc. in order not to crowd the
  notation. The $T_{\eta}$ and the $\bG_{m,\eta}$-action on
  $L_{\eta}$ give the action of a split torus $\bT_{\eta}$ with the
  character group $\bX=\bZ\oplus X$. Now $T=\Spec\cR[X]$ and
  $\bT=\Spec\cR[\bX]$ will denote the corresponding split tori over
  $\cR$.
\vglue12pt

  We will start the proof of the theorem with the case when $P_{\eta}$
  is geometrically irreducible, i.e.\ when the pair
  $(P_{\eta},L_{\eta})\otimes \overline{\cK}$ corresponds to a single
  lattice polytope $Q\subset X_{\bR}$. 
\vglue12pt

{\it Reduction Step} 2.
  After another \'etale base change we can assume that $P_{\eta}$
  admits a $\cK$-point in the ``interior'', i.e., a point that over
  the algebraic closure $\overline{\cK}$ induces a point in the main
  torus orbit.
\vglue12pt

  Consider the $K$-algebra
  $R_{\eta}=\mathbold{\oplus}_{d\ge0}H^0(P_{\eta},L_{\eta}^d)$. The $\cK$-point
  gives an embedding $T_{\eta}\to P_{\eta}$ and, consequently, the
  embedding $R_{\eta}\subset \cK[\bX]$, which allows us to fix an
  isomorphism $R_{\eta}\simeq \cK[\Cone\, Q]$. We have an element
  $\theta_{\eta}\in H^0(P_{\eta},L_{\eta})$ which can be uniquely written as$$ 
    \theta_{\eta}=\sum\xi_{\eta,x}
    =    \sum_{x\in X\cap Q} c(x)\, \zeta^{(1,x)}, \qquad
    c(x)\in \cK .
  $$ 
  As before, let $C=\{x \,|\, c(x)\ne0 \}$.  The discrete valuation on
  $\cK$ defines an integral-valued function on $C$: $\psi: x\mapsto
  \val\, c(x)$, and now
$$ 
    c(x)= c'(x)\cdot s^{\psi(x)}  \qquad \hbox{with }
    \val\, c'(x)=0.
  $$ 
  We will denote by $c_0(x)$ the residue of $c'(x)$ in $k$.

\demo{Definition {\rm 2.8.2}}
  Consider the convex hull ${\rm ConvHull}_{\psi}$ in
  $X_{\bR}\oplus\bR $ of the rays$$ 
    \{(x,h)
    \,|\, x\in C, \,    h\ge\psi(x)\}     .
  $$ 
  The graph of the lower envelope of this convex hull is a
  piecewise-linear function $g_{\psi}: Q\to\bR$.  Thinking of $X$ here
  as the subset $(1,X)$ of $\bX$ of degree 1, we can extend this
  function to a function on $\Cone\, Q$ by dilations: $g_{\psi}(\chi)= d
  g_{\psi}(\chi/d)$, where $d=\deg\chi$. The fields of linearity of
  $g_{\psi}$ define a cell decomposition $\Delta$ of the polytope $
  Q$. For each cell $\delta_i$ in this decomposition let $C_i$ be the
  subset of $X\cap\delta_i$ of the points with
  $\psi(1,x)=g_{\psi}(1,x)$, i.e.\ the points $x$ for which the
  corresponding point in $\bX\oplus\bZ$ lies on the lower boundary of
  ${\rm ConvHull}_{\psi}$.
\enddemo

  A pointed decomposition $(\Delta_0,C_0)=\{(\delta_i,C_i)\}$ defined
  by the above lifting property is usually called {\it regular} or
  {\it coherent}, see \cite[Ch.\ 7]{GKZ}.

\demo{Definition {\rm 2.8.3}}
  Define the $\cR$-subalgebra $R\subset\cK[\bX]$ as the
  $\cR$-module with the generators $\{s^h\zeta^{\chi} \,|\, \chi\in
  \Cone\, Q, \, h\ge g_{\psi}(\chi)\}$.  Although $R$ is not a semigroup
  algebra, it formally corresponds to the semigroup of integral
  elements in $\bZ\oplus\bX$ lying inside the cone over
  ${\rm ConvHull}_{\psi}$. Since this is a rational finitely
  generated cone, the algebra $R$ is finitely generated over $\cR$.
\enddemo

  Now we wish to look at the residue algebra $R_0=R/\fm R$.

\specialnumber{2.8.4}\proclaim{Lemma}
  The algebra $R_0$ is nilpotent\/{\rm -}\/free if and only if the function
  $g_{\psi}$ is integral\/{\rm -}\/valued on $\bX\cap\Cone\, Q${\rm .} In this case{\rm ,}
  $R_0$ is isomorphic to the standard algebra $R[\Delta,1]$ over the
  field $k${\rm .}
\endproclaim

\demo{Proof}
  For each $\chi\in\bX\cap\Cone\, Q$ the corresponding eigenspace in
  $R_0$ is a 1-dimensional $k$-space with the generator
  $\bar\zeta_{\chi}$ which is the residue of $s^{\ulcorner
    g_{\psi}(\chi) \urcorner}\zeta^{\chi}$ ($\ulcorner\urcorner$
  denotes the upper integral part). It is clear that for the points
  $(1,x)$ corresponding to the edges of the lower envelope,
  $g_{\psi}(\chi)$ is integral, but in general this may not be true.
  If for some $x$ the value $g_{\psi}(\chi)$ is not integral then a
  power of $\bar\zeta{\chi}$ corresponds to an element lying in the
  interior of ${\rm ConvHull}_{\psi}$ that can be represented as
  the sum of $\bar\zeta$'s for the edges plus $(0,n)$ with $n>0$. Thus,
  in the residue algebra this power is zero, and we have a nilpotent element.
  
  Now assume that all the points on the lower envelope of the cone
  over ${\rm ConvHull}_{\psi}$ are integral and consider the
  product $\bar\zeta_{\chi_1}\bar\zeta_{\chi_2}$. Geometrically, it
  corresponds to the sum of the corresponding two elements on the
  lower envelope of $\Cone({\rm ConvHull}_{\psi})$. If
  $\chi_i$ belong to the same cone, we get a point on the same face of
  the envelope, and so
  $\bar\zeta_{\chi_1}\bar\zeta_{\chi_2}=\bar\zeta_{\chi_1+\chi_2}$.
  Otherwise, we get a point in the interior of
  $\Cone({\rm ConvHull}_{\psi})$ which corresponds to an
  element on the lower envelope plus $(0,n)$ with $n>0$. In this case,
  $\bar\zeta_{\chi_1}\bar\zeta_{\chi_2}=0$. This is precisely the
  description of the algebra $R_0[\Delta,1]$ in Section
  2.3. By the results of
  section~2.3 the algebra $R_0$
  in this case is nilpotent-free.
\enddemo

{\it Reduction Step} 3.
  The condition that $g_{\psi}$ is integral-valued on
  $\bX\cap\Cone\, Q$ can be achieved by a ramified base change
  $s=t^n$. Make this base change.
\vglue12pt

  We are now ready to complete the proof of Theorem
  \ref{thm:STPs_1paramfamily} in the case of geometrically irreducible
  $P_{\eta}$. Simply consider $P=\Proj_{\cR}R$ for the above-defined
  algebra $R$. This algebra is a free module over $\cR$; hence $P$ is
  flat over $\Spec\, \cR$. We have established that the central fiber is
  an STV. The line bundle $L=\cO(1)$ comes with a global section
  $\theta=\sum c(x)\zeta^{(1,x)}$. By the construction, on the central
  fiber it restricts to a section of $H^0(P_0,\cO(1))$ with the
  equation$$ 
    \theta_0= \sum_{x\in C} c_0 \bar\zeta^{(1,x)}.
  $$ 
  It does not contain $0$-dimensional toric orbits since for each
  vertex of a polytope $\delta_i$ the corresponding coefficient
  $c_0(x)$ is nonzero. This proves the existence.

\vglue12pt 2.8.5.
  To summarize, if the generic fiber is geometrically irreducible then
  the central fiber $(P_0,\Theta_0)$ corresponds to the following
  data:
  \begin{itemize}
  \ritem{1.} The decomposition $(\Delta_0,C_0)$ defined by the
    integral-valued function $\psi$,
  \ritem{2.} The trivial 1-cohomology class $\tau_0=1\in
    Z^1(\Delta_0,\uhbX)$,
  \ritem{3.} The element $\hf_0=(c_0(x))\in
    C_0(\Delta_0,\uhFun[\Delta_0,C_0])$.
  \end{itemize}
  The pair $(\tau_0,\hf_0)$ is an element of the group
  $Z^1(\hbM^*[\Delta_0,C_0])$ as   defined earlier.
\vglue12pt

  Let us prove the uniqueness of the pointed cell decomposition for
  the central fiber $(P_0,L_0)$.  Let $f:(P,\Theta)\to\Spec\, R$ be an
  extended flat family.  As above, after an additional base change
  (which does not change the type of the central fiber) we can make the
  tori split and put the generic fiber in the standard form. By the
  results of   Section~2.5
  on each geometric fiber the higher cohomologies of $L^d$, $d>0$,
  vanish.  Therefore, by cohomology and base change the
  $\cO_{\cS}$-algebra $f_*L^d$ is locally free. The $\bT$-action gives
  a canonical splitting of this algebra into the direct sum of
  $\bX$-eigenspaces which therefore also have to be locally free and
  hence invertible. Let $\theta=\sum\xi_x\in R$ be the equation of the
  divisor $\Theta$. For each $f=(n_x)\in\Fun(C,\bN_{\ge0})$
  define $\xi^f=\prod \xi_x^{n_x}$.
  
  Our condition on the divisor $\Theta_0$ implies that the algebra $R$
  is finite over the subalgebra $\Rp$ generated by $\xi^f$'s, in fact
  for every $\bX$-homogeneous element $r\in R$ we have
  $r^n=\hbox{unit}\cdot \xi^f$ for some $n>0$. The cell decomposition
  $\Delta$ depends only on  whether for each pair of
  homogeneous elements $r_1,r_2$ the residue of $(r_1r_2)^n$ is zero
  or not, i.e.\ whether $\xi^{f_1+f_2}$ equals zero or not.  On the
  other hand, by finiteness for every $\xi^f$ one has
  $\xi^{nf}=c(nf,g)\xi^g$ for some $g$ such that the residue of
  $\xi^g$ is not zero. Finally, the element $c(nf,g)\in \cK$ can be
  computed on the generic fiber as $\xi^{nf}_{\eta}/\xi^g_{\eta}$.
  This proves the uniqueness of the decomposition.  It also shows the
  uniqueness of the $\cO_S$-subalgebra $\oplus f_*L^d$ of
  $\cK[\bX]$, and therefore of its $\Proj$, i.e.\ the family
  $P$. Clearly, the divisor $\Theta$, which is the closure of
  $\Theta_{\eta}$, is unique as well.

  In the general case, the pair
  $(P_{\eta},\Theta_{\eta})\otimes\bar\cK$ corresponds to a complex
  $(\Delta,C)=\{(Q_j,C_j)\}$.  After the \'etale base change, we
  represent the algebra $R_{\eta}$ in the standard form
  $R_{\eta}[\Delta,\tau]$, with $\tau\in \Hom(C_1(\ubX),\cK^*)$ and
  the equation $\theta$ as the system of equations $\theta_j$ which
  give an element $\hf\in \Hom\left(C_0(\uFun[\Delta,C]),\cK\right)$.
  The pair $(\hf,\tau)$ has to be an element of $Z^1(\hbM^*)$ as
  defined earlier.
  
  For each of the algebras $R_{\eta,j}$ repeat the above procedure,
  defining a finitely generated $\cR$-subalgebra $R_j$ in each. One
  then checks that the gluing functions between neighboring
  $\zeta_{ji}$ are units in $\cR$.  Therefore, the image of the
  collection $R_j$ in $R_{\eta}$ is an $\cR$-algebra of the type
  $R[\Delta_0,C_0]$ for the pair $(\tau',\hf')\in Z^1(\hbM^*)$ with the
  coefficients in $\cR$, as opposed to simply in $\cK$. The central
  fiber will then correspond to $(R_0,\theta_0)$ for the pair
  $(\tau_0,\hf_0)$ which is the residue of $(\tau',\hf')$ in $k$.

  The uniqueness of the decomposition and the subalgebra $\oplus
  f_*L^d$ follows because by the special case we have the uniqueness
  of decomposition and the subalgebras for each $Q_j$. The theorem is
  now completely proved.

\vglue12pt 2.8.6.
  To summarize, the central fiber $(P_0,\Theta_0)$ corresponds to the
  following data:
  \begin{itemize}
  \ritem{1.} The sub decomposition $(\Delta_0,C_0)$ of $(\Delta,C)$ defined
    by the integral-valued functions~$\psi_j$,
  \ritem{2.} The 1-cohomology class $\tau_0\in Z^1(\Delta_0,\uhbX)$ that is
    the residue of the image of the class $\tau'$ under the natural
    morphism $Z^1(\Delta,\uhbX) \to Z^1(\Delta_0,\uhbX)$,
  \ritem{3.} The element $\hf_0\in C^0(\Delta_0,\uhFun[\Delta_0,C_0])$
    which is the residue of the image of the element $\hf'$ under the
    natural morphism $C^0(\uhFun[\Delta,C])\to
    C^0(\uhFun[\Delta_0,C_0])$.
  \end{itemize}
  The pair $(\tau_0,\hf_0)$ is an element of the group
  $Z^1(\hbM^*[\Delta_0,C_0])$ as was defined earlier.
\vglue12pt

2.9. {\it The modified complex $\bM_*$}.
   Starting with a pointed complex $(\Delta,C)$ we are now going to
  define the ``semigroup'' analogs of the chain complexes $C_*(\ubX)$,
  $C_*(\uFun)$ and $\bM_*$. We will only need them in degrees $0,1$ and
  $2$.

\demo{Definition {\rm 2.9.1}}
  We will define $SC_0(\ubX\lge)=\oplus \bX_i$. The natural
  projections $\Cone\,\delta_i\to \Cone|\Delta|$ extend by linearity to
  a map $\varrho:\oplus\Cone\,\delta_i \to \Cone|\Delta|$.  In degree
  $n>0$ we will define $SC_n(\ubX\lge)$ as a free abelian semigroup on
  the generators $(\bar\chi_0,\chi_1,\dots,\chi_n)$ with \pagebreak
  \begin{eqnarray*}
    & \bar\chi_0\in \oplus\Cone\,\delta_i,
    \quad  \chi_1,\dots,\chi_n\in \cup\Cone\,\delta_i 
    \subset \oplus\Cone\,\delta_i, \\
    & \hbox{such that } 
    \varrho\bar\chi_0=\varrho\chi_1=\dots=\varrho\chi_n ,
  \end{eqnarray*}
  modulo a certain natural equivalence relation. We will not even
  spell it out in degrees $n\ge2$ since we do not need it. In degree 1
  it has the form
  \begin{eqnarray*}
    &(\bar\chi_1+\bar\chi_2,\chi)=(\bar\chi_1,\chi_1)+
    (\bar\chi_2,\chi_2)+(\chi_1+\chi_2,\chi), \\
    &(\chi,\chi)=0,
  \end{eqnarray*}
  which reduces to the linearity condition$$ 
    (\chi'_1+\chi'_2,\chi_1+\chi_2)=(\chi'_1,\chi_1)+
    (\chi'_2,\chi_2)
  $$ 
  when $\chi_1,\chi_2$ belong to the same cone.
  
  We will define the differential $SC_1(\ubX\lge)\to SC_0(\ubX\lge)$
  by $(\bar\chi_0,\chi_1)\mapsto \bar\chi_0-\chi_1$. The image of the differential from $SC_2(\ubX\lge)$ to
$SC_1(\ubX\lge)$ will be
  interpreted as an equivalence relation on $SC_1(\ubX\lge)$ generated 
  by$$ 
    (\chi_1,\chi_2) + (\bar\chi_0,\chi_1) = (\bar\chi_0,\chi_2).
  $$ 
  Hence, $SC_1/SB_1(\ubX\lge)$ will denote the quotient of the
  semigroup $SC_1(\ubX\lge)$ by this equivalence relation.  
\enddemo

{\it Remark} 2.9.2.
  There is a natural ring homomorphism $k[SC_1/SB_1(\ubX\lge)] \to
  k[C_1/B_1(\ubX)]$. It sends every $x^{(\sum \bar\chi_i,\chi)} $ with
  not all $\bar\chi_i$ belonging to a common cone to 0. This defines
  an open embedding on the $\Spec$'s of these algebras in the opposite 
  direction.

\demo{Definition {\rm 2.9.3}}
  For the cosheaf $\uFun$ we define $SC_0(\uFun\lge)=C_0(\uFun\lge)$.
  $SC_1(\uFun\lge)$ will be the free abelian group with the generators
  $(f_0,f_1)$ with $\varrho\phi f_0=\varrho\phi f_1$.
\enddemo

{\it Definition} 2.9.4.
  We define the complex $\bSM_*$ as follows:
  \begin{eqnarray*}
    && \bSM_0=SC_0(\ubX_{\ge0}), \quad
    \bSM_1=SC_0(\uFun\lge) \oplus SC_1(\ubX_{\ge0}), \\
    && \bSM_2=SC_1(\uFun\lge) \oplus SC_2(\ubX_{\ge0}).
  \end{eqnarray*}
  The differential $\bSM_1\to\bSM_0$ is defined as one should expect,
  $(f_0,(\bar\chi_0,\chi_1))\mapsto \phi f_0+\bar\chi_0-\chi_1$. The
  image of the differential $\bSM_2\to\bSM_1$ again will be
  interpreted as an equivalence relation generated by$$ 
    f_0 = f_1 + (\phi f_0,\phi f_1) \hbox{ and }
    (\chi_1,\chi_2) + (\bar\chi_0,\chi_1) = (\bar\chi_0,\chi_2).
  $$ 
 
\specialnumber{2.9.5}\proclaim{Lemma}
  The semigroups $SC_1/SB_1(\ubX)$ and $\bSM_1/B_1(\bSM_*)$ are
  finitely generated and finitely presented{\rm .}
\endproclaim \pagebreak

\demo{Proof}
  We start with $SC_1/SB_1(\ubX)$. If $\chi_1,\chi_2$ are in the same
  $\Cone\,\delta_i$ then the relations imply that
  $(\bar\chi_1+\bar\chi_2,\chi_1+\chi_2)=
  (\bar\chi_1,\chi_1)+(\bar\chi_2,\chi_2)$.  Hence the symbols
  $(\bar\chi,\chi)$ for a fixed $\delta_i$ form a semigroup.  This is
  the kernel of a homomorphism from a semigroup of integral elements
  in a finitely generated rational cone in the space
  $\left((\oplus_j\Cone\,\delta_j)\oplus\Cone\,\delta_i\right)\otimes \bR$
  to $\bXR$.  Therefore, it is finitely generated. The semigroup
  $SC_1/SB_1(\ubX)$ is generated by these semigroups for different
  cells $\delta_i$, hence it is finitely generated as well. 
  
  Adding $C_0(\uFun)$ together with our relations does not add more
  new generators than the number of generators in $H_0(\uFun)$, which
  is definitely finite.  Finally, by a standard result of the
  semigroup theory every finitely generated commutative semigroup is
  finitely presented.
\enddemo

\specialnumber{2.9.6} \proclaim{Lemma}\label{lem:semigroups_finite_action}
  Let $H$ and $L$ be a finitely generated commutative semigroup and
  group respectively{\rm ,} and $f:H\to L$ a homomorphism defining an action
  on the scheme $V=\Spec\, H$ by the dual split group of multiplicative
  type $T=\Spec\bZ[L]${\rm .} If the image of the invertible elements of $H$
  is a subgroup of finite index in $L$ then the action morphism
  $T\times V\to V\times V$ is finite{\rm .} In other words{\rm ,} in this case the
  $T$\/{\rm -}\/action is proper with finite stabilizers{\rm .} If the image of the
  invertible elements is all of $L$ then the $T$-action is free{\rm ;}
  i.e.\ $T\times V\to V\times V$ is a closed embedding{\rm .}
\endproclaim

\demo{Proof}
  The morphism $T\times V\to V\times V$ is defined by the semigroup
  homomorphism $H\times H\to L\times H$, $(h_1,h_2)\mapsto (fh_1,
  h_1+h_2)$. Take any $(l,h)\in L\times H$. If $nl=fh_1$ with $h_1$
  invertible then $(h_1,nh-h_1)\mapsto n(l,h)$.
\enddemo

\specialnumber{2.9.7}\proclaim{{C}orollary}\label{cor:main_action_finite}
  The action of the torus $\Spec\, \bZ[\bSM_0]=\Spec\, \bZ[C_0(\ubX)]$ on
  the scheme $\Spec\, \bZ[\bSM_1/B_1(\bSM_*)]$ is proper with finite
  stabilizers{\rm .}
\endproclaim

\demo{Proof}
  Indeed, in this case already the image of the invertible elements of 
  $C_0(\uFun_{\ge0})$ in $C_0(\ubX)$ has finite index.
\enddemo

{\it Definition} 2.9.8.
  For each $\Delta'\ge\Delta$ the ideal $I[\Delta']\subset
  \bZ[\bSM_1/B_1(\bSM_*)]$ is   generated by
  $(\bar\chi,\chi)=(\sum\chi_i,\chi)$ such that there is no common
  polytope $\delta$ so that all $\chi_i,\chi\in \Cone\,\delta$.
  Similarly, for each $(\Delta',C')\ge(\Delta,C)$ we will denote by
  $I[\Delta',C']$ the ideal generated by these elements plus
  $f_{x,i}\in \Fun_i$ for all $x\in X\cap \delta \setminus C_i$.

\vglue12pt 2.10. {\it Moduli of stable toric pairs}.
 
\specialnumber{2.10.1} \proclaim{Lemma}\label{lem:type_is_constant}
  Let $f:(P,\Theta)\to S$ be a flat family of stable toric pairs with
  a linearized action of a split torus $T/S$ over a connected locally
  Noetherian base $S${\rm .}  If one geometric fiber $(P_s,\Theta_s)$
  corresponds to a cell decomposition of a polytope $Q\subset\bR$ then
  any other geometric fiber corresponds to a cell decomposition of the
  same polytope $Q${\rm .}
\endproclaim

\demo{Proof}
  Since on each geometric fiber the higher cohomologies of $L^d$,
  $d>0$, vanish, by cohomology and base change the sheaves $f_*L^d$ are 
  locally free. The action by $\bT=\bG_m\oplus T$ defines a
  decomposition of $\oplus_{d\ge0}f_*L^d$ into a direct sum of
  eigenspaces over all $\chi\in\bX$. They also must be locally free,
  and since $S$ is connected -- of a constant rank. By looking at the
  fiber $(P_s,\Theta_s)$ we see that each of these eigenspaces is
  1-dimensional if $\chi\in\Cone\, Q$ and $0$-dimensional
  otherwise. Therefore, the same is true for any other fiber, and this 
  property defines the type of the fiber completely.
\enddemo

{\it Definition} 2.10.2.
  We will refer to the above as the {\it convex $1$\/{\rm -}\/sheeted case}.

\demo{{R}emark {\rm 2.10.3}}
  For the purposes of moduli in this paper we will restrict ourselves
  to this case only (and its analogs when the action on $L$ is not
  linearized or when there is an abelian part). We call that a
  general linearized stable toric pair corresponds to a complex
  $\rho:|\Delta|\to\XR$ which is quite general: the image in $\XR$
  need not be convex and the map can have several sheets. The moduli
  of such pairs are more complicated and will be left for another
  occasion.
\enddemo

{\it Definition} 2.10.4.
  Fix a lattice $X$ and the corresponding split torus $T=\Spec\bZ[X]$.
  Denote by $\cTP^{\framed}$ the moduli stack on the category of
  locally Noetherian schemes associating to each scheme $S$ the
  groupoid of flat families $(P,\Theta)/S$ together with a linearized
  action of $T_S$. This groupoid is equivalent to the groupoid of flat
  families $(P,L,\theta)$ of linearized polarized toric varieties
  together with a section $\theta\in \Gamma(S,f_*L)$.
  
  By the above, for a fixed polytope $Q\subset\XR$ we have a connected
  substack $\cTP^{\framed}[Q]$ for which every geometric fiber in a
  family corresponds to a cell decomposition of $Q$.

\demo{Definition {\rm 2.10.5}}
  For a pointed cell decomposition $(\Delta,C)$ let
  $\cTP^{\framed}[\Delta,C]$ be a substack in which for each geometric
  fiber the corresponding cell decomposition is $\ge(\Delta,C)$.
\enddemo

\specialnumber{2.10.6}\proclaim{Lemma}
  $\cTP^{\framed}[\Delta,C]$ is open in $\cTP^{\framed}[Q]${\rm .}
\endproclaim

\demo{Proof}
  For a geometric fiber the condition that the cell decomposition is
  $\ge(\Delta,C)$ is the condition that for the generators of
  eigenspaces $R_{\chi_i}$ with $\chi$'s in the same $\Cone\,\delta$
  ($R=\oplus_{d\ge0}f_*L^d$) the products are not zero and that
  residues of homogeneous components of $\theta$ for all $x\in C$ are
  not zero. These all are open conditions.
\enddemo

\specialnumber{2.10.7}\proclaim{Lemma}
  The union of all $\cTP^{\framed}[\Delta,C]$ covers
  $\cTP^{\framed}[Q]${\rm .}
\endproclaim

\demo{Proof}
 This is obvious since every fiber induces some $(\Delta,C)$.
\enddemo

\specialnumber{2.10.8}\proclaim{Theorem}
  $\cTP^{\framed}[\Delta,C]$ is a separated Artin stack of finite
  presentation over $\bZ$ with finite stabilizers{\rm .} It has a coarse
  moduli space $\TP^{\framed}[\Delta,C]$ which is a separated scheme
  over $\bZ${\rm .}
\endproclaim

\demo{Proof}
  We claim that $\cTP^{\framed}[\Delta,C]$ coincides with the quotient
  stack$$ 
    \left[\Spec\, \bZ[\bSM_1/B_1(\bSM_*)] /
      \Spec\, \bZ[\bSM_0]
      \right]
  $$  
  for the action defined in the previous subsection where we checked
  that the stabilizer of the action is finite. The existence of a
  coarse moduli space for such a stack as a separated algebraic space
  is the main result of [KM].
  
  Indeed, let $(P,\Theta)/S$ be any flat family of stable toric
  pairs. First of all, a family like this is equivalent to a family
  $(P,L,\theta)$, $\theta\in H^0(P,L)$ satisfying the necessary
  conditions on geometric fibers. Secondly, we claim that it is
  equivalent to giving a locally free graded $\cO_S$-algebra
  $(R,\theta)$ with an element of degree 1, again with appropriate
  conditions on geometric fibers. Given $(P,L,\theta)$, we take
  $R=\oplus_{d\ge0} f_*L^d$. Vice versa, given $R$, the pair $(P,L)$
  is recovered as $(\Proj\, R,\cO(1))$. Hence, the moduli of our pairs
  are essentially the moduli of algebras $(R,\theta)$.

  In an \'etale neighborhood $\Spec\, D\to S$ of any point
  locally in \'etale topology we can find the generators in the
  subalgebras $R_i\subset R=\oplus_{d\ge0}f_*L^d$ corresponding to
  each cell $\delta_i\in \Delta$ and identify the $D$-algebra $R$ with
  the union of semigroup algebras $R_i$. One sees that the relation 
  among these elements is precisely a homomorphism from the
  semigroup algebra $\bZ[SC_1/SB_1(\ubX)]$ to $D$, and that the pair
  $(R,\theta)$ uniquely defines a homomorphism from the semigroup
  algebra $\bZ[\bSM_1/B_1(\bSM_*)]$ to $D$. Hence, locally we get a
  $D$-point of the scheme $\Spec\, \bZ[\bSM_1/B_1(\bSM_*)]$. To get the
  moduli stack $\cTP^{\framed}[\Delta,C]$ we have to divide this scheme 
  by the choices  made, and that amounts to the torus
  action by $\Spec\, \bZ[\bSM_0]$.
  
  By [KM] there is a coarse moduli space as an algebraic
  space. However, for a proper torus action with finite
  stabilizers the quotient of an affine scheme is simply the affine
  scheme corresponding to the subring of invariants. Therefore, in the
  case at hand it is the scheme corresponding to the semigroup algebra
  $\bZ\left[K[\Delta,C]\right]$, where $K=\ker (\bSM_1/B_1(\bSM_*) \to
  \bSM_0)$.
\enddemo

{\it Remark} 2.10.9.
  Note the role of the ideals $I[\Delta']$ and $I[\Delta',C]$. The
  first one defines a closed subscheme over which the fibers in our
  family are varieties of type $\le\Delta'$ and $\ge\Delta$. The second
  defines a subscheme over which the pair $(P,\Theta)$ have the type
  $\le(\Delta',C')$ and $\ge(\Delta,C)$.

\specialnumber{2.10.10}\proclaim{Theorem}\label{thm:moduli_finite_case}
  $\cTP^{\framed}[Q]$ is a proper Artin stack of finite presentation
  over $\bZ$ with finite stabilizers{\rm .} It has a coarse moduli space
  $\TP^{\framed}[Q]$ which is a proper scheme over $\bZ${\rm .}
\endproclaim

{\it Proof}.
  By the result of Section 2.8, any 1-parameter family over a complete
  discrete valuation field $\cK$ can be completed to a family over the
  DVR after a finite base change. Therefore, $\cTP^{\framed}[Q]$ is
  complete. Moreover, the decomposition $(\Delta_0,C_0)$ for the
  central fiber is uniquely defined, so  that the limit is contained in a
  unique minimal $\cTP^{\framed}[\Delta_0,C_0]$. Since the latter is
  separated, $\cTP^{\framed}[Q]$ is separated. The coarse moduli space
  is obtained by gluing the coarse moduli spaces of all the
  $\cTP^{\framed}[\Delta,C]$'s.
\hfill\qed\vglue4pt

2.11. {\it Approximation of the moduli space}.  Here, we are going to show that the just constructed moduli space
  $M=\TP^{\framed}[Q]$ admits a finite morphism to an easily
  describable projective scheme which can be considered to be its
  simplification. Consequently, $\TP^{\framed}[Q]$ is itself
  projective. Moreover, both schemes admit toric stratifications in
  which every stratum in one scheme corresponds in a 1-to-1 way to an
  isogenic stratum in the second. We may say that our moduli space and
  its simplification differ in only finitely many ways.

\vglue4pt {\it Definition} 2.11.1.
  Denote $C_{\max}=X\cap Q$ and let $\phi$ be as before the natural
  semigroup homomorphism from $\Fun(C_{\max},\bZ_+)$ to $\Cone\, Q
  \subset \bX$ sending $1_x$ to $(1,x)$. Take a rational point $q\in
  (1,Q)$ and let $N$ be the minimal natural number with $Nq$ integral.
  Let $H_q$ be the semigroup $\phi\inv(\bZ_+ Nq)$.
\vglue4pt

{\it Definition} 2.11.2.
  On the other hand, consider the map $\Fun(C_{\max},\bR_+)\cap \{(f_x)
  \,|\, \sum f_x=1\}$ to $(1,Q)$. This is a map of polytopes: from the
  simplex on the vertices $1_x$ of dimension $\# C_{\max} -1$, which
  we will denote $\sigma$ to $Q$. Let $\sigma_q$ be the fiber of this
  polytopal map over $q$.
\vglue4pt

  It is clear that $H_q$ is a saturated sub-semigroup in
  $\Fun(C_{\max},\bZ)$ consisting of all integral points in the cone
  over $\sigma_q$.

\specialnumber{2.11.3}\proclaim{Lemma}\label{lem:coords_vertices_fiber}
  Faces of $\sigma_q$ are in a one\/{\rm -}\/to\/{\rm -}\/one correspondence with lattice
  subpolytopes $\delta$ of $Q$ containing $q$ in their interior{\rm .}
  Vertices correspond to simplices with this property{\rm .} If
  $\delta=\langle x_0,\dots, x_r\rangle $ is such a simplex and $(1,q)=\sum c_i
  (1,x_i)$ then the coordinates of the corresponding vertex are
  $(c_i)${\rm .} 
\endproclaim

\demo{Proof} 
  This is truly obvious.
\enddemo

{\it Definition} 2.11.4.
  $F_q=\Proj\, H_q$.
\vglue12pt

2.11.5.
  So, $F_q$ is a projective normal toric scheme corresponding to the
  lattice polytope obtained by a dilation of $\sigma_q$ (we are not
  interested in the ample linearized sheaf on $F_q$, only in the
  scheme itself). If two fibers $\sigma_{q_1},\sigma_{q_2}$ are
  combinatorially equivalent and have parallel faces  then
  $F_{q_1}=F_{q_2}$ -- since the dual fans are the same. If
  $\sigma_{q_2}$ is obtained by ``degenerating'' $\sigma_{q_1}$, i.e.\ by contracting some of the faces, then there is a
natural morphism
  $F_{q_1}\to F_{q_2}$; this corresponds to a natural map of fans.
  It is birational if $\sigma_{q_i}$ have the same dimension, for
  example if both $q_i$ lie in the interior of $Q$.
  
  The polytope $Q$ is naturally subdivided into locally closed
  polytopal domains $D_{\alpha}$ such that for $q$'s in the same
  $D_{\alpha}$ the fibers are equivalent, and for $D_{\alpha_2}\subset
  \overline{D_{\alpha_1}}$, fibers $\sigma_{q_2}$ are degenerations of
  fibers $\sigma_{q_1}$. This stratification is defined by
  intersecting the projections of faces of $\sigma$. Indeed, for two
  points in $D_{\alpha}$ the cells $\delta$ containing these points in 
  their interior are exactly the same.
  
  Fix one scheme $F_{\alpha}$ for every stratum $D_{\alpha}$. We have
  obtained a partially ordered set $J={\alpha}$, ordered by the
  ``degeneration relation'' and a system of projective schemes
  $F_{\alpha}$ indexed by $J$ with morphisms $F_{\alpha_1}\to
  F_{\alpha_2}$ for each $\alpha_1\ge \alpha_2$.
\vglue12pt

{\it Definition} 2.11.6.
  Define $M_{\simp}={\displaystyle\lim_{{\llar}}}_J F_{\alpha}$.

\vglue12pt

{\it Remark} 2.11.7.
  Note that in the category of schemes over a fixed base inverse
  limits over any finite category $J$ exist. Indeed, by a standard
  result from   category theory it would suffice for finite products
  and equalizers to exist, and both are particular cases of fiber
  products. It is also easy to show (although we do not need it) that
  ${\displaystyle\lim_{{\llar}}}_J F_{\alpha}$ is a closed subscheme of
  $\prod_{J_{\max}} F_{\alpha}$. We also note that this scheme was
  previously considered in \cite[Sec.4]{KSZ}.
 
\vglue12pt
{\elevensc Theorem 2.11.8.}
  {\it There are natural morphisms $M\to F_{\alpha}$ compatible with the
  morphisms $F_{\alpha_1}\to F_{\alpha_2}$ and{\rm ,} hence{\rm ,} a morphism}
  $f:M\to M_{\simp}${\rm .}
\vglue12pt

{\it Proof}.
  Fix $q\in Q$. For any family $(P,L,\theta)/S$ of stable toric pairs
  we will construct a map $S\to F_q$ in a functorial way. That will
  induce the morphism $M\to F_q$ from the moduli space. For any point
  $s\in S$ we will construct an affine neighborhood $\Spec\, R\ni s$
  and a map $\Spec\, R\to F_q=\Proj\, \bZ[H_q]=\Proj\, \bZ[H_q]^{(d)}$,
  where $\bZ[H_q]^{(d)}$ is the subring of the graded ring $\bZ[H_q]$
  of elements of degrees divisible by $dNq$, for some $d>0$.  This will be
  induced by a ring homomorphism $\bZ[H_q]^{(d)}\to R$. Write
  $\theta\in H^0(P,L)$ as the sum of its homogeneous components:
  $\theta =\sum_{x\in C_{\max}} \xi_x$.  For $f=(f_x)$, $g=(g_x)$ with
  $\phi f =\phi g$ the sections $\xi^f=\prod \xi_x^{f_x}$ and $\xi^g$
  both lie in the same eigenspace of $H^0(P,L^d)$, $d=\deg f$. Let
  $s\in S$ be a point and $\kappa(s)$ be its residue field. Then on
  the fiber $P_s$ one of the residues of $\xi^h$, $h\in H_q$, $\deg
  h>0$, must be nonzero: $q$ lies in a cell $\delta$ of the
  decomposition $\Delta$ induced by the fiber $(P_s,L_s,\theta_s)$,
  and for this cell for any $\chi\in\Cone\, \delta$ there is an $h$ with
  $\phi h = k\chi$ and $\xi^h\ne0$: this is another way of saying that
  the group $\bK_{\delta}$ is finite. Let $d$ be minimal such that
  there is a $\xi^h$ with $\phi h=dNq$ and nonzero residue. Choose a
  small affine neighborhood $\Spec\, R'\ni s$ and a generator $e$ in
  $H^0(\Spec\, R',L)_{dNq}=R'$ so that we can write $\xi^f= c(f)e^k$.
  Further, for one fixed $h$ let $\Spec\, R$ be a smaller neighborhood
  so that for any $s'\in \Spec\, R$ the residue of $c(h)$ in
  $\kappa(s')$ is not zero. Then the rule $f\mapsto c(f)$ defines the
  required homomorphism $\bZ[H_q]\to R$ inducing the morphism $\Spec\, R\to F_q$.  This morphism clearly
\pagebreak does not depend on a choice of
  generator $e$, hence is the same on intersections, hence defines a
  unique morphism $S\to F_q$, as required. Compatibility with maps
  $F_{\alpha_1}\to F_{\alpha_2}$ is easily shown.
\hfill\qed

\vglue6pt
  Now, let us fix a field $k$ and describe the $k$-points of
  $M_{\simp}(k)$.

\specialnumber{2.11.9}\proclaim{Lemma}
  $M_{\simp}(k)$ is the disjoint union of $H^0(\Delta,\uhbL)$ going
  over all subdivisions of $(Q,C_{\max})${\rm .}
\endproclaim

{\it Proof}.
  By the definition of $k$-points and inverse limit, $M_{\simp}(k)=
  {\displaystyle\lim_{{\llar}}} F_{\alpha}(k)$. A point $p$ of the toric variety $F_q$
  lies in one of the toric orbits which, in turn corresponds to a face
  of the polytope $\sigma_q$. It is easy to see that faces of
  $\sigma_q$ are in a 1-to-1 correspondence with cells $\delta\subset
  Q$ with $q$ in the interior, and that the corresponding toric orbit
  is $\hbL_{\delta}$. If $\delta_2$ is a face of $\delta_1$ and
  $q_2\in \delta_2^0$ is a point then the image of $p$ under
  $F_{q_1}(k) \to F_{q_2}(k)$ is given by the homomorphism $\hbL_1\to
  \hbL_2$. Hence, we see that a point of ${\displaystyle\lim_{{\llar}}} F_{\alpha}(k)$
  is given precisely by fixing a cell decomposition $\Delta$ and an
  element of ${\displaystyle\lim_{{\llar}}} \hbL_i = H^0(\Delta,\uhbL)$.
\hfill \qed

\specialnumber{2.11.10}\proclaim{{C}orollary}
  $M$ and $M_{\simp}$ have natural stratifications with strata in a
  $1$\/{\rm -}\/to\/{\rm -}\/$1$ correspondence with decompositions $\Delta$ of $Q${\rm .} For every 
  $\Delta$ these strata are isogenous{\rm .}
\endproclaim

{\it Proof}.
  By Lemma \ref{2.2.7} the map $H^1(\hbM^*) \to H^0(\uhbL)$ has
  finite kernel and cokernel.
\hfill\qed

\specialnumber{2.11.11}\proclaim{{C}orollary}\label{cor:moduli_projective_finite_case}
  The morphism $M\to M_{\simp}$ is finite{\rm .} $M$ is projective over the
  base scheme{\rm .}
\endproclaim 

{\it Proof}.
  Indeed, $M_{\simp}$ is clearly projective: it is obtained from
  several projective schemes by taking fiber products. Hence, both
  schemes are proper and so the morphism is proper (our base scheme is
  separated). Since all geometric fibers are finite, it is finite, and
  so $M$ is also projective.
\hfill\qed \vglue6pt

{\it Definition} 2.11.12.
  For every decomposition $(\Delta,C)$ of $(Q,Q\cap X)$ there is a
  natural open affine subscheme in $M_{\simp}$ which is the spectrum
  of a semigroup algebra. Indeed, in each $F_{q}$ the cell $\delta\in
  (\Delta,C)$ containing $q$ in its interior defines a saturated sub-semigroup $K_{\delta,q}$ in $\bL$ defined by the
differences
  $\chi_1-\chi_2$ with $\chi_2\in\Cone\, \delta$. This gives an open
  affine subscheme $\Spec\,\bZ[K_{\delta,\alpha}]$ in $F_{\alpha}$. In
  the inverse limit we get an open affine subscheme \smallbreak
\centerline{${\displaystyle
    U[\Delta,C] = 
    {\displaystyle\lim_{{\llar}}} \Spec\, \bZ[K_{\delta,\alpha}] =
    \Spec\, \bZ[ {\displaystyle\lim_{\longrightarrow}}\,  K_{\delta,\alpha} ].
 }$ }
\smallbreak\noindent
  As $(\Delta,C)$ goes over all triangulations of $(Q,Q\cap X)$ we
  obtain an open cover of $M_{\simp}$. 
\pagebreak

2.12. {\it Generalized secondary polytopes for finite complexes}.

\vglue9pt 2.12.1.
  All polytopes $\sigma_{\alpha}$ lie in spaces which are parallel to
  $\bL=$ \linebreak $ \ker (\Fun(C_{\max},\bZ)\to \bX)$. Consider the
  Minkowski sum $\Sigma(Q)=\Sigma(Q,C_{\max})=\sum \sigma_{\alpha}$.
  The projective toric variety $F_{\Sigma(Q)}$ corresponding to this
  polytope and the lattice $\bL$ maps to every $F_{\alpha}$ since
  every $\sigma_{\alpha}$ can be obtained from the Minkowski sum by
  degenerating and collapsing some faces. Therefore, it maps to the
  inverse limit, i.e.\ $M_{\simp}$.

\specialnumber{2.12.2}\proclaim{Theorem}
  $F_{\Sigma(Q)}$ maps isomorphically onto an irreducible component of 
  $M_{\simp}${\rm .}
\endproclaim

{\it Proof}.
  In   semigroup terms, the morphism $f:f\inv U[\Delta,C] \to
  U[\Delta,C]$ corresponds to a semigroup homomorphism ${\displaystyle\lim_{\longrightarrow}}\,  K_{\delta,q} \to \bL$ and since every $K_{\delta,q}$ is saturated in
  $\bL$, the image is also saturated. The semigroup image is what
  describes $\im\, f$.  We note that $\im\, F_{\Sigma(Q)}$ maps to an irreducible
  component of $M_{\simp}$ since both contain an open dense subset
  $\Spec\, \bZ[\bL]$.
\hfill\qed

\specialnumber{2.12.3}\proclaim{{C}orollary}\label{cor:main_comp_moduli_finite_case}
  $F_{\Sigma(Q)}$ maps isomorphically onto an irreducible component of
  $M${\rm .}
\endproclaim

\demo{Proof} 
  Let $M'=M[\Delta_{\max},C_{\max}]$ be the reduced part of the main
  irreducible component of $M$. Both $M'$ and $F_{\Sigma(Q)}$ have
  finite morphisms to the same component $M'_{\simp}$ of $M_{\simp}$.
  Since $F_{\Sec(Q)}$ is normal, it is the normalization of $M'$. On
  the other hand, the composition $F_{\Sec(Q)}\stackrel{\psi}{\to}
  M'\to M'_{\simp}$ is an isomorphism; hence $\psi$ is also an
  isomorphism. 
\enddemo

  The Minkowski average of the fibers of $\sigma\to Q$ is precisely the
  {\it secondary polytope} of Gelfand-Kapranov-Zelevinsky. There is a
  more ``balanced'' version of it which is combinatorially equivalent
  to it and has parallel faces (hence, corresponds to the same toric
  variety) given by the following explicit formula:

\vglue12pt
{\it Definition} {\rm 2.12.4}.
  For each triangulation $\Tr=(\Delta,C)$ of $(Q,X\cap Q)$, i.e.\ a
  pointed decomposition for which all the cells are simplices and $C$
  consists precisely of the vertices of these simplices, let
  $\phi_{\Tr}$ be an element of $\Fun(C,\bZ)$ such that
  \begin{equation}\label{eqn:vertices_of_secpol}
    \phi_{\Tr}(x)= \sum_{\delta_i: x\in\delta_i} \Vol\, \delta_i. \speqnu{2}
  \end{equation}
  
  $\Sec(Q)=\Sec(\Delta_{\max},C_{\max})$ is the convex hull in the
  space $\Fun(C,\bR)$ of the points $\phi_{\Tr}$ as ${\Tr}$ goes over
  all triangulations.
 \pagebreak

{\it Definition} 2.12.5.
  Alternatively, $\Sec(Q)$ is the Minkowski
  integral, dilated by $(r+1)$, of the fibers of all fibers $\sigma_q$, $q\in Q$.  Indeed,
  the coordinates of the vertices above can be obtained by averaging
  the coordinates of the vertices in Lemma
  \ref{lem:coords_vertices_fiber}.  This is a particular case of the
  fiber polytope construction of Billera and Sturmfels,
  \cite{BS}.

\demo{Definition {\rm 2.12.6}}
  Choose a system of heights which is an arbitrary real-valued
  function $\psi:X\cap Q\to\bR$, i.e.\ an element of
  $\Fun(C_{\max},\bR)$.  Consider the lower envelope of the
  convex hull of the rays $\{(x,h) \,|\, x  \in X\cap Q,h\ge\psi x\}$
  in $\XR\oplus\bR$. The projections of the faces of this lower
  envelope define a subdivision $\Delta$ of $Q$. The set $C$ is the
  subset of $Q\cap X$ of the points where the corresponding point
  $(x,\psi x)$ lies on the boundary of
  ${\rm ConvHull}_{\psi}$, i.e.\ on the lower envelope. A
  subdivision $(\Delta,C)$ of this type has been called, in
  [GKZ], {\it regular} or
  {\it coherent}.
\enddemo

These definitions have already appeared in our description of the 1-para\-meter
degenerations of stable toric pairs in Section
2.8.

\vglue12pt 2.12.7.
  For a polytope $\sigma_{\alpha}$ denote by $\Fan_{\alpha}$ its dual
  fan in $\Hom(\bL_{\bR},\bR)$. Recall that its 1-dimensional cones
  are the interior normals of codimension faces of $\sigma_{\alpha}$
  and several such rays form a cone if and only if the corresponding
  faces intersect. One of the basic facts about Minkowski sums is that
  the fan $\Fan_{\Sec(Q)}$ is obtained by intersecting all
  $\Fan_{\alpha}$.
  
  Also basic  is how the faces of the Minkowski sum (and
  therefore the cones in the intersection of fans) are obtained. For
  every $\bar h\in \Hom(\bL_{\bR},\bR)$, i.e a restriction of a height
  function $h$, for each $\sigma_{\alpha}$ one picks a face
  $F_{\alpha,\bar \psi}$ on which $h$ takes the minimum. Then the face
  corresponding to $h$  in the Minkowski sum $\Sec(Q)$ is the
  Minkowski sum $\sum F_{\alpha,\bar \psi}$.
  
  This implies that faces of the secondary polytope correspond
  bijectively to the {\it regular} subdivisions of $Q$. In
  particular, the vertices bijectively correspond to the regular
  triangulations. 
\vglue12pt

  Below is our generalization of the above combinatorial
  constructions.  These generalized polytopes describe other
  irreducible components in our moduli space and its simplification.
  
  We now fix an arbitrary subdivision $(\Delta,C)=\{(Q_i,C_i)\}$ of
  $(Q,X\cap Q)$.

\vglue12pt {\it Definition} 2.12.8.
  The generalized secondary polytope $\Sec(\Delta,C)$ is the image in
  $C_0(\uFun_{\bR})/B_0(\ubL_{\bR})$ of the direct product of the
  secondary polytopes $\Sec(Q)_i$ in $C_0(\uFun_{\bR})$. Alternatively, 
  consider the image of the direct product $\prod\sigma_i$ in
  $C_0(\uFun_{\bR})/B_0(\ubL_{\bR})$. This maps surjectively to the
  product of polytopes $Q_i$ in $C_0(\bX_{\bR})$. The generalized
  secondary polytope is the (dilated by $(r+1)$) Minkowski integral of
  the fibers of this polytope map.
\pagebreak

{\it Example} 2.12.9.
  In the easiest situation, assume that all the codimension one faces
  of all $(Q_i,C_i)$ are simplices. Then
  $C_1(\ubL_{\bR})=0$. Consequently, $B_0(\ubL_{\bR})=0$ and 
  $\Sec(\Delta,C)$ is simply the direct product of the secondary
  polytopes $\Sec(Q)_i$. 
\vglue8pt

{\it Definition} 2.12.10.
  A {\it regular} sub-decomposition of $(\Delta,C)$ is defined by the
  system of heights $\psi_i$ in $C_0(\uFun_{\bR})/B_0(\ubL_{\bR})$. On
  each $Q_i$ it is defined as before, and it is easy to see that these
  fit on the intersections. In other words, we consider the system of
  heights $\psi_i$ which on the intersections of $Q_i$'s may differ,
  but only by linear functions (nonhomogeneous on $\XR$, which is the
  same as homogeneous linear on $\bXR$).
\vglue8pt

  We note that by the result of Section 2.8 these are precisely the
  subdivisions of $(\Delta,C)$ that appear for the 1-parameter
  degenerations of a pair of type $(\Delta,C)$. The same argument as
  above about the faces of the Minkowski sum gives the following:
\vfill
{\elevensc Lemma 2.12.11}.  
 {\it  There is a bijective correspondence between the faces of the
 generalized secondary polytope $\Sec(\Delta,C)$ and regular
  subdivisions of}  $(\Delta,C)${\rm .}
\vfill

 {\it Remark} 2.12.12.
  Every regular decomposition of $(Q,X\cap Q)$ which is a sub-decomposition of $(\Delta,C)$ is a regular
decomposition of
  $(\Delta,C)$, because we can take the same heights. The opposite is
  not true. Every regular decomposition of $(\Delta,C)$ defines a
  collection of regular decompositions of $(Q_i,C_i)$. The opposite is 
  not true as well.
\vglue8pt

  2.12.13. 
  Give the vector space
  $C_0(\Delta,\uFun_{\bR})/B_0(\Delta,\ubL_{\bR})$ the integral
  structure by the lattice
  $\left(C_0(\Delta,\uFun)/B_0(\Delta,\ubL)\right)/\hbox{torsion}$. We
  can consider the toric scheme $F_{\Sigma(\Delta,C)}$ corresponding
  to this polytope and the above lattice. Again, there is a natural
  morphism $F_{\Sigma(\Delta,C)}\to M_{\simp}$. This time, however, we 
  can only claim that it is finite to image because for a face $F$ the 
  lattice $X_F$ may contain the corresponding lattice $\bL_i$ for one
  of the polytopes $\sigma_{\alpha}$ as a sublattice of finite
  index. Therefore, $F_{\Sigma(\Delta,C)}$ and an appropriate closed
  subset of $M$ differ ``in finitely many ways''.
\vglue4pt

2.13. {\it The moment map for the moduli space}.
  Consider now the case of the base field $\bC$. For any semigroup $H$
  let $X=\Spec\, \bC[H]$. Then the moment map $\mom: X(\bC)\to \mom\, X$
  is defined to be the quotient map by the equivalence relation: two
  points $p_i:\bC[H]\to \bC$ coincide if $|p_1(h)|=|p_2(h)|$ for all
  $h\in H$. If a variety $Y$ is glued from several varieties $\Spec\, \bC[H_i]$'s and the gluing maps are given by monomials then the
  equivalence relation is preserved and the moment map is globally
  defined. This is the case, for example, for toric varieties, and it
  is well-known that the image of a projective toric variety can be
  naturally identified with a lattice polytope.\pagebreak

 2.13.1.
  Both $M=\TP^{\framed}[Q]$ and $M_{\simp}$ were defined in semigroup
  terms. They are covered by affine schemes corresponding to semigroup
  algebras, and the gluing maps are given by monomials. Therefore, we
  have moment maps defined. Moreover, since all strata differed only
  by a finite group of multiplicative type, i.e.\ $\mu_{n_1}\times
  \cdots \times \mu_{n_k}$, the images of the moment maps of $\mom\, M(\bC)$ and $\mom\, M_{\simp}(\bC)$ are exactly the same.
 
\specialnumber{2.13.2}\proclaim{Lemma}
  $\mom\, M = {\displaystyle\lim_{{\llar}}} \mom\, F_{\alpha} = {\displaystyle\lim_{{\llar}}}
  \sigma_{\alpha}${\rm .} 
\endproclaim

{\it Proof}.
  Indeed, it is clear that ${\displaystyle\lim_{{\llar}}}$ commutes with taking the
  quotient by the above equivalence relation.
\hfill\qed

\specialnumber{2.13.3}\proclaim{Theorem}
  $\mom\, M$ is obtained by gluing   the generalized polytopes $\Sigma(\Delta,C)$
  {\rm (}\/some of these polytopes may appear as faces of others{\rm ,} others may
  be truly new\/{\rm ).}
\endproclaim

{\it Proof}.
  This follows at once from the results of the previous subsection.
\hfill\qed\vglue4pt

2.14.  {\it Relation to families of Gelfand-Kapranov-Zelevinsky}.
\vglue4pt

 2.14.1. Consider a pair $(P,\Theta)$ over a closed field $k$ that
  corresponds to a decomposition $(\Delta,C)$. As before, a section
  $\theta\in H^0(P,L)$ defining the divisor $\Theta$ can be written
  uniquely as the sum of its homogeneous components, i.e.,
  eigenfunctions$$ 
    \theta=\sum_{x\in X\cap Q} \xi_x.
  $$ 
  Therefore, we have $n=\#(Q\cap X)$ sections $\xi_x$ and the
  corresponding rational map $f_{\theta}:P-\to \bP(H^0(P,L)^*)$.

\specialnumber{2.14.2}\proclaim{Lemma}
  $f_{\theta}$ is a finite morphism{\rm .} The image is a cycle with
  multiplicities $\sum d_i V_i${\rm ,} where $d_i$ is the order of the group
  $\bK_i${\rm .}
\endproclaim

{\it Proof}.
  Look at the restriction of $f_{\theta}$ to an irreducible component
  $P_i$. The condition that this restriction is a morphism is
  equivalent to the condition that $\xi_x\ne0$ for all the vertices of 
  the polytope $\delta_i$; this is precisely our condition on the
  divisor $\Theta$. The computation of the multiplicity is standard.
\hfill\qed\vglue4pt

  Fixing an origin in the main orbit of the $\bT$-action on $L_i$
  identifies $P_i$ with the standard projective variety $(P,L)[Q]$,
  and $H^0(P_i,L_i)$ -- with the standard vector space
  $\Fun(\delta_i\cap X, k)=k^{n_i}$.  Now choosing a different origin
  will change the isomorphism $H^0(P_i,L_i)^*\,\isoto\, k^{n_i}$. However,
  it will not change the image $f_{\theta}(P)$ since it is $T$- and
  $\bT$-invariant. Therefore, the cycle $\sum d_i V_i$ does not depend 
  on the choice of origins made.

\vglue4pt 2.14.3.
  In their book [GKZ] and in the
  previous papers referenced there Gelfand-Kapranov-Zelevinsky
  consider a possibly nonnormal toric variety $V$ which is the image
  of the standard toric variety $P[Q]$ as described above, the cycles
  obtained from it by the action of the big torus $(k^*)^n$ and the
  cycles which are obtained as toric degenerations of these. They show
  that the extreme toric degenerations correspond to the Chow forms of
  weight $\phi_{\Tr}$ (with respect to the character group of the big
  torus, i.e.\ $\Fun(X\cap Q,\bZ)$), for all  the regular decompositions
  $\Tr$. Moreover, they show that the parameter space of the cycles
  that appear (as the result of the toric action
  and toric degenerations) is a possibly nonnormal toric variety
  corresponding to the polytope $\Sec(Q)$.
\vglue12pt

2.15. {\it Other families}.  Because of combinatorial restrictions on lattice polytopes $\delta$ --
integral points in some generate the lattice and the semigroup of
integral points in $\Cone\,\delta$, and for some it is not true -- it
seems inevitable that we should either have
\vglue4pt

 1.  nice fibers in the families but a somewhat complicated moduli
  space,~or 

  2.  a nice parameter space but somewhat complicated fibers.
\vglue4pt

In this paper we have chosen the first solution, in part because we can 
compute the cohomologies $H^p(P,L^d)$ and it just seems to be more
reasonable for the objects we are trying to parametrize to be as nice 
as possible. But, taking the second point of view one can construct at 
least two types of families, of schemes and not just the cycles, over
the projective toric scheme $P_{\Sec(Q)}$ corresponding to the
secondary polytope. 

\vglue12pt 2.15.1. (\/{\it Family of the first type}\/).
  Take a regular decomposition $(\Delta,C)$ of $(Q,Q\cap X)$. There is
  a corresponding face $F[\Delta,C]$ of $\Sec(Q)$ and the corresponding
  fiber semigroup $K[\Delta,C]$ consisting of all differences
  $\chi_1-\chi_2\in \bL$ with $\chi_1\in \Cone\Sec(Q)$, $\chi_2\in
  \Cone\, F[\Delta,C]$. Since $\Sec(Q)$ is a Minkowski sum of fibers
  $\sigma_q$, for every $q\in Q$ there is a fiber group $K_{\delta,q}$
  contained in $K[\Delta,C]$. Now consider the semigroup algebra over
  $\bZ[K[\Delta,C]]=\oplus_{k\in K[\Delta,C]} \bZ \zeta^k$ generated
  by variables $\xi_x$, $x\in Q\cap X $   with the following
  relations:$$ 
    \prod \xi_x^{m_x} = \prod \xi_x^{n_x} \times
    \zeta^{(m_x)-(n_x)} 
  $$ 
  whenever $(m_x)-(n_x)\in K[\Delta,C]$. Define the family as the
  $\Proj$ of this algebra. Intersection of two families over
  $K[\Delta_1,C_1]$ and   $K[\Delta_1,C_1]$ is identified with the
  one over $K[\lub((\Delta_1,C_1),(\Delta_2,C_2))]$ so they glue to 
  a global family over the projective scheme.
\vglue12pt

2.15.2. (\/{\it Family of the second type}\/).
  This is defined to be the normalization of the first family.
\vglue12pt

Both of these familes are different from the ones we are considering. The
generic fiber of the first need not be normal, and the flatness of the
second is somewhat unclear. Both of these families may and do have
some nonreduced fibers. We note that in \cite{Nam1}
Y. Namikawa  \pagebreak has constructed families of first and second type over the
toroidal compactification of $A_g^{(2n)}$ for the 2nd Voronoi
decomposition with the principal level structure, $n\ge3$. These are
analogs of the above two families.

\vglue12pt 2.16. {\it First nontrivial example}.
  A well-known (the first, really) example when there are nonregular
  decompositions is the following configuration: $Q$ is a triangle and
  $C$ is a 6-point set as in Figure~1.  There are
  exactly 14 nonregular decompositions of $(Q,C)$. Two of them are
  triangulations as shown in Figure~1 (they
  appear e.g.\ on p.~219 of
  [GKZ]).  We call them ``Escherian
  staircases'': one of them always goes up, and the other always goes
  down.
\begin{center}
\hskip-.5in\BoxedEPSF{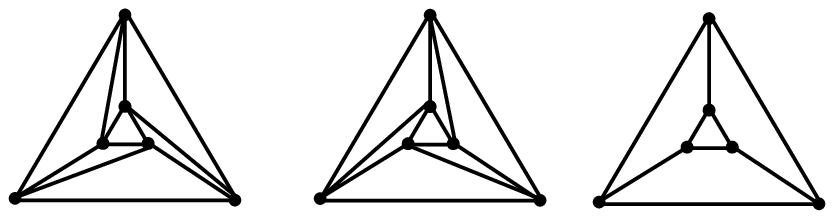 scaled 1000}\hglue.5in
\end{center}
\vglue-24pt
 \centerline{\ninepoint Figure 1. Some decompositions of the 6-point configuration.}
  \vglue12pt

  The secondary polytope $\Sec$ can be computed with some 
  work, and it is combinatorially equivalent to the first polytope in
  Figure~2.  For all regular decompositions
  $\Delta$ the secondary polytopes $\Sec(\Delta)$ are the faces of
  $\Sec(Q,Q\cap X)$ with one exception: when $\Delta$ is the third
  decomposition  in Figure~1.  The interior triangle
  $\delta_0$, being a simplex, has a point as its secondary polytope.
  The secondary polytopes for the quadrangulars
  $\delta_1,\delta_2,\delta_3$ are closed intervals. The space
  $H_0(\Delta,\bL)=\bL_1\oplus \bL_2\oplus \bL_3$ is 3-dimensional
  since $\bL_{\delta}=0$ when $\delta$ is an interval, so that
  $C_1(\ubL)=0$. The secondary polytope $\Sec(\Delta)$ in $\LDel$ is
  the direct product of three intervals, i.e.\ a cube, somewhat
  unrealistically shown in Figure~2.  The
  projection of this cube to $ \bL_{C,\bR}$, however, is
  2-dimensional, and it is the hexagon lying in the foundation of
  $\Sec$. Every cell decomposition of $(Q,C)$ corresponds to a face of
  one or both of these secondary polytopes.

  The moduli interpretation of all this is as follows. The space
  $\MQAred$ has two irreducible components. Both are projective toric
  varieties for the polytopes in Figure 2, the
  second being simply $\bP^1\times\bP^1\times\bP^1$. The first
  irreducible component parametrizes ``smooth'' (in fact only normal)
  and ``smoothable'' stable toric pairs $(P,\Theta)$.
  
  A point $(a,b,c)\in \bP^1\times\bP^1\times\bP^1$ corresponds to a
  smoothable pair if and only if $abc=1$ (or $x_0y_0z_0=x_1y_1z_1$ in projective
  coordinates). The nonregular decompositions are the ones that
  violate this formula. For example, for the two triangulations in
  Figure~1 one has $0\cdot0\cdot0\ne1$ and
  $\infty\cdot\infty\cdot\infty\ne1$ (on the other hand,
  $0\cdot0\cdot\infty=1$, so to say, and thus gives a regular
  triangulation).
\pagebreak

\begin{center}
\BoxedEPSF{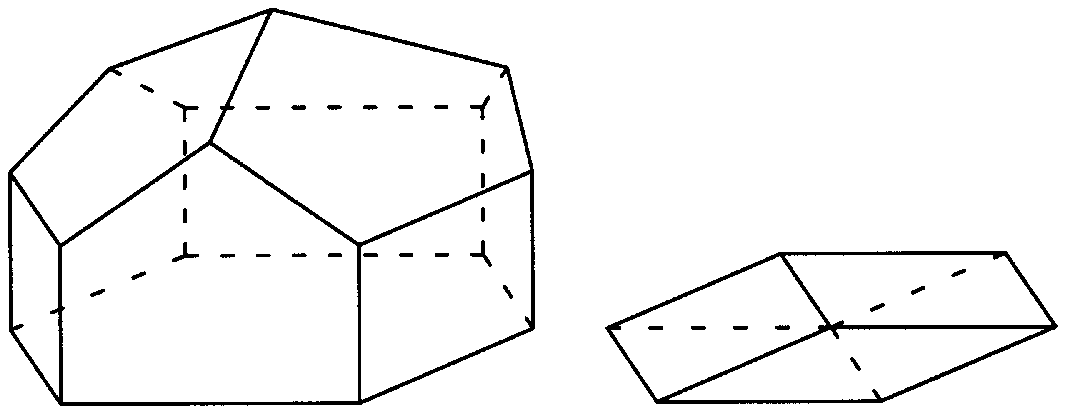 scaled 1100}\end{center}
\vglue12pt
\centerline{\ninepoint Figure 2. Two secondary polytopes of the 6-point example.}
    \vglue12pt

  Note also that the intersection of two components is a blowup
  of $\bP^2$ at three points.

  The image of the moment map $\mom\, M$ in this case is the
  intersection of the above two polytopes, with the hexagonal
  foundation of the first polytope sitting in a ``curved way'' in the
  second.

\vglue12pt 2.16.1.
  Prof. W. J. Whiteley had pointed out to me that if we rotate the
  interior triangle slightly in Figure 1 then the third decomposition
  will no longer be regular. The secondary polytope in this case is
  the deformation of the one above obtained by adding a point under
  the foundation and connecting it to three of the six vertices above it.
  The second polytope is still a cube. However, the intersection of
  the polytopes in this case will consist of three square faces lying
  under the old foundation. We can adopt the above arithmetics to this
  case, too: a point $(a,b,c)$ corresponds to a smoothable pair if and
  only if $abc=\infty$.

\section{Abelic pairs}
\label{sec:Abelic pairs}

\specialnumber{3.0.1}\proclaim{Lemma}\label{lem:polarization_ca_to_a}
  Let $(A\acts P,L)/S$ be a polarized abelic scheme over $S$ of degree
  $d${\rm .} Then it induces a polarization $\lambda:A\to A^t=\bPic^0_{A/S}$
  of degree $d$ on the abelian scheme $A${\rm .}
\endproclaim

\demo{Proof}
  The invertible sheaf $L$ and the action $A\acts P$ define a
  polarization homomorphism $\lambda(L):A\to\bPic^0_{P/S}$; cf.\
  4.1.15.  After an fppf (even \'etale) base
  change $S'\to S$ we have an isomorphism $A'\simeq P'$ so that it is is
  obvious that $\deg\lambda=d^2$. The result now follows from the next
  theorem.
\enddemo

{\it Definition} 3.0.2.
  Let $A$ be an abelian scheme over $S$ and let $B$ and $C$ be two
  schemes over $S$ with an $A$-action. $B\stackrel{A}{\times}C$
  will denote the contracted product, i.e.\ the quotient
  $(B\mathop{\times}\limits_{S}C)/A$ in the fppf topology with the normal
  action of $A$ on $B$ and the opposite to normal action on $C$.  In
  particular, we can take $C=P$ and $B=\bPic_{A/S}$ on which $A$
  naturally acts by translations: $(a,L)\mapsto T^*_aL$.
 
\specialnumber{3.0.3}\proclaim{Theorem}\proclaimtitle{Raynaud \cite{Ray1}, XIII.1.1}
  \label{thm:torsors}
  Let $A$ be an abelian scheme over $S$ and $P$ be an $A$\/{\rm -}\/torsor{\rm .}
  Then
  \begin{itemize}
  \ritem{(i)} There exists a canonical isomorphism$$ 
      \uPic_{A/S}\stackrel{A}{\times}P\simeq \uPic_{P/S};
    $$ 
  \ritem{(ii)} The action of $A$ on $A^t=\uPic^0_{A/S}$ is trivial and it
    induces canonical isomorphisms$$ 
      A^t\simeq\uPic^0_{P/S} \hbox{ and } 
     \uNS_{A/S}\simeq \uNS_{P/S}.
    $$ 
  \end{itemize}
\endproclaim

\demo{Definition {\rm 3.0.4}}
  All abelic pairs over $S$ of dimension $g$ and degree $d$ form
  a category in which we let the arrows be isomorphisms only{\rm .} In this
  way{\rm ,} we get a groupoid $\cAP_{g,d}$ over the category of schemes{\rm .} The
  automorphism classes of abelic pairs over $S$ form a set{\rm .} This
  gives a functor
  $\APgd:{\rm Schemes}\to{\rm Sets}$. 
\enddemo

  By Lemma \ref{lem:polarization_ca_to_a} there exists a natural
  ``forgetful morphism'' from $\APgd$ to the functor $\Agd$ of
  polarized abelian schemes of dimension $g$ and degree $d$, and the
  same for the groupoids.

\demo{Construction {\rm 3.0.5}}
  Consider now the following object. The polarization $\lambda$
  defines a distinguished component $\bPic^{\lambda}_{A/S}$ of the
  scheme $\bPic_{A/S}$ whose geometric points correspond to line
  bundles on $A_{\bar s}$ of numerical type $\lambda$. This component
  is an $A^t$-torsor and is a smooth scheme over $S$. Since $A/S$
  has a section, there is a universal Poincar\'e bundle $\cF_{\lambda}$
  on $A{\displaystyle\mathop{\times}_S}\bPic^{\lambda}_{A/S}$ with a
  rigidification $\cF_{\lambda}|_{0\times\bPic^{\lambda}_{A/S}}\simeq
  \cO_{\bPic^{\lambda}_{A/S}}$. We have a $\bP^{d-1}$-bundle
  $\bP(p_{2*}\cF_{\lambda})$ over $\bPic^{\lambda}_{A/S}$ which
  classifies Cartier divisors on $A$ of numerical type $\lambda$.
  There is a relative Cartier divisor $\Theta_{\lambda}$ on it whose
  geometric points are divisors $D_{\bar s}$ on $A_{\bar s}$ with
  $0\in D_{\bar s}$.  Indeed, this condition is locally given by one
  equation, evaluation at the zero section.
\enddemo

\specialnumber{3.0.6}\proclaim{Lemma}\label{lem:coatprsor_maps_to_picl}
  For every abelic pair $(A,P,\Theta)/S$ there is a natural
  $A$\/{\rm -}\/equivariant morphism $\phi:P\to \bP(p_{2*}\cF_{\lambda})$ such
  that $\phi^*(\Theta_{\lambda})=\Theta${\rm .} When $S=\Spec\, k${\rm ,} $k=\bar
  k${\rm ,} for every $A/k$ there is a $1$\/{\rm -}\/to\/{\rm -}\/$1$ correspondence between
  abelic pairs $(P,\Theta)$ over $k$ and the $A$\/{\rm -}\/orbits of
  $\bP(p_{2*}\cF_{\lambda})${\rm .}
\endproclaim

\demo{Proof}
  Let $p'$ be a functorial point of $P$. By this we mean, as usual, an
  $S'$-point of $P$ for a scheme $S'/S$, i.e., a section of $P'/S'$. It
  defines an isomorphism $A'\to P'$ sending the zero section to $p'$.
  Via this isomorphism we have a relative Cartier divisor $\Theta'$ on
  $A'$ which gives a morphism
  $P'\to\bP(p_{2*}\cF_{\lambda}){\displaystyle\mathop{\times}_S} S'$. It is
  clearly $A$-equivariant.  This gives a morphism of functors and, by
  Yoneda's lemma, a morphism of schemes representing them. Also, the
  points mapping to $\Theta_{\lambda}$ are exactly the points
  contained in $\Theta$ and the preimage of the local equation of
  $\Theta_{\lambda}$ is the local equation of $\Theta$. 
\enddemo

\specialnumber{3.0.7} \proclaim{{C}orollary}\label{cor:Ag=APg}
  In the principally polarized case the groupoids $\cAP_g$ of
  abelic pairs and $\cA_g$ of abelian varieties are naturally
  equivalent{\rm .}
\endproclaim

\demo{Proof}
  Indeed, in this case
  $\bP(p_{2*}\cF_{\lambda})=\bPic^{\lambda}_{A/S}$ is an
  $A=A^t$-torsor, and the above morphism is an isomorphism.
\enddemo

\demo{{R}emark {\rm 3.0.8}}
  Recall that $\cA_g$ is a smooth Mumford\/{\rm -}\/Deligne stack{\rm .}
\enddemo

\specialnumber{3.0.9}\proclaim{{C}orollary}\label{cor:APg_exists}
  There exists a coarse moduli space $\APg$ of principally polarized
  abelic pairs and it is isomorphic to ${\rm A_g}${\rm .}
\endproclaim

3.0.10.
  Let us emphasize again how this equivalence of stacks works. For
  any family $(A\acts P\supset \Theta)/S$ we already have the abelian
  scheme $A/S$, and the principal polarization $\lambda:A\to A^t$ is
  given by Lemma~\ref{lem:polarization_ca_to_a}. Vice versa, given a
  family $(A,\lambda:A\to A^t)/S$ we take $P=\bPic^{\lambda}_{A/S}$
  and by the above construction it comes with a natural relative
  Cartier divisor.

\demo{{R}emark {\rm 3.0.11}}
  A choice of a section $p:S\to P$ gives an isomorphism $\phi_p:A\to
  P$ sending the zero section of $A$ to $p$. Since in general a
  section $p$ need not exist, $A$ and $P$ need not be isomorphic. Even
  when a section does exist (for example, $S=\Spec\, k$, $k=\bar k$), there is
  {\it no canonical isomorphism.\/} On the other hand, note that by
  \cite[6.2]{Mum2} for every $\lambda$ there is a canonical
  section of $\bPic^{2\lambda}_{A/S}$.
\enddemo

\specialnumber{3.0.12}\proclaim{Theorem}\label{thm:APgd}
  The stack $\cAP_{g,d}$ is a separated Artin stack admitting a
  coarse moduli space which is a separated scheme{\rm .} This is a
  Mumford\/{\rm -}\/Deligne stack over $\bZ[1/d]${\rm .}
\endproclaim

\demo{Proof}
  By the above lemma $\cAP_{g,d}$ is the quotient stack
  $[\bP(p_{2*}\cF_{\lambda})/A]$ by the group action, where we
  consider both $\bP(p_{2*}\cF_{\lambda})$ and $A$ to be the stacks
  over $\cA_{g,d}$. Going to an \'etale cover, we obtain $\cAP_{g,d}$
  as the quotient stack of a scheme by an appropriate equivalence
  relation. Since $A$ is proper over $\cA_{g,d}$, this is a proper
  stack, and it also has finite stabilizer. Existence of the coarse
  moduli space as a separated algebraic space in this situation is
  established in [KM]. Outside characteristics dividing
  $d$, the stabilizers are \'etale, and the quotient stack therefore is
  Mumford-Deligne.

  The association $(A,P,\Theta)/S \mapsto (A,\lambda)/S$ defines a
  morphism of stacks $\cAPgd\to \cAgd$ and their coarse moduli
  spaces. Since $\bP(p_{2*}\cF_{\lambda})$ is a projective bundle of
  relative dimension $d-1$, the morphism $\APgd\to \Agd$ is projective 
  of relative dimension $d-1$. Since $\Agd$ is a scheme, so is
  $\APgd$. 
\enddemo

Over $\bZ[1/d]$ the moduli space $\Agd$, and hence also $\APgd$, is the
disjoint union of components corresponding to different types
$(d_1,\dots, d_g)$ with\break $\prod d_i=d$. Recall that this means that on
every geometric fiber, $\ker \lambda:A\to A^t$ together with the
associated symplectic form is isomorphic to $H\times \widehat H$,
$H\simeq \bZ_{d_1}\times \cdots \times \bZ_{d_g}$. We can formally
associate with this component a complex $\Delta$ consisting of one
cell $\oXR/Y$ with $X/Y\simeq H$.
  
\vglue-4pt
\section{Linearization of torus action}
\label{sec:Linearization of torus action}
\vglue-10pt

  Let $H$ be a connected algebraic group acting on a proper variety
  $P$ and let $L$ be an invertible sheaf on $P$. One of standard
  problems of geometric invariant theory is to find a
  linearization of $L$, i.e.\ to extend the $H$-action to $L$. A simple
  answer is that if $H$ does not have any characters, i.e.\ homomorphisms to $\bG_m$, and if $P$ is normal then some power $L^n$
  is uniquely linearizable, see \cite[Ch.\ 1, \S3]{Mum2}.
  
  Here, we are concerned with the opposite situation: $H$ is a torus
  $T$, so that it is essentially built of characters, and $P$ is not
  necessarily normal (but is proper). The answer turns out to be quite
  interesting and nontrivial.  First of all, when it is not empty,
  the set of linearizations is a torsor over the Cartier dual $X$ of
  $T$ (over a closed field, $X\simeq\bZ^r$). Secondly, the
  linearization need not exist at all.  However, we prove that if $P$
  is seminormal then, in a canonical way, there is an infinite \'etale
  $T$-invariant cover $\wP\to P$ on which the pullback $\wL$ of $L$ is
  canonically linearizable. It comes with a properly discontinuous  action of $X$ in
  Zariski topology, and one has $P=\wP/X$, $L=\wL/X$. If
  $L$ is linearizable (e.g.\ $P$ is normal) then $\wP$ is not  new:
  it is simply a disjoint union of many copies of $P$. However, if $L$
  is not linearizable then $\wP$ contains a connected component which
  is not of finite type!
  
  To the author's knowledge, covers of this sort have made their first
  appearance in   Mumford's construction
  \cite{Mum1} where they are presented as {\it ad
    hoc\/} examples. Some schemes $\wP$ are given first, and then
  schemes $P$ are obtained as quotients.  We have arrived at the
  answers of this section largely by trying to understand the exact
  nature of examples in \cite{Mum1}. This
  understanding will enable us to invert   Mumford's construction
  and go from the scheme downstairs to its infinite cover. 
  
  The main result to establish is the equivalence of categories:
\vglue-20pt
\phantom{hi}

\begin{itemize}
\item[1.] Families $(T\acts P,L)$ over some base $S$, with the action by
    a split torus $T/S=\hX/S$, with proper $P/S$ and relatively
    ample $L$, and with seminormal geometric fibers;
 
\item[2.] Families $(T\acts \wP,T\acts \wL)$ with a properly
    discontinuous in Zariski topology action $X\acts \wP$ commuting
    with the $T$-action in a special way, such that
    $(P,L)=(\wP,\wL)/X$. 
  \end{itemize}

4.1. {\it The polarization morphism and theorem of square}.

\vglue4pt\hglue22pt  A. {\it The case of a closed field}.
\vglue4pt

{\it Setup} 4.1.1. 
  Consider 
  \begin{itemize}
  \item[1.] A proper variety $P$ defined over an algebraically closed
    field~$k$. Recall that all our varieties are connected and reduced
    but not necessarily irreducible.
  \item[2.] A connected smooth commutative group variety $G$ over $k$
    acting on $P$ via $\sigma:G\times P \to P$.
  \item[3.] An invertible sheaf $L$ on $P$.
  \end{itemize}

\demo{Definition {\rm 4.1.2}}
  We will say that {\it theorem of square\/} holds if for any
  $a,b\in~G(k)$$$ 
    T^*_{a+b} L \otimes L \simeq T^*_a L \otimes T^*_b L.
  $$ 
\enddemo

{\it Definition} 4.1.3.
  The {\it polarization morphism\/} $\lambda(L):G\to\bPic\, P$ is
  defined by$$ 
    \lambda(L):a\mapsto T^*_{a}L \otimes L^{-1} .   
  $$ 
\vglue12pt

  Obviously, theorem of square holds if and only if the polarization morphism is
  a homomorphism.

  Here $\bPic\, P$ is a group scheme locally of finite type over $k$.
  The connected component of identity $\bPic^0 P$ is of finite type,
  as well as the group scheme $\bPict P$ classifying sheaves $M$ whose
  finite powers are algebraically equivalent to zero. $(\bPic^0
  P)_{\red}$ is a group variety over $k$. In characteristic 0 one of
  course has $\bPic^0 P=(\bPic^0 P)_{\red}$.  Since we assumed $G$ to
  be connected and reduced, $\lambda(L)$ factors through $(\bPic^0
  P)_{\red}$.

  The action $G\acts P$ induces the translation action $G\acts\bPic\, P$. If $a\in G'(S')$, $G'=G{\displaystyle\mathop{\times}_k}S'$ is a functorial
  point of $G$ for a $k$-scheme $S'$ and $M$ is an invertible sheaf on
  $P'=P{\displaystyle\mathop{\times}_k}S'$ then$$ 
    T^*_a M = a^*(\sigma^* M).
  $$ 
  Note that this action preserves the group law on $\bPic\, P$ and that
  $\bPic^0 P$, $\bPict P$ and $(\bPic^0 P)_{\red}$ are all
  $G$-invariant.

\vglue12pt 4.1.4.
  Let us review the situations in which theorem of square is known:
  \begin{itemize}
  \item[1.] $G$ is an abelian variety: because any morphism from an
    abelian variety to a group variety taking zero to zero is a
    homomorphism.  Moreover, this works in families, too:
    \cite[Cor.\ 6.4]{Mum2}.
  \item[2.] $(\bPic^0 P)_{\red}$ is an abelian variety (which is the case
    if $P$ is nonsingular in codimension 1): because any morphism from
    a group variety to an abelian variety taking zero to zero is a
    homomorphism as well; see \cite[II.1,
    Thm.4]{Lan}.
  \item[3.] \cite[III.3]{Lan} also contains the following
    version of theorem of square: Let $D$ be a Weil divisor on an
    irreducible but not necessarily proper $P$. Then for any $a,b$ the
    cycle$$ 
      T^*_{a+b}D - T^*_a D - T^*_b D +D
   $$ 
   is rationally equivalent to zero.
   
   This implies the version we want provided there are no rational
   maps from $\bP^1$ to $(\bPic^0 P)_{\red}$. This is the case, again,
   if $(\bPic^0 P)_{\red}$ is an abelian variety.
   \end{itemize}

  We would like to have theorem of square in the next best case: when
  $G$ is semiabelian. Naively, one would hope that theorem of square
  always holds here, too.

\demo{Counterexample {\rm 4.1.5}}
  Let $P$ be the cuspidal curve $y^2z=x^3$ with the parametrization
  $p_t=(x,y,z)=(t^{-2},t^{-3},1)$, $t\in\bP^1$, together with the
  following torus action: for $a\in\bG_m$, $T_a(p_t)=p_{a^{-1}t}$.
  
  The Picard group $\bPic^0P$ is $\bG_a$ and   is identified with
  $P\setminus\hbox{cusp }\{t=\infty\}$ via $t\mapsto [p_t-p_0]$.  It
  is easy to see that $T^*_a\cO(p_{t_0})=\cO(p_{at_0})$, so that
  $\lambda(\cO(p_{t_0}))(a)= at_0-t_0\in\bG_a$. This is not a
  homomorphism unless $t_0=0$.
\enddemo

  However, we will prove the following:

\specialnumber{4.1.6}\proclaim{Theorem}\label{thm:thmsq_field}
  Let $G$ be a semiabelian variety and $P,L$ be as in
  {\rm 4.1.1.}  Assume that the unipotent rank of
  $\bPic^0P$ is zero{\rm ,} i.e.\ it does not contain a copy of $\bG_a${\rm .}  Then
  theorem of square holds{\rm .}
\endproclaim

\specialnumber{4.1.7}\proclaim{Theorem}\label{thm:seminormal_Picard}
  Let $P$ be a proper seminormal variety{\rm .} Then the unipotent rank of
  $\bPic^0P$ is zero\/{\rm ;} moreover{\rm ,} any morphism $\bA^1\to \bPic^0P$ is
  constant{\rm .}
\endproclaim

\specialnumber{4.1.8}\proclaim{Lemma}
  Let $R$ be a reduced ring {\rm (}\/commutative with {\rm 1,} of course\/{\rm )} and $T$ be
  a free indeterminate{\rm .} Then the subring of invertible elements
  $(R[T])^*$ coincides with $R^*${\rm .}
\endproclaim

\demo{Proof}
  This is obvious for an integral ring, and hence for direct sums of
  integral rings. Since any reduced ring can be embedded into a direct
  sum of integral rings, we are done.
\enddemo

\specialnumber{4.1.9}\proclaim{{C}orollary}  
  For any reduced scheme $X$ over $k$ one has
  $p_{1*}\cO^*_{X\times\bA^1}\break =\cO^*_X${\rm .}
\endproclaim

\specialnumber{4.1.10}\proclaim{Lemma}
  Let $X$ be a reduced scheme over $k$ and $F$ be an invertible sheaf
  on $X\times\bA^1${\rm .} Assume that there exists an open cover $\{U_i\}$
  such that the restriction of $F$ to each $U_i\times\bA^1$ is the
  pullback of an invertible sheaf on $U_i${\rm .} Then $F$ is the pullback
  of an invertible sheaf on $X${\rm .}
\endproclaim

\demo{Proof}
  By the Leray spectral sequence we have the following sequence of
  cohomologies computed in Zariski topology:$$ 
    0\ \to H^1(p_{1*}\cO^*_{X\times\bA^1}) \to 
    H^1(\cO^*_{X\times\bA^1}) \stackrel{\phi}{\to }
    H^0(R^1p_{1*}\cO^*_{X\times\bA^1}).
  $$ 
  The condition of the lemma says exactly that the image of the class
  $[F]\in\Pic(X\times\bA^1)$ under $\phi$ is zero. On the other hand,
  by the previous corollary the kernel is
  $H^1(p_{1*}\cO^*_{X\times\bA^1})=H^1(\cO^*_X)=\Pic\, X$.
\enddemo

\demo{Proof of Theorem {\rm \ref{thm:seminormal_Picard}}}
  A morphism $\bA^1\to\bPic\, P$ induces a sheaf $F$ on $P\times\bA^1$,
  the pullback of the universal bundle on $P\times \bPic\, P$. 
  The morphism is constant if and only if $F$ is a pullback from~$P$. By the
  previous lemma we only need to know this locally.  The local
  statement is proved by Traverso, \cite[3.6]{Tra}.
\enddemo

\specialnumber{4.1.11}\proclaim{Lemma}\label{lem:trivaction_to_thmsq}
  If the translation action of $G$ on $(\bPic^0P)_{\red}$ is trivial
  then theorem of square holds{\rm .}
\endproclaim

\demo{Proof}
  Indeed,
  \begin{eqnarray*}
     T^*_{a+b}L\otimes (T^*_{a}L)^{-1}\otimes 
    (T^*_{b}L)^{-1} \otimes L 
    &=& T^*_{a} (T^*_{b}L)\otimes (T^*_{a}L)^{-1}\otimes 
    (T^*_{b}L\otimes L^{-1})^{-1}  \\
    &=& T^*_{a} (T^*_{b}L\otimes L^{-1})\otimes 
    (T^*_{b}L\otimes L^{-1})^{-1}. 
  \end{eqnarray*}
  Since we assumed in 4.1.1 that $G$ is connected,
  $T^*_{b}L\otimes L^{-1}$ gives a point of $(\bPic^0P)_{\red}$.
\enddemo

\specialnumber{4.1.12}\proclaim{Lemma}\label{lem:triv_action_Picz_field}
  Assume that $G$ is divisible{\rm ,} i.e.\ for every integer $n$ the
  homomorphism $n:G\to G$ of multiplication by $n$ is surjective{\rm .}  If
  the unipotent rank of $(\bPic^0P)_{\red}$ is zero then the
  translation action of $G$ on $(\bPic^0P)_{\red}$ is trivial{\rm .}
\endproclaim

\demo{Proof}
  It is well-known that $H=(\bPic^0P)_{\red}$ in this case is a
  semiabelian variety itself. Indeed, every commutative group variety
  over $k=\bar k$ is an extension of an abelian variety by an affine
  group, which in turn is obtained by extending a torus by a unipotent
  group, which in turn has a filtration with every quotient isomorphic
  to $\bG_a$.
  
  For an integer $m$ denote $_mH=\ker(m:H\to H)$.  Since the
  $G$-action preserves the group law, each $_mH$ is $G$-invariant.
  The action $G\acts \,_mH$ is trivial. Indeed, it is equivalent to a
  homomorphism $G\to\Aut(_mH)$ and the latter is a finite group scheme
  over $k$. If we take $m$ to be coprime to $\chr k$ then $_mH$ and so
  also $\Aut(_mH)$ are reduced. Every point of it is annihilated by
  some $n$. On the other hand $n:G\to G$ is surjective. Hence
  $G\to\Aut(_mH)$ must be trivial. This argument works as well for any
  $m$, with geometric points replaced by functorial points.
  
  Since the family $\{{_m}H\}$ is dense in $H$ (where we can only take
  $m$ that are powers of an arbitrary integer coprime to $\chr k$)
  the action $G\acts H$ is trivial.
\enddemo

{\it Proof  of {\rm 4.1.6}}.
  Since a semiabelian group variety is divisible, we are done by
  \ref{lem:triv_action_Picz_field}, \ref{lem:trivaction_to_thmsq}.
\hfill\qed\vglue4pt

  For the future, we will need the following:

\specialnumber{4.1.13}\proclaim{Lemma}\label{lem:action_Picz_to_Pict}
  Assume that $G$ is divisible{\rm .} Then the action
  $G\acts(\bPic^0P)_{\red}$ is trivial if and only if the action $G\acts\bPict P$
  is trivial{\rm .}
\endproclaim

{\it Proof}.
  Indeed, $\bPict P/(\bPic^0P)_{\red}$ is finite. By the same argument
  as in the previous lemma the action of the divisible group on it
  must be trivial.
\hfill\qed\vglue6pt

B. {\it The case of families}.
  Our next step is prove the 
  result for families  corresponding to 4.1.6. For this, we have to deal with questions about
  the relative Picard functor and multiplicative tori over arbitrary,
  locally noetherian schemes. When discussing these topics, one
  usually works with one of the Grothendieck topologies that is more
  flexible than the Zariski topology. The standard choice is the fppf
  (faithfully flat of finite presentation) topology.  Indeed, by
  \cite[X]{SGA3} every torus becomes split on an fppf cover, and the
  relative Picard functor is defined by sheafifying the ordinary
  Picard functor in the fppf topology.
  
  Hence, for the rest of this section when we say {\it locally}, it
  means {\it locally in the {\rm fppf} topology,\/} and the cohomologies are
  those with respect to  fppf topology. Note that for the statements above  fppf can be replaced by \'etale.  For
smooth, locally of
  finite type, group schemes the cohomologies computed in the fppf and
  \'etale topologies coincide as well:
  \cite[IV.11.7]{DE}.

\vglue6pt {\it Setup {\rm 4.1.14}}.
 1. $P\to S$ is a proper flat morphism with geometrically reduced,
    connected {\it seminormal\/} fibers.
\begin{itemize}
  \item[2.] $G$ is a smooth commutative group scheme over $S$ with
    connected fibers acting on $P$ via the morphism
    $\sigma:G{\displaystyle\mathop{\times}_S}P\to P$.
  \item[3.] $L$ is an invertible sheaf on $P$.
  \end{itemize}

  $\bPic_{P/S}$ will denote the relative Picard functor.  If
  irreducible components of fibers of $P\to S$ are geometrically
  irreducible then by Mumford's theorem $\bPic_{P/S}$ is represented
  by a not necessarily separated scheme locally of finite presentation
  over $S$.  On the other hand, if we instead assume that $P$ and $S$
  are of finite type over a field or an excellent Dedekind domain
  then, according to M. Artin \cite[7.3]{Art1}, the
  relative Picard functor is represented by an algebraic space,
  locally of finite presentation over $S$. We refer to
  \cite[Ch.\ 8]{BLR} or \cite{FGA} for more
  details about this.

  The group scheme $G$ acts on $\bPic_{P/S}$ by translations: if $S'$
  is an $S$-scheme, and $a$ and $M$ are local sections of
  $G'=G{\displaystyle\mathop{\times}_S}S'$ and
  $\bPic_{P'/S'}=\bPic_{P/S}{\displaystyle\mathop{\times}_S}S'$ respectively, then$$ 
    a.M=  T^*_a M =  a^*(\sigma^* M).
  $$ 
  This action preserves the group structure.

\demo{Definition {\rm 4.1.15}} 
  Let $L$ be an invertible sheaf on $P$. The {\it polarization
    morphism\/} $\lambda(L):G\to\bPic_{P/S}$ is defined by$$ 
    \lambda(L):a\mapsto T^*_{a}L \otimes L^{-1}    
  $$ 
  for a local section $a$ of $G'$.
\enddemo

  Assume that we are in the one of the two basic situations above, so
  that $\bPic_{P/S}$ is a scheme (resp.\ algebraic space) locally of
  finite type over $S$.  If $G$ is of finite type over $S$ and has
  connected fibers then $\lambda(L)$ factors through $\bPict_{P/S}$
  which is of finite type over $S$, and the above-defined action
  preserves $\bPict_{P/S}$.  Recall that the latter is the open
  subscheme (resp.\ algebraic space) of $\bPic_{P/S}$ whose geometric
  points correspond to sheaves $M_{\bar s}$ on $P_{\bar s}$ with
  $M_{\bar s}^n\in \Pic^0 P_{\bar s}$ for some $n\in\bN$.

\demo{Definition {\rm 4.1.16}} 
  We say that {\it theorem of square\/} holds for $L$ if $\lambda(L)$
  is a homomorphism. Equivalently, {\it locally\/} on $S$ we must have$$ 
    T^*_{a+b} L \otimes L \simeq T^*_a L \otimes T^*_b L.
  $$ 
\enddemo

{\it Remark} 4.1.17.
  This condition is weaker than asking the above two sheaves to be
  isomorphic or differ by an invertible sheaf on $S$.
\vglue12pt

  As   remarked in 4.1.4, if $G$ is an
  abelian scheme over $S$ then this more general form of theorem of
  square holds as an application of the rigidity lemma; see
  \cite[Cor.\ 6.4]{Mum2}.

\specialnumber{4.1.18}\proclaim{Theorem}\label{thm:thm_square} 
  Let $G$ be a semiabelian group scheme over $S${\rm ,} and $P,L$ be as in
  {\rm 4.1.14.}  Then the theorem of square holds for
  any invertible sheaf $L$ on~$P${\rm .}
\endproclaim

{\it Remark} 4.1.19.
  We use the assumption that every geometric fiber $P_{\bar s}$ is
  seminormal only to conclude that the unipotent rank of $\bPic^0
  P_{\bar s}$ is zero.

\demo{Proof}
  By the same exact argument as in \ref{lem:trivaction_to_thmsq},
  substituting geometric points by functorial points, we see that it
  suffices to prove that the action $G\acts\bPict_{P/S}$ is trivial.
  We already know this statement for an algebraically closed field
  $\bar k$: \ref{lem:triv_action_Picz_field},
  \ref{lem:action_Picz_to_Pict}. Since the question is local, this
  implies the statement for a nonclosed field $k$, too, so we have it
  for every closed fiber.
  
  The main point is to prove the statement for a local artinian base
  $S$. It then gives the general case by a standard argument: Artinian
  case implies the case of a complete local ring, and since the base
  change from a local noetherian ring to its completion is faithfully
  flat, this implies the statement over the local ring $\cO_{S,s}$,
  which clearly suffices.

  Another reduction step is that, since the question is local, we can
  change to an fppf cover on which all irreducible fibers of $P_s$ are
  geometrically irreducible. Hence, we can assume that $\bPict_{P/S}$
  is a scheme.
  
  Now, the artinian case follows from \ref{lem:polmap_artinian}.
\enddemo

\specialnumber{4.1.20}\proclaim{Lemma}\label{lem:polmap_artinian}
  Let $S$ be the spectrum of a local artinian ring with a closed point
  $S_0${\rm .} Let $G$ be a semiabelian group scheme over $S$ and let $Q$ be
  a scheme of finite type over $S$ on which $G$ acts via the morphism
  $G{\displaystyle\mathop{\times}_S}Q\to Q${\rm .} If the action $G_0\acts Q_0$ is
  trivial then the action $G\acts Q$ is trivial{\rm .}
\endproclaim

\demo{Proof}
  $G$ is an extension of an abelian scheme $A/S$ by a torus $T/S$. For
  the torus action we quote \cite[IX.3.8]{SGA3} where it is proved for
  any group scheme of multiplicative type.
  
  Therefore, we are reduced to the case of abelian action $A\acts Q$.
  Let $s$ be a section of $Q$ over $S$. By restricting the action to
  $s$ we have a morphism $A{\displaystyle\mathop{\times}_S}s \to Q$. It is trivial
  on the central fiber, i.e.\ factors through
  $A{\displaystyle\mathop{\times}_S}s\simeq A\to S$ and $s:S\to Q$. By the
 rigidity lemma \cite[Prop.6.1]{Mum1} there is the same
  factorization globally.  Since we also have this for any section
  after any local artinian base change $S'/S$, it follows that the
  morphism $G{\displaystyle\mathop{\times}_S}Q\to Q$ is the projection to the
  second component.
\enddemo

  Next, we would like to analyze the kernel of the polarization
  morphism.

\specialnumber{4.1.21}\proclaim{Lemma}\label{lem:polmap_DVR}
  Assume that $S$ is the spectrum of a {\rm DVR} with the closed point $0$
  and the generic point $\eta${\rm ,} and $G,P,L$ are as in
  {\rm 4.1.14.}  Let $H_{\eta}$ be a subgroup of
  $G_{\eta}$ and let $H$ be its closure in $G${\rm .}  Denote by $\lambda_H$
  the restriction of $\lambda$ to $H${\rm .}  Then $\lambda_H(L_{\eta})=0$
  implies $\lambda_H(L)=0${\rm .}
\endproclaim

\demo{Proof}
  Let $a\in H(S)$.  The sheaf $M$ induces a morphism
  $\lambda(a):S\to\bPict_{P/S}$ and we want to show that it factors
  through the zero section. Combination of  this with all base changes
  $S'/S$ which can be assumed to be again DVRs,   would prove the
  statement.
  
  The image of $\lambda(a)$ is contained in the closure $S_1$ of the
  zero section $S_0$. Note that $S_{1,\eta}=S_{0,\eta}$ and that $S_1$
  is reduced.  We claim that the fiber of $S_1\to S$ over $0$ is a
  finite set. Indeed, for any affine open set $U\subset
  \bPic^{\tau}_{P/S}$ the closure of $S_{0,\eta}$ in $U$ has no more
  than one point over $0$ since $U$ is separated, and
  $\bPic^{\tau}_{P/S}$ is covered by finitely many such $U$'s.
  Because $G_0$ is connected, $\lambda(a)(0)=0$ lies in $S_0$.
  Therefore, $\lambda(a)$ factors through $S_0$.
\enddemo

When we say a {\it group scheme of multiplicative type} we mean that
locally in the fppf topology it is isomorphic to a product of several
copies of $\bG_m$ and $\mu_n$ and we do not necessarily assume that
$n$ is coprime to the characteristics of the residue fields.

\specialnumber{4.1.22}\proclaim{Theorem}\label{thm:imgker_polm_torus}
  Let $G$ be a torus $T$ over $S${\rm ,} and $P,L$ be as in
  {\rm 4.1.14.}  Then $\ker\lambda(L)$ is a flat and
  closed sub group scheme of $T$ of multiplicative type{\rm .} The quotient
  $T/\ker\lambda(L)$ is a torus{\rm ,} and the polarization morphism factors
  through an immersion $T/\ker\lambda(L)\to \bPict_{P/S}${\rm .}
\endproclaim

{\it Proof}.
  Since the question is local, we can assume that $S$ is connected,
  $\bPict_{P/S}$ is a scheme and that $T$ is diagonalizable:
  $T=\bHom_{S\hbox{-gps}}(\bZ^r_S,\bG_{m,S})$. For every $s\in S$ the
  subgroup $\ker\lambda_s\subset T_s$ is closed and is also
  diagonalizable, i.e.\ of the form $\bHom(M_s,\bG_{m,s})$, where $M_s$
  is a quotient of $\bZ^r$.  Denote by $H_s$ the corresponding
  diagonalizable subgroup of $T$.
  
  Since we already know that $\lambda(L)$ is a homomorphism, by
  \cite[IX.5.2]{SGA3} there is an open subset $U\ni s$ of $S$ such
  that $H_s\subset\ker \lambda_U$.  This implies that $\ker\lambda$ is
  flat.  Applying \ref{lem:polmap_DVR} to a base change by a DVR we
  see that it is also closed (if $\bPict_{P/S}$ were separated that
  would be trivial). Hence, $\ker\lambda(L)=H$ for a diagonalizable
  subgroup of $T$.
  
  The quotient $T/H$ is also diagonalizable, and so is a torus. The
  monomorphism $T/H\to \bPict_{P/S}$ is an immersion by
  \cite[XVI.1.4]{SGA3}.
\hfill\qed

\specialnumber{4.1.23}\proclaim{{C}orollary}\label{cor:polmapzero}
  Under the conditions of the previous theorem{\rm ,} if $S$ is connected
  and for one geometric fiber $\lambda(L_{\bar s})=0$ then
  $\lambda(L)=0$ globally{\rm .}
\endproclaim

\vglue-6pt
4.2. {\it Tori in Picard groups and infinite covers}.
 
\vglue2pt {\it Example {\rm 4.2.1}}.
  Let $P$ be a projective rational curve with a single node defined
  over an algebraically closed field $k$. Let $L$ be an ample
  invertible sheaf on $P$, say of degree 1. There is a standard action
  on $P$ of the multiplicative group $\bG_{m,k}$ with two orbits: the
  node and its complement.  It is easy to see that the polarization
  morphism $\lambda(L)$ in this case is a group isomorphism from $G$
  to $\bPic^0P$.
  
  Consider the scheme $\wP$ which is the infinite chain of copies
  $C^{\nu}_i\simeq\bP^1$, $i\in\bZ$, with $C_i$ intersecting $C_{i-1}$
  and $C_{i+1}$, as on the picture. The group $X=\bZ$ naturally acts
  on $\wP$ sending $C_i$ to $C_{i+x}$ and the singular points  to
  singular points. It is easy to see that $P$ is the quotient $\wP/X$.

\begin{center}
\BoxedEPSF{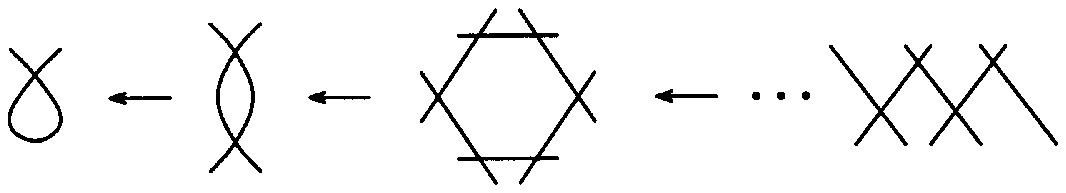 scaled 1000}
\end{center}
\vglue12pt

  Quite interestingly, $\wP$ is not of finite type. However, it is
  locally of finite type and the morphism $\wP\to P$ is \'etale and
  Galois. Another series of examples is provided by   Mumford's
  construction \cite{Mum1}.  We would like to
  understand the nature of such infinite covers.

  The following statement and its corollary are well-known; see   \cite[6.2]{Ray2} (also
  \cite[III.4.16]{Mil}). Sheaves and cohomologies are in the fppf
  topology:

\specialnumber{4.2.2}\proclaim{Theorem}\label{thm:finite_covers}
  Let $f:P\to S$ be a proper morphism with geometrically reduced and
  connected fibers{\rm ,} and let $G$ be a finite flat commutative group
  scheme over $S${\rm .} Then there is a canonical isomorphism
$$ 
    R^1f_* \hG_P = \uHom_{S{\rm -gps}}(G,\bPic_{P/S}),
  $$ 
  where $\hG=\bHom_{S{\rm -gps}}(G,\bG_{m,S})$ is the Cartier
  dual of $G${\rm ,} and $\hG_P=\hG{\displaystyle\mathop{\times}_S}P${\rm .}
\endproclaim

\specialnumber{4.2.3}\proclaim{{C}orollary}\label{cor:finite_covers}
  Let $P$ be a proper connected reduced scheme over an algebraically
  closed field $k${\rm .}  Then there is a $1$\/{\rm -}\/to\/{\rm -}\/$1$ correspondence between
  \begin{itemize}
  \ritem{(i)} finite commutative subgroup schemes $G$ of $\bPic\, P${\rm ,}
  \ritem{(ii)} connected $\hG$-torsors $\wP$ over $P${\rm ,} i.e.\ connected schemes
    $\wP$ over $S$ with a free $\hG$-action and with $\wP/\hG=P$.
  \end{itemize}

\endproclaim

\demo{Proof}
  Fppf sheaves which are $\hG$-torsors are classified by $H^1(P,\hG)$;\break
  cf.\ \cite[III.4]{Mil}. On the other hand, $\hG$ is affine, and
  affine sheaf torsors are always representable:
  \cite[VIII.7.9]{SGA1}.
\enddemo

  The example 4.2.1 is an application of
  \ref{cor:finite_covers}, but with $G$ a multiplicative torus instead
  of being finite. Checking the proof of the above theorem, we see
  that it goes through for this case as well, with obvious changes.

\specialnumber{4.2.4}\proclaim{Theorem}\label{thm:infinite_covers}
  Let $f:P\to S$ be a proper morphism with geometrically reduced and
  connected fibers{\rm ,} and let $T$ be a torus over $S${\rm ,} a twist of
  $\bG_{m,S}^r${\rm ,} and let
  $X=\hT=\bHom_{S{\rm -gps}}(T,\bG_{m,S})$ be its Cartier
  dual{\rm ,} a twist of the constant group $\bZ_{S}^r${\rm .} Then there is a
  canonical isomorphism$$ 
    R^1f_* X_P = \uHom_{S{\rm -gps}}(T,\bPic_{P/S}).
  $$ 
\endproclaim

\demo{Proof}
  The same proof that is given in
  \cite[6.2]{Ray2} for the finite case (it is
  repeated in \cite[4.16]{Mil} as well) works here verbatim.
  Aside from the Cartier duality, the only fact specific to the group
  $T$ is that $\uExt^1(T,\bG_{m,S})=~0$.  This, in fact, is well-known.
  Let us recall the argument. We have to show that for any faithfully
  flat scheme $S'$ over $S$ an extension $H$ of $\bG_{m,S'}$ by
  $T_{S'}$ is locally split. But since both are of multiplicative
  type, by \cite[X]{SGA3} $H$ is also of multiplicative type, and
  therefore all three are locally diagonalizable. By taking the
  Cartier dual we obtain an exact sequence of constant groups, and
  that is obviously split on every connected component of the base
  scheme.
\enddemo

  As for the corollary, we do not even have to assume that $k$ is closed:

\specialnumber{4.2.5}\proclaim{{C}orollary}\label{cor:infinite_covers}
  Let $P$ be a proper connected reduced scheme over a field $k$ which
  has a $k$\/{\rm -}\/rational point{\rm .} Then there is a $1$\/{\rm -}\/to\/{\rm -}\/$1$ correspondence
  between
  \begin{itemize}
  \ritem{(i)} subgroups $T\subset\Pic^0_{P/k}$ which are tori{\rm ,}
  \ritem{(ii)} connected $X$-torsors $\wP$ over $P${\rm ,} i.e.\ connected schemes
    $\wP$ over $S$ with a free $\hT$\/{\rm -}\/action and with $\wP/\hT=P${\rm .}
  \end{itemize}

\endproclaim

\demo{Proof}
  An fppf sheaf which is an $X$-torsor over $P$ is again automatically
  representable. Indeed, this is proved in \cite[X.5.4]{SGA3} for a
  torsor with any twisted constant structure group. The sheaf torsors
  are classified by the first cohomology, as before.  So, we have to
  prove the isomorphism $H^1(P,X) \simeq \Hom(T,\bPic_{P/k})$.
  
  In the case of an algebraically closed field we are already done by
  \ref{thm:infinite_covers}.  Over an arbitrary field $k$ we have the
  following exact sequence:
$$ 
    H^1(k,X)\to H^1(P,X) \to H^0(k,R^1f_*X) \to
    H^2(k,X)\to H^2(P,X),
  $$ 
  and by \ref{thm:infinite_covers},
  $H^0(k,R^1f_*X)=\Hom(T,\bPic_{P/k})$.  Hence, to prove the theorem
  it suffices to show that:
  \begin{itemize}
  \item[(i)] $H^1(k,X)=0$.
  \item[(ii)] The obstruction in $H^2(k,X)\to H^2(P,X)$ vanishes.
  \end{itemize}
  For the first assertion we have$$ 
    H^1(k,X)=\Hom(\Gal(\bar k/k),X)=0
  $$ 
  because $X$ is torsion-free and $\Gal(\bar k/k)$ is profinite.
  
  As for the second statement, by a standard argument (cf.\ proofs of \cite[6.2]{Ray2} or
  \cite[4.16]{Mil}) obstructions of this type always vanish when
  there is a section, which in our case means a $k$-rational point.
\enddemo

{\it Remark} 4.2.6.
  Note that neither object in \ref{cor:infinite_covers} exists when
  $P_{\bar k}$ is normal.  Indeed, in this case every connected
  \'etale cover of $P$ is quasifinite (obviously since it is
  irreducible), and $(\Pic^0 P)_{\red}$ is an abelian variety (see
 \cite[6.2.1]{FGA}).

\demo{{R}emark {\rm 4.2.7}}
  Existence of infinite connected \'etale $\bZ^r$-Galois covers of $P$
  is closely related to existence of twists of the constant group
  $\bZ^r$ over $P$ and the $r$-dimensional tori over $P$ which do not
  become trivial after any {\it finite\/} \'etale cover (by
  \cite[X.5.9]{SGA3} the latter two categories are anti-equivalent).
  Two examples in \cite[X.1.6]{SGA3} of a situation when this occurs
  are exact  for the rational nodal curve as above and also for its
  \'etale double cover.
\enddemo

{\it Example} 4.2.8.
  To generalize 4.2.1, let $C=\cup_{i=1}^N C_i$ be a
  proper reduced connected curve with nodes only, for simplicity over
  an algebraically closed field $ k$. It is well-known that $\Pic^0C$
  is a semiabelian variety which fits into the following exact
  sequence$$ 
    0\to T\to\Pic^0C\to\Pic^0C^{\nu}\to 0,
  $$ 
  where $C^{\nu}$ is the normalization of $C$ and $T=( k^*)^r$. By
  assigning to each irreducible component $C_i$ a vertex and to each
  node an edge we obtain the dual graph $\Gamma$ of $P$. Then the
  dimension $r$ of the toric part is the cyclomatic number of
  $\Gamma$, the rank of $H_1(\Gamma)$ considered as a topological
  space, i.e., the number of ``independent loops''.
  
  Consider the universal abelian cover $\wGamma$ of $\Gamma$
  corresponding to$$ 
    \pi_1^{ab}=\pi_1/[\pi_1,\pi_1]=\bZ^n=\hT=X,
  $$ 
  where $\pi_1=\pi_1(\Gamma)$.  This is an infinite graph with
  vertices marked $C^{\nu}_1$ through $C^{\nu}_N$.  In an obvious way
  one constructs a scheme $\wC$ corresponding to it  which is a
  connected scheme locally of finite type with a morphism to $C$, so
  that $\wC\to C$ is \'etale Galois with the group $X$.

\demo{{\rm 4.3.} Infinite covers arising from the torus action}
Let $S$ be a locally noetherian scheme.
The aim of this section is to prove the following:
\specialnumber{4.3.1}\proclaim{Theorem}\label{thm:torus_action_equiv}
  The following categories are naturally equivalent\/{\rm :}
  \begin{itemize} 
  \ritem{(1)} triples $(T,P,L)/S${\rm ,} where
    \begin{itemize}
    \ritem{(a)} $T$ is a split torus over $S$ with a Cartier dual
      $X_S=\bZ_S^r${\rm .}      
    \ritem{(b)} $P$ is a proper scheme over $S$ with geometrically
      connected{\rm ,} reduced and seminormal fibers{\rm ,} with a $T$\/{\rm -}\/action{\rm .}
    \ritem{(c)} $L$ is a relatively ample invertible sheaf on $P$.
    \end{itemize}
  \ritem{(2)}   triples $(T,\wP,\wL)/S${\rm ,} where
    \begin{itemize}
    \ritem{(a)} $T$ is a split torus over $S$ with Cartier dual $X${\rm .}     
    \ritem{(b)} $\wP$ is a scheme{\rm ,} locally of finite type{\rm ,} over $S$ with
      geometrically reduced and seminormal fibers{\rm ,} with commuting
      actions of $T$ and~$X${\rm .}
    \ritem{(c)} $\wL$ is a relatively ample invertible sheaf on $P${\rm ,} with
      the actions of $T$ and $X$ extending those on $\wP${\rm .}
    \end{itemize}
    such that 
    \begin{itemize}
    \ritem{(a)} $X$\/{\rm -}\/action is properly discontinuous in the Zariski
      topology{\rm ,} and $P=\wP/X$ is proper over $S$ with geometrically
      connected fibers{\rm .}
    \ritem{(b)} The actions of $T$ and $X$ on $\wL$ commute in the following
      way\/{\rm :} if $a$ and $x$ are functorial points of $T$ and $X$ {\rm (}\/i.e.\ sections after a base change $S'/S${\rm ),}
then
$$ 
        T^*_a S^*_x = x(a) \cdot S^*_x T^*_a.
      $$ 
    \end{itemize}
  \end{itemize}

\endproclaim

{\it Definition} 4.3.2.
  A sheaf $\wL$ on a scheme $\wP$ which is not necessarily of finite
  type over $S$ is called {\it relatively ample\/} if there exists an
  open cover $\{S_i\subset S\}$ such that for each $i$ the sections in
  $H^0(\wP_i,\wL_i^d)$, $d\gg0$, form a basis of topology on $\wP_i$.

\demo{Example {\rm 4.3.3}} 
  Let $P$ be a projective $r$-dimensional toric variety over a closed
  field $k$. The polarization morphism $\lambda(L):T\to \Pic\, P$ is
  constant: $\lambda(L)=0$.  An ample line bundle $L$ on $P$ can
  be linearized in infinitely many ways, each equivalent to giving a
  lattice polytope $Q$ in $X_{\bR}\simeq\bR^r$.  The difference
  between any two such polytopes is a shift by a lattice vector $x\in
  X$: $Q_2=Q_1+x$.  Speaking invariantly, we have $\wP$ as in the
  above theorem which is simply a disjoint union of $X$ copies of $P$
  on which the pullback $\wL$ is canonically linearized.
\enddemo

  Of course, this example is rather trivial, and the situation is a
  lot more interesting when the cover $\wP$ is not of finite type as
  in Example~4.2.1. However, we would like to show
  that these examples are both special cases of a general phenomenon
  provided by the above theorem. We would also like to understand the
  number of the connected components of $\wP$.

\demo{Proof of Theorem {\rm \ref{thm:torus_action_equiv}}}
  Given the pair $(\wP,\wL)$ we define $P$ and $L$ as $\wP/X$ and
  $\wL/X$. Because the action is properly discontinuous in Zariski
  topology, for a subgroup $X'\subset X$ of finite index there exists
  a system of neighborhoods $\{U_i\}$ of $\wP$ such that for any
  $x'\ne0$ in $X'$ one has $U_i\cap T_{x'}U_i=\emptyset$. So, the
  quotient $\wP/X'$ is defined as a scheme, and $P$ is defined at
  least as an algebraic space. But it comes with an ample sheaf $L$,
  and so $P$ is actually a scheme.  The opposite direction, going from
  $(P,L)$ to the cover $(\wP,\wL)$ is more complicated.
  
  By \ref{thm:imgker_polm_torus} we have a subtorus
  $T_Y=\im\lambda(L)(T)$ in $\bPict_{P/S}$. Let $Y$ be its Cartier
  dual. Dually, we have a closed subscheme $Y_S\subset X_S$. Theorem
  \ref{thm:infinite_covers} provides that {\it locally} the subtorus
  $T_Y$ gives a $Y$-torsor $P^{\dag}$ over $P$. In other words, there
  exists an fppf cover $\{S_i\to S\}$ such that the subtori
  $T_{Y,i}=T_Y{\displaystyle\mathop{\times}_S}S_i\subset \bPict_{P_i/S_i}$ induce
  the $Y$-torsors $P^{\dag}_i$ over $P_i=P{\displaystyle\mathop{\times}_S}S_i$.
  First, we are going to show that the action of $T$ lifts to each
  $P^{\dag}_i$ and that each $L^{\dag}_i$ is linearizable.
  
  Let $G^{\dag}_i$ be the fppf sheaf of isomorphisms of the
  $Y_i$-torsor $P^{\dag}_i$ covering the $T_i$-action on $P_i$. We
  have a subgroup $Y_i=Y{\displaystyle\mathop{\times}_S}S_i \subset G^{\dag}_i$
  and we would like to show that $G^{\dag}_i\to T_i$ is surjective.
  Let $a$ be a local section of $T_i=T{\displaystyle\mathop{\times}_S}S_i$.
  Because $\lambda(L)$ is a homomorphism, the translation action by
  $a$ on $\bPict_{P_i/S_i}$ induces the trivial action on $T_{Y,i}$.
  Indeed, locally for every $M=T^*_bL\otimes L^{-1}$ we have
  $T^*_aM\simeq M$.  This implies that the torsors $P^{\dag}_i$ and
  $T_a^*P^{\dag}_i$ are isomorphic, so that the above homomorphism is
  surjective. Therefore, we have the exact sequence$$ 
    0\to Y_i \to G^{\dag}_i \to T_i \to 1.
  $$ 
  On the other hand, $\Ext^1(T_i,Y_i)=0$ which  follows because 
  $\Ext^1(\bG_m,\bZ)=0$ (cf.\ the  proof of \cite[3.4]{SGA71}).
  Therefore, $G^{\dag}_i/S_i$ is smooth, and locally there is a
  splitting defined by mapping $T_i$ to the connected component of the
  identity $(G^{\dag}_i)^0$. Clearly, it glues to a global splitting.
  Hence, $G^{\dag}_i=Y_i\times T_i$ and $\{P^{\dag}_i\}$ is a
  $T$-invariant cover of $P$.
  
  Consider a finite flat subgroup $H_i$ of $T_{Y,i}$. Dually, there is
  a subgroup $Y_i'$ of $Y_i$ of finite index. The scheme
  $P_i'=P^{\dag}_i/Y_i'$ is proper over $S_i$ and it is an
  $\hH_i$-torsor over $P_i$, $\hH_i=Y_i/Y_i'$. Denote by $L_i'$ the
  pullback of $L_i$ to $P_i'$. The polarization morphism
  $\lambda(L_i'):T_i\to\bPict_{P_i'/S_i}$ factors through
  $\bPict_{P_i/S_i}\to\bPict_{P_i'/S_i}$, and the kernel of the latter
  homomorphism contains $H_i$. As a conclusion, the action of $H_i$ on
  $L_i'$ and its pullback $L^{\dag}_i$ is  linearizable. In particular, we
  can apply this to the family $\{_nT_{Y,i}\}$ of $n$-torsion
  subgroups of $T_{Y,i}$. Since this family is schematically dense
  (\cite[IX.4]{SGA3}) in $T_{Y,i}$, the action of $T_{Y,i}$ on
  $L^{\dag}_i$ is linearizable.
  
  The same argument works for the bigger torus $T_i$. The
  $T_i$-linearization is induced by the $T_{Y,i}$-linearization but
  now we can make it canonical: defining $\wP_i$ to be a disjoint
  union of $X/Y$ copies of $P_i^{\dag}$, on for the pullback $\wP_i$
  we obtain a canonical $T_i$-linearization. We have
  $(P_i,L_i)=(\wP_i,\wL_i)/X$. Because sheaf $\wL_i$ is ample, $\wP_i$
  is identified with a $\Proj$ of the graded algebra
  $\cO_{S_i}$-algebra $R(\wL_i)=\oplus_{d\ge0}f\wL_i^d$. This algebra
  is graded by $\bZ_{\ge0}$ and the $T_i$-linearization provides an
  additional grading by $X$ as this algebra splits into the
  eigenspaces. Finally, the actions \pagebreak of $T$ and $X$ on $\wL_i$ commute
  in the way that is claimed in the formulation of the theorem because
  that is how the actions of $X$ and $T$ on the algebra $R(\wL_i)$
  commute: $x\in X$ takes the $x'$-eigenspace to the
  $(x+x')$-eigenspace.
  
  This solves the problem locally. Now we want to glue $\wP_i$ into an
  $X$-torsor $\wP$ over $P$. Without the $T$-action, if we just had
  $\wP_i$'s, we would proceed in the following way. We would glue
  $\wP_i$'s on the ``intersections'' $S_i{\displaystyle\mathop{\times}_S}S_j$, and
  we would have a choice: we could change an identification by an
  element $x_{ij}\in X$.  But then on the triple ``intersections'' we
  may run into the incompatibility problem.  Exactly as in the exact
  sequence in the proof of \ref{cor:infinite_covers} there would be an
  obstruction in $H^2_{\rm fppf}(S,X)$ which has to vanish in order
  for the gluing to exist.  In the presence of the $T$-action,
  however, this obstruction is zero. Because each algebra $R(\wL_i)$
  is $X$-linearized, on the ``intersections'' they glue in a canonical
  way, according to the $X$-eigenspaces. This completes the proof of
  the theorem.
\enddemo

\proclaimtitle{of the proof}
\specialnumber{4.3.4 }\proclaim{{C}orollary}
  The connected components of $\wP$ are parameterized by $X/Y$ which
  is the Cartier dual of the kernel of the polarization morphism
  $\lambda(L):T\to\bPic_{P/S}${\rm .}
\endproclaim

{\it Remark} 4.3.5.
  Note that the actions of $T$ and $X$ on $\wL$ combine into an action
  of the {\it infinite Heisenberg group\/} $\cG = \bG_{m,S}
  {\displaystyle\mathop{\times}_S}  T {\displaystyle\mathop{\times}_S} X$ 
  with the following group law: $\bG_m$ is the center, and if $a$ and
  $x$ are functorial points of $T$ and $X$ (i.e.\ sections after a base
  change $S'/S$), then$$ 
    T^*_a S^*_x = x(a) \cdot S^*_x T^*_a.
  $$ 
  The surprising element here is that we are doing in Zariski topology
  things that   normally only work in the classical topology over
  $\bC$.
 
\section{Stable semiabelic varieties and pairs}
\label{sec:Stable semiabelic varieties and pairs}

  We are now going to classify stable semiabelic varieties and pairs
  over an algebraically closed field $k$ in the remaining cases: when
  there is no linearization on $L$ and also when $G$ has a nontrivial
  abelian part. We will also consider the case of a quotient field of
  a complete normal ring.

\demo{{\rm 5.1. (}\/Co\/{\rm )}\/sheaves on general cell complexes}
Let $\Delta=\wDelta/X =\Deltad/Y =\Deltadd/\pi_1(\Delta)$ be a complex
referenced by $\oXR/Y$. Let $\{\delta_{i}, \, i\in \Delta_{\min}\}$ be
the   minimal  order cells under the reverse to inclusion, and choose
their representatives $\wdelta_i$, $\deltad_i$, $\deltadd_i$ in the
complexes $\wDelta$, $\Deltad$, $\Deltadd$ resp.\ These latter
cells are either lattice polytopes or pullbacks of lattice polytopes
under a projection to a lower dimension.

We will go through the definition in the case of the cover
$\Deltad$. The other two cases, which will also be useful, are
obtained by replacing the group $Y$ everywhere formally by $X$,
$\pi_1(\Delta)$ resp.\ For simplicity, we will only consider the case
when $\Deltad$ is a join semilattice; i.e., intersection of two cells
is again a single cell.

To define (co)chain groups $C_p$ and $C^p$ consider the nonempty
intersections$$ 
  \delta_{[i_0| y_1\delta_{i_1}| \dots |y_p\delta_{i_p}]} =
  \delta_{i_0} \cap y_1.\delta_{i_1} \cap \dots y_p.\delta_{i_p}
$$ 
such that we do not have both $i_{k-1}=i_k$ and $y_k=0$. For each of
these we can take the groups $\bX_{\delta}$, $X_{\delta}$,
$\Fun(X\cap\delta,\bZ)$, $\Fun(C\cap\delta)$ and define
$C_p(\Delta,G)$ for the cosheaf $\ubX,\uX,\uFun$ resp.\ as the
direct sum of these. The differential $C_p\to C_{p-1}$ is defined as
$\sum (-1)^k d_k$, where$$ 
  d_k: G_{[i_0| y_1\delta_{i_1}| \dots |y_p\delta_{i_p}]} \to
  G_{[i_0| y_1\delta_{i_1}| \dots \hat{y_k\delta_{i_k}}
    \dots|y_p\delta_{i_p}]} 
$$ 
is the natural restriction homomorphism if $k\ne0$; for $k=0$, it is
the composition of restriction and the homomorphism$$ 
   S_{-y_1}: G_{[y_1\delta_{i_1}|y_2\delta_{i_2}| \dots
     |y_p\delta_{i_p}]} \to 
   G_{[\delta_{i_1}| y_2\delta_{i_2}|\dots |y_p\delta_{i_p}]} 
$$ 
given by the group action. We give the structure of the $Y$-module on
the above groups in the following way:
\begin{itemize}
\item[1.] on $\Fun(X,\bZ)$: $y.1_x=1_{x+y}$,
\item[2.] on $\bX$: $y.(d,x)=(d,x+dy)$,
\item[3.] on $X$: $y.x=x$, i.e.\ a trivial action.
\end{itemize}

As before, we define in this more general situation:
\begin{itemize}
\item[1.] groups of cochains $C^p(\hG)= \Hom(C_p(G), \bG_m(S)$,
\item[2.] (co)homologies,
\item[3.] $\bK$ and $\hbK$,
\item[4.] complexes $\bM_*=\Cone(C_*(\uFun)\to C_*(\ubX))$ and
  $\hbM^*=\Cone(C^*(\uhFun)\to C^*(\uhbX))$,
\item[5.] all the exact sequences of (co)sheaves and associated long exact 
  sequences of (co)homologies as before,
\item[6.] variants of these definitions with $\uFun_{\ge0}$.
\end{itemize}

When it makes a difference, we will specify which cover we use in the
above definition by adding symbols $\dag$, $\ddag$ or $\tilde{}$
respectively. 

\demo{Example {\rm 5.1.1}}
  If $Y=0$, we recover the previous Cech (co)homologies on a finite
  complex of lattice polytopes.
\enddemo

{\it Example} 5.1.2.
  If $Y\ne0$ but $\Deltad$ is a complex of lattice polytopes, we get
  (co)homologies of just slightly more complicated constructible
  (co)sheaves on the topological space $|\Delta|$.

\demo{Example {\rm 5.1.3}} 
  Consider the complex $\Delta$ consisting of one big cell $\XR/Y$,
  where $Y$ is a sublattice in $X$ of finite index. In this case, we
  obtain the usual group (co)homologies of the corresponding
  $Y$-modules.
\enddemo

\specialnumber{5.1.4}\proclaim{Lemma}\label{lem_homs_one_huge_cell}
 {\rm 1.}  $H_1(\ubX)=(Y\otimes X)/\Lambda^2Y\supset (Y\otimes
    Y)/\Lambda^2Y = \Gamma^2(Y)$.
\begin{itemize}
\ritem{2.}  $B_1(\ubX)\subset C_1(\ubX)=\oplus_{y\ne0}[y]\bX$ is the
    subgroup generated by the elements $([y_1+y_2]-[y_1]-[y_2])(0,x)$
    and $([y_1+y_2]-[y_1]-[y_2])(1,0)-[y_2](0,y_1)$.
  \end{itemize}

\endproclaim

{\it Proof}.
  From the long exact sequence of homologies for the short exact
  sequence$$ 
    0\to \uX \to \ubX \to \bZ \to 0 
  $$ 
  we get $H_1(\ubX)= H_1(X)/H_2(\bZ)$. On the other hand,
  $H_p(\bZ)=H_p(|\Delta|,\bZ) = \Lambda^pY$. Also, since $X$ is a
  trivial $Y$-module, $H_p(X)=H_p(\bZ)\otimes X$. This proves Part 1.
  Part 2 is a direct calculation.
\hfill\qed

\specialnumber{5.1.5}\proclaim{Lemma}\label{lem:homologies_of_L} For $\Delta=\XR/Y${\rm :}
  \begin{itemize}
  \ritem{1.} $C_0(\ubL)=\bL$ is the group generated by the elements
    $1_{x_1+x_2}-1_{x_1}-1_{x_2}+1_0${\rm .}
  \ritem{2.} $B_0(\ubL)$ is the subgroup of $\bL$ generated by the elements$$ 
      (1_{x_1+x_2+y}-1_{x_1+y}-1_{x_2+y}+1_y)-
      (1_{x_1+x_2}-1_{x_1}-1_{x_2}+1_0).
    $$ 

    Assume in addition that $Y=X${\rm .} Then\/{\rm :}
  \ritem{1.} $H_0(\ubL)=C_0/B_0(\ubL)$ is the lattice of dimension
    $r(r+1)/2$ isomorphic to $\Sym^2X$ via
    $1_{x_1+x_2}-1_{x_1}-1_{x_2}+1_0\mapsto x_1\otimes x_2$. The dual
    to it is the lattice of homogeneous quadratic functions on $X${\rm .}
  \ritem{2.} $C_0(\uFun)/B_0(\ubL)$ is the lattice of dimension
    $r(r+1)/2+(r+1)$. The dual to it is the lattice of nonhomogeneous
    quadratic functions on $X${\rm .}
  \end{itemize}

\endproclaim

\vglue-12pt
{\it Proof}.
  This is obtained directly from definitions.
\hfill\qed
\vglue4pt
 
5.2. {\it Linearized varieties with a nontrivial abelian part}.
 
\specialnumber{5.2.1}\proclaim{Lemma}
  Let $(P,L)$ be an irreducible variety together with an action of a
  semiabelian variety $G${\rm ,} $1\to T\to G\to A\to 0${\rm ,} such that the
  $T$\/{\rm -}\/action on $L$ is linearized and the condition  on orbits is
  satisfied{\rm .} Then $P$ is fibered over $A$ {\rm (}\/more precisely{\rm ,} over an
  $A$\/{\rm -}\/torsor\/{\rm ),} all fibers are isomorphic and the normalization of each
  is an irreducible projective toric variety with a linearized line
  bundle corresponding to a lattice polytope $\delta\subset \XR${\rm .} 
\endproclaim

{\it Proof}.
  Pick an arbitrary point $p$ in the interior of the main orbit and
  let $P_0$ be the closure of the orbit $Tp$. For each coset $T_a$ of
  $T\subset G$ the closure of $T_ap$ is a subvariety isomorphic to
  $P_0$. It is clear that $P_0$ satisfies our condition on the orbits
  as well, for the $T$-action, so that its normalization together with the
  pullback of $L$ is a projective toric variety corresponding to a
  lattice polytope $\delta\subset\XR$. The $0$-dimensional orbits of
  $T\acts P_0$ correspond to the vertices of this polytope. If $q_x$
  is one of these orbits then by our condition on the orbits,
  $Gq_x=Aq_x\simeq A$. The restriction of $L$ to $Aq_x$ is
  $T$-linearized and is of weight $x$; therefore for $x\ne y$ one has
  $Aq_x\cap Aq_y=\emptyset$. This implies that different varieties
  $P_a$ do not intersect and that they define a fibration of $P$ over $A$.
\hfill\qed

\specialnumber{5.2.2}\proclaim{Lemma}
  Let $(P,L)$ be a stable semiabelic variety with a $T$\/{\rm -}\/linearized
  line bundle{\rm .} Then every irreducible component of $P$ is normal and
  is also a stable semiabelic variety{\rm .}
\endproclaim

{\it Proof}.
  It is clear for any point $p\in P$ that the above fibration on the
  irreducible components defines a fibration of a neighborhood with
  isomorphic slices. The seminormality of $P$ implies the
  seminormality of each of these slices, and, by our results in the
  toric case, the normality of each irreducible component.
\hfill\qed\vglue4pt

  Recall that the semiabelian variety $G$ is defined, via the negative
  of pushout, by a homomorphism $c:X\to A^t$. If $\cO_x=c(x)$ is the
  corresponding rigidified sheaf then$$ 
    G= \Spec_A \oplus_{x\in X} \cO_{x}.
  $$

\demo{Definition {\rm 5.2.3}}
  Let $\cM$ be an ample invertible sheaf on $A$. For $\chi=(d,x)\in
  \bX$ we will denote by $\cM_{\chi}$ the  
  sheaf $\cM^d\otimes \cO_x$ rigidified at the origin.
\enddemo

\specialnumber{5.2.4}\proclaim{Lemma}
  Let $(P,L)$ be an irreducible semiabelic variety with a linearized
  line bundle{\rm .} Fix a point $p\in P$ in the main $G$\/{\rm -}\/orbit{\rm .} With this
  choice{\rm ,} $P$ can be identified with the variety$$ 
    \Proj_A \oplus_{\chi\in \Cone\,\delta} \cM_{\chi},
  $$ 
  and $L$ with the sheaf $\cO(1)${\rm .}
\endproclaim

\demo{Proof}
  With the fixed projection $f:P\to A$ consider the $\cO_A$-algebra
  $\cA=\oplus_{d\ge0}L^d$. Since the fibers of $f$ are irreducible
  stable toric varieties, we know that the higher direct images
  $R^if_*L^d$ are zero, and the algebra $\cA$ is locally free. The
  $T$-action on $L$ defines a splitting of $\cA$ into a direct sum of
  eigenspaces, and since on each $f^{-1}(a)$ the eigenspaces are
  1-dimensional, these $\cO_A$-sheaves are invertible.  Say, the
  $(1,x)$-eigenspace is $\cM\otimes\cO_x$ for some $\cM$. Then the
  prescribed form is dictated by the fact that $\cA$ is an algebra and
  $P$ has $G$-action. Finally, $\cM$ has to be ample since from the
  growth of $h^0(L^d)$ we see that $h^0(\cM^d)$ has to grow as $g^d$,
  $g=\dim A$, and on an abelian variety a big sheaf is ample. Clearly, 
  $(P,L)=(\Proj\cA,\cO(1))$.
\enddemo

{\it Definition} 5.2.5.
  We will denote the polarized variety in this standard form by
  $(P,L)[\delta,c,\cM]$.

\specialnumber{5.2.6}\proclaim{Lemma}\label{thm:cohs_irred_SSAVs}
  Let $(P,L)=(P,L)[\delta,c,\cM]${\rm .} Then $H^p(P,L)=0$ for $p>0$ and
  there is a canonical isomorphism$$ 
    H^0(P,L)=\oplus_{x\in X\cap\delta} H^0(A,\cM_{(1,x)}).
  $$ 
\endproclaim

\demo{Proof}
  For the global sections this is obvious from the above, since\break
  $H^0(P,L)=H^0(A,f_*L)$. The higher direct images $R^pf_*L$ vanish
  since we have vanishing on the fibers which are projective toric
  varieties. Therefore, $H^p(P,L)=H^p(A,f_*L)$ and these are zero
  since higher cohomologies of ample sheaves on abelian varieties
  vanish. 
\enddemo

\specialnumber{5.2.7}\proclaim{Theorem} \hglue-8pt
  Let $\theta\in H^0(P,L)$ be a section corresponding to a divisor~$\Theta${\rm .} Then $(P,\Theta)$ is a stable
semiabelic pair if and only
  if the homogeneous components for the $x\in\Ver\, \delta$ are not
  zero{\rm .} 
\endproclaim

\demo{Proof}
  Indeed, these homogeneous components define the restrictions of
  $\Theta$ to each minimal $G$-orbit.
\enddemo

{\it Definition} 5.2.8.
  We will denote the subset of $H^0(A,\cM_{(1,x)})$ with the above
  condition by $\hFun_{\ge0,\delta}$.

\demo{{\rm 5.3.} Arbitrary polarized {\rm SSAV}s}
  Now, let $(P,L)$ be an arbitrary polarized stable semiabelic
  variety.  By the main result of Section~\ref{sec:Linearization of
    torus action} the pair $(P,L)$ is equivalent to a pair $(\wP,\wL)$
  together with an $X$-action, where $\wL$ is $T$-linearized but the
  scheme $\wP$ is only {\it locally} of finite type over $k$.  For
  schemes which are locally of finite type over a field the
  seminormality is a formal property, and it is certainly preserved
  under an \'etale morphism.  Therefore, any finite union of
  irreducible components of $\wP$ together with the restriction of
  $\wL$ is a linearized polarized stable semiabelic variety.
  Therefore, the pair $(\wP,\wL)$ defines a complex $\wDelta$ of
  lattice polytopes which is {\it locally} of finite type, together
  with an $X$-action compatible with the reference map
  $\tilde\rho:|\wDelta|\to X_{\bR}$.
  
  The $X$-action on $\wDelta$ is not properly discontinuous in the
  order topology. However, it {\it is} if we work instead with the
  classical topology on $|\wDelta|$. Introduce the complex $\Delta$ as
  the quotient $\wDelta/X$. If $\Deltad$ is a connected component of
  $\wDelta$ and $Y\subset X$ the subgroup leaving it invariant then
  $\Delta= \Deltad/Y$. The complex $\Delta$ is naturally referenced
  with a map to $\oX_{\bR}/Y$, where $\oX$ is an $X$-torsor, i.e., the
  group $\bZ^r$ but without the origin fixed. Moreover, for every
  $X$-linearized cosheaf $G$, resp., sheaf $F$, on $|\wDelta|$ we can
  consider the corresponding quotient cosheaf, resp., sheaf, by the
  free and properly discontinuous $X$-action. We will denote these
  quotient cosheaves and sheaves on $|\Delta|$ by the same letters.
\enddemo

{\it Definition} 5.3.1.
  We introduce one more complex $\Deltadd$, the universal cover of
  $\Deltad$. It comes with the reference map to $\XR$ and corresponds
  to a polarized scheme $(\Pdd,\Ldd)$ locally of finite type such that
  $(P,L)=(\Pdd,\Ldd)/\pi_1(|\Delta|)$.  We will assume for simplicity
  that the complex $\Deltadd$ is a poset; i.e., the intersection of two
  polytopes is again a polytope (the general case is only a little
  harder).  
 
\demo{Definition {\rm 5.3.2}}
  Pick an arbitrary point $p$ in a minimal orbit $Ap\subset P$. For
  each cycle in $\pi_1(|\Delta|)$ we can trace a cycle of $T$-orbits
  on our variety along the fibrations to come back to another point
  which is a shift of $p$ by the $A$-action. This way, the variety $P$ 
  defines a homomorphism from $\pi_1\Delta$ to $A$ which we will
  denote by $c^t$.  
\enddemo

5.3.3.
  Choose a system $\{\delta_{i}, \,|\, i\in I\}$ of representatives in
  $\Deltad$ of all the minimal cells in $\Delta$. Also, choose a point
  in $\Pdd$. With these choices made, we can identify the subvarieties
  corresponding to $\delta_i$ with the standard polarized varieties
  $P[\delta_i,c,\cM]$ for the same $c,\cM$.  We can do so uniquely:
  since there are no cycles in $\Deltadd$, any tracing out gives the
  same origin in $A$.  All other irreducible components of $\Pdd$
  correspond to the cells $z.\delta_i$ for some $z\in\pi_1\Delta$.
  Hence, we can also identify all the irreducible components of $\Pdd$
  with the translations of these standard varieties. To give the
  variety $(P,L)$ is the same as to give the scheme $(\Pdd,\Ldd)$
  together with the action of $\pi_1\Delta$. Therefore, we see that we
  have to describe the gluing of the shifted algebras
  $T^*_{c^tz}R_{i}$ for the cells $z.\delta_i$ on the intersections.
  To this end, we see first of all that these corresponding
  invertible sheaves have to be isomorphic. This translates to the
  condition that
  \begin{equation}\label{eqn:commutin_condition}
     \lambda(\cM) \circ c^t=   c \circ i: \,  \pi_1\to A^t  \speqnu{3} 
  \end{equation}
  where we denoted by $i$ the projection $\pi_1\Delta\to X$ given by
  the reference map $\Deltadd\to\XR$.
 
\demo{Definition {\rm 5.3.4}}
  Assuming $c,c^t$ and $\cM$ given, we will define the cochain groups
  $C_{\ddag}^p(\Delta,\uhbX)$ for $p=0,1$ as functions on
  $C^{\ddag}_p(\Delta,\ubX)$ linear in $\bX$, where the function on
  \begin{itemize}
  \item[1.] $[\,]\chi$ takes values in $\bG_m(k)$,
  \item[2.] $[z]\chi$ takes values in the $\bG_m(k)$-torsor 
    $\Iso(\cM_{\chi}, T^*_{-c^tz}\cM_{z.\chi})$ which is the same as
    the fiber $\cM_{\chi}^{-1}(c^tz)$ without the zero.
  \end{itemize}

\enddemo

\specialnumber{5.3.5}\proclaim{Lemma}
  The functions on $B_1(\ubX)$ take value canonically in $\bG_m(k)${\rm .}
\endproclaim

\demo{Proof}
  We have to check this in two cases: $\chi=(0,x)$ and $\chi=(1,0)$.
  In the first case   this reduces to a canonical identification of
  $\bG_m$-torsors $\cO_x(c^tz_1+c^tz_2)=\cO_x(c^tz_1)\otimes \cO_x(
  c^tz_2)$ which is a consequence of the so-called bi-extension
  structure on the Poincar\'e bundle on $A\times A^t$. In the second
  case   this reduces to the identity$$ 
    \cM(c^tz_1+c^tz_2)\otimes
    \cM(c^tz_1)^{-1} \otimes 
    \cM(c^tz_2)^{-1}=\cO_{iz_1}(c^tz_2)    
  $$ 
  which is to say that the symmetric bi-extension on $A\times A$
  coincides with the pullback of the Poincar\'e bi-extension on $A\times
  A^t$ via $(\id,\lambda(\cM))$. This follows from our commuting
  condition 3.
\enddemo

{\it Definition} 5.3.6.
  We will denote by $Z^1(\uhbX)$ the functions on $C_1/B_1(\ubX)$,
  i.e.\ the functions on $C_1$ that are identically 1 on $B_1$.

\demo{Definition {\rm 5.3.7}}
  A framing for a polarized stable semiabelic variety will have the
  following characteristics:
  \begin{itemize}
  \item[1.] isomorphism $X\isoto\, \bZ^r$,
  \item[2.] the connected component $\Deltad$,
  \item[3.] choice of the representatives for the minimal cells
    $\delta_i$, 
  \item[4.] projections to $A$ of the varieties for these selected cells. 
  \end{itemize}
\enddemo

  Under   minor modifications of the definitions of cohomologies,
    Theorem \ref{thm:polarized_linearized_STVs} holds literally:

\specialnumber{5.3.8}\proclaim{Theorem}\label{thm:classif_arbitrary_polarized_SSAV}
  The groupoid $M^{\framed}[\Delta,c,c^t,\cM](k)$ is equivalent to$$ 
    [  Z_{\ddag}^1(\Delta,\uhbX) / C_{\ddag}^0(\Delta,\uhbX) ].
  $$ 
\endproclaim

\demo{Proof}
  Indeed, elements of $Z^1$ describe the compatible gluings, and
  $C^0$ describes the choices of an origin in each of the  $P[\delta_i,c,\cM]$
  that do not change the projection to $A$.
\enddemo

  Concerning the singularities, we have the following, which is the same as  the
  toric case.

\specialnumber{5.3.9}\proclaim{Theorem}\label{thm:CM_SSAV}
  Let $(P,L)$ be a stable semiabelic variety corresponding to a
  complex $\Delta$ plus all the other data{\rm .} If $\Delta$ is locally {\rm CM}
  {\rm (}\/for example $|\Delta|$ is a manifold with a boundary\/{\rm )} then $P$ is
  {\rm CM. }
\endproclaim

\demo{Proof}
  For the schemes which are locally of finite type over a field,
  the property of being Cohen-Macaulay is a formal property. Therefore, we can go over
  to the \'etale cover $\wP$ for this question. Additionally, for any
  point $p\in P$ the completion of the local ring splits as the direct 
  product of the abelian part (smooth) and the toric part, the latter
  being the same as for a stable toric variety for the same
  complex. The result now follows.
\enddemo

5.4. {\it Arbitrary stable semiabelic pairs}.
 
\specialnumber{5.4.1}\proclaim{Theorem}\label{thm:cohs_polarized_SSAVs}
  For any polarized stable semiabelic variety $(P,L)$ one has
  $H^p(P,L^d) =0$ for $p,d\ge1${\rm .}
\endproclaim

\demo{Proof}
  The irreducible linearized case was established in Theorem
  \ref{thm:cohs_irred_SSAVs}. 
  
  It is sufficient to prove the statement  for a finite cover with the
  pullback of $L$. This finite cover can be chosen to be again an SSAV
  but with the condition that the intersections of irreducible
  components are irreducible, i.e., that intersection of two polytopes
  in $\Delta$ is a single polytope. As in the proof for the toric
  case, consider the cover of $P$ by the varieties
  $P_{i_0}$ corresponding to $\Umin$. Denoting
  $P_{i_0i_1}=P_{i_0}\cap P_{i_1}$ and so forth (they are all also SSAVs of
  lower dimension), we get a resolution of $\cO_P$ by the complex
  having $\oplus\cO_P{i_0\dots i_p}$ in degree $p$. Computing the
  hypercohomologies of this complex and using   induction on the
  dimension we complete  the proof.
\enddemo

{\it Definition} 5.4.2.
  Keeping the analogy with the linearized toric case, for an SSAV
  described in the previous section we   define
  \begin{itemize}
  \item[1.] $C^0(\Delta,\uhFun\lge)= \oplus \hFun_{\ge0,i}$ (see
    Definition~5.2.8),
  \item[2.] $\hbM^0(\Delta)=C^0(\Delta,\uhbX)$,
  \item[3.] $\hbM^1(\Delta)= C^0(\Delta,\uhFun\lge)\oplus
    C^1(\Delta,\uhbX)$.
\end{itemize}
 Note that the homogeneous part $\tau_{i_1i_0}^{[z]\chi}$ of
    $\tau\in C^1(\Delta,\uhbX)$ which is the value on $[z]\chi$ is an
    element of $\Iso(\cM_{\chi}, T^*_{-c^tz}\cM_{z.\chi})$ and it
    induces an isomorphism from $H^0(P_{i_0},\cM_{\chi})$ to
    $H^0(P_{i_1},T^*_{-c^tz}\cM_{z.\chi})=   T_{-c^tz}
    H^0(P_{i_1},\cM_{z.\chi})$.  We will denote by $Z^1(\hbM^*)$ the
    set of pairs $(f,\tau)\in \hbM^1(\Delta)$ such that
    $f_{i_1}=T_{-c^tz}\left(\tau_{i_1i_0}(f_{i_0})\right)$.
\vglue12pt

  With these straightforward generalizations,
  Theorem~\ref{thm:linearized_STpairs} remains literally true:

\specialnumber{5.4.3}\proclaim{Theorem}\label{thm:arbitrary_polarized_SSApairs}
  The groupoid $\MP^{\framed}[\Delta,c,c^t,\cM](k)$ is equivalent to$$ 
    [  Z_{\ddag}^1(\Delta,\hbM^*) / \hbM_{\ddag}^0(\Delta) ].
  $$ 
\endproclaim

{\it Proof}.
  With the previous theorem at hand, the proof goes exactly as in the
  polarized toric case.
\hfill\qed
\pagebreak

5.4.4.
  It is now time to describe the ``unframed'' groupoid of all stable
  semiabelic pairs. We will start with a complex $\Delta$ referenced
  by $\oX_{\bR}$ and a fixed abelian variety $A$ with a polarization
  $\lambda$.  First, we have the space $\Hom'(\pi_1\Delta\times X,
  A\times A^t)$ of homomorphisms $(c^t,c)$ satisfying the commuting
  condition \ref{eqn:commutin_condition}. Let $\Pic^{\lambda}A$ be the 
  component of $\Pic\, A$ containing the sheaves with
  $\lambda(\cM)=\lambda$. Then over this space we have the family of
  $Z^1(\Delta,\hbM^*)$. To arrive at $\MP[\Delta,A,\lambda]$ we have to
  divide this family by three consecutive equivalence relations:
  \begin{itemize}
  \item[1.] by the action of $\hbM^0(\Delta)$ -- this one is the same for
    each fiber,
  \item[2.] by the action of $A$ which corresponds to a change of origin
    in $A$,
  \item[3.] and, finally by the groups of symmetries of the complex
    $\Delta$ and the pair $(A,\lambda)$.
  \end{itemize}
  Here, the first one is given by a torus action but the action is
  proper with finite stabilizers, by the semigroup
Lemma~\ref{lem:semigroups_finite_action}. The second one is a
  quotient by a proper action of a proper group. And finally the third 
  one is a quotient by a finite group.

  Moreover, we can do the same over a stack of all abelian varieties
  of dimension $g$ and a polarization of degree $d$. By the result of
  Keel-Mori [KM] this gives a separated stack
  $\MP[\Delta,g,d]$ over $\bZ$ which admits a coarse moduli space as a 
  separated algebraic space.

\demo{{\rm 5.5.} Stable semiabelic pairs over $\bC$}
To describe stable semiabelic varieties and pairs over an arbitrary
closed field we used a complex of lattice polytopes $\Delta'$ (we
called it simply $\Delta$ before to simplify notation)
referenced by $\oXR'/Y'$, where $X'\simeq \bZ^r$, $r\le g$, was the
character group of the toric part. The formulas were a little
complicated since they involved cohomologies with values in certain
$\bG_m$-torsors. 

A variety or a pair over the field of complex numbers can be
described much more directly, by using a complex $\Delta$ with cells
of full dimension which are not necessarily polytopal. The reason for
this is that a $g$-dimensional semiabelian variety $G$ can be,
noncanonically, represented as $(\bC^*)^g/Y$ and so a variety $P$ with
$G$-action can be formally considered to be ``toric'' (this action is
not algebraic, of course). An abelian variety corresponds in this way 
to one big cell $\oX/Y$. A general semiabelic variety corresponds to a
general cell complex. As in the toric case, a variety and its complex
are related through the moment map $\mom:P(\bC)\to |\Delta|$. A fiber
of this map over a point in the interior of a $k$-dimensional cell
$\delta$ is isomorphic to $(S^1)^k$.

It turns out that the classification theorem for stable toric
varieties and pairs hold here verbatim, with one small difference:
the varieties and pairs are classified not by entire groups
$H^1(\uhbX)$, resp.\ $H^1(\hbM^*)$, but by open subsets of points that
satisfy the Riemann positivity condition.  The dimension of the
corresponding stratum is, therefore, still $h^1(\uhbX)$, resp.\ $h^1(\hbM^*)$. Since our strata were constructed over $\bZ$ and the
dimension does not depend on a particular field, the same dimension
formula holds over an arbitrary field $k$.

\vglue6pt 5.5.1.
  We will start with a single abelian variety $A=\bC^g/\bZ^{2g}$ with
  a polarization $\lambda$ of type $\bZ_{d_1}\times \dots \bZ_{d_g}$.
  Recall from the classical theory that via exponentiation $\exp(2\pi
  i\cdot)$ of half of the periods $(A,\lambda)$ can be represented in
  the following way:
  \begin{itemize}
  \item[1.] $A=T/i(Y)$, $T=(\bC^*)^g$, $Y=\bZ^g$. If we denote by $X$ the
    group of characters of $T$ then we get a bilinear function
    $b:Y\times X\to \bC^*$, $b(y,x)= x(y)$;
  \item[2.] $Y\subset X$ is a sublattice with $X/Y=\bZ_{d_1}\times \dots
    \bZ_{d_g}$ such that $b|_{Y\times Y}:Y\times Y\to \bC^*$ is
    symmetric and positive definite in the following sense:
    $|b(y,y)|<1$ for all $y\ne 0$.
  \end{itemize}
  There is a group $\SP_{2g}(\diag(d_1,\dots,d_g), \bZ)$ acting
  properly discontinuously in classical topology on the set of all
  such representations, and the quotient is the component of $\Agd$
  corresponding to polarizations of type $(d_1,\dots,d_g)$. The
  subgroup of $\SP_{2g}(\diag(d_1,\dots,d_g), \bZ)$ that preserves the
  structure of a ``toric variety'' on $A$, i.e.\ the presentation
  $A=\hX/Y$ is a subgroup $\GL(X,Y) \subset \GL(X)$ that sends $Y$ to
  itself.  The remaining equivalence relation  by which we have to divide in order to arrive at the
  quotient by $\SP(\diag(d_1,\dots,d_g), \bZ)$ corresponds to a
  different choice of a maximal isotopic subspace in $\bZ^{2g}$.

  The datum for an ample sheaf $L$ on $A$ with $\lambda(L)=\lambda$ is
  a function\break $a:Y\to \bC^*$ such that
  $a(y_1+y_2)=a(y_1)a(y_2)b(y_1,y_2)$.

\specialnumber{5.5.2}\proclaim{Lemma}
  The datum above is precisely the same as an element of  
  $Z^1(\oXR/Y,\uhbX)=\Hom(C_1/B_1(\oXR/Y,\ubX),\bC^*)$ whose
  restriction to the subgroup $\Sym^2 Y$ is positive definite{\rm .}
\endproclaim

{\it Proof}.
  Indeed, by Lemma~\ref{lem_homs_one_huge_cell} an element of
  $\Hom(C_1/B_1(\oXR/Y,\ubX),\bC^*)$ is a function $\tau(y,d,x)$ on
  $Y\times \bX$ which is linear in $\bX$ and which satisfies two types
  of relations. If we denote $\tau(y,0,x)=b(y,x)$ and
  $\tau(y,1,0)=a(y)$ then by linearity we can write
  $\tau(y,d,x)=b(y,x)a(y)^d$ and the relations then take form:
  \begin{itemize}
  \item[1.] $b(y_1+y_2,x)=b(y_1,x)b(y_2,x)$,
  \item[2.] $a(y_1+y_2)=a(y_1)a(y_2)b(y_1,y_2)$,
  \end{itemize}
  which is exactly what we have above.
\hfill\qed\pagebreak

Hence, formally, $(\bC^*)^g\to A$ plays the role of the cover
$P^{\dag} \to P$ that we had before, and the line bundle
$(\bC^*)^g\times \bC$ plays the role of $L^{\dag}$ with $L=L^{\dag}/Y$. 
We still have the  algebra $\Rdag= \bC \oplus_{d>0}H^0(\Pdag,(\Ldag)^d)$
and $\Pdag =\Proj\, \Rdag$, even though it is no longer a projective
variety since $\Rdag_+$ is not finitely generated.

Note that to classify abelian torsors we have to divide by a choice of
the origin in the above representation, which is the same as dividing
by the action of $C^0(\uhbX)=\hbX=\bT$. As a result, polarized abelian
torsors over $\bC$ are classified by an open subset of
$H^1(\oXR/Y,\uhbX)$ modulo $\Sym(\oXR/Y)$, which contains $\GL(X,Y)$
as a subgroup of finite index. This set is precisely the set of
bilinear functions $b(y,x)$ whose restriction on $Y\times Y$ is
symmetric and positive definite.

Next, consider sections in $H^0(A,L)$. We start with functions
$\zeta^{(1,x)}$ on $(\bC^*)^g\times \bC =\Spec\, \bC[\bZ_{\ge0} \times
X]$ corresponding to elements $(1,x)\in \bZ_{\ge0}\times X$. These are
sections of $H^0(\wP,\wL)=H^0(T,\cO)$.  According to the classical
theory of theta functions, a section of $H^0(A,L)$ can be represented
as a $Y$-invariant section of $H^0(T,\wL)$ given by a power series$$ 
  \sum_X \xi_x = \sum_X f_x (1,x) = \sum_X f_x \zeta^{(1,x)}
$$ 
with $f_x\in \bC$ satisfying $f_{x+y}=a(y)b(y,x)f_x$.  The convergence
of this power series reduces to the fact that the radius of
convergence of the series $\sum \rho^{n^2} z^n$ is $+\infty$ if
$|\rho|<1$.

\specialnumber{5.5.3}\proclaim{Lemma}
  Abelic pairs over $\bC$ with a given $T$\/{\rm -}\/action are classified by an
  open subset in $H^1(\oXR/Y, \hbM^*)=H^0(\oXR/Y, \uhbL)$ whose
  restriction to $\Sym^2 Y$ is positive\/{\rm -}\/definite{\rm ,} modulo
  $\Sym(\oXR/Y)${\rm .}
\endproclaim

\demo{Proof}
  First of all, in this case $\bK=0$, so that the exact sequence of
Lemma~\ref{2.2.7}.  (4) gives $H^1(\oXR/Y,
  \hbM^*)=H^0(\oXR/Y, \uhbL)$. An element of $Z^1(\hbM^*)$ is by definition
  a pair $(\tau,f)$ with $\tau$ as before, satisfying
  $f_{x+y}=\tau(y,1,x)f_x$. But $\tau(y,1,x)=\tau(y,1,0)\tau(y,0,x) =
  a(y) b(y,x)$. This proves that pairs $(A,\Theta)$ are classified by
  the open subset of $Z^1(\hbM^*)$. Dividing by the choice of the
  origin, i.e.\ by the action of $C^0(\hbM^*)=C^0(\uhbX)=\bT$ gives
  $H^1(\hbM^*)$. Finally, we have to divide by automorphisms of the
  pairs $(X,Y)$ and by the choice of origin in $X$.
\enddemo

5.5.4.
  Next, consider an (irreducible) semiabelic variety $P_{\delta'}$
  with an action by a semiabelian variety $G$, $1\to T'\to G\to A\to
  0$, $T'=\hX'$, $X'\simeq\bZ^r$, and with a $T'$-linearized ample
  sheaf $L_{\delta'}$. As we saw, this variety corresponds to a
  polytope $\delta'$ referenced by the lattice $X'$. Now make a choice
  of a decomposition $G=T/Y=\hX/Y=(\bC^*)^g/Y$, so that the morphism
  $T'\to T/Y$ extends to an injection $T'\to T$ and, dually, to a
  surjection $X\to X'$. Let $X_1=\ker(X\to X')$. If the sheaf $\cM$ on
  the abelian part $A$ defines a polarization of type
  $(d_1,\dots,d_a)$, then choose a sublattice $Y_1\subset X_1$ with
  $X_1/Y_1=\bZ_{d_1}\times\dots \bZ_{d_a}$ and consider the complex
  $\Delta$ referenced by $\XR/Y_1$ consisting of one cell
  $\delta/Y_1$, where $\delta$ is the pullback of $\delta'$ under the
  projection $\XR\to \XR'$. The standard exact sequence $0\to \uX\to
  \ubX\to \bZ\to 0$ together with an embedding of the constant cosheaf
  $Y_1$ into $\uX$ give, as before, an embedding $\Sym^2Y_1\to
  H^1(\delta/Y_1,\uhbX)$.
 
\specialnumber{5.5.5}\proclaim{Lemma}
  Polarized semiabelic varieties $(P,L)$ with linearized $L$ with a
  fixed $T$\/{\rm -}\/action are classified by an open subset in
  $H^1(\delta/Y_1, \uhbX)$ of positive\/{\rm -}\/definite classes {\rm (}\/i.e.\ those
  that restrict to positive-definite functions on $\Sym^2 Y_1${\rm ),}
  modulo the group of symmetries of $\delta/Y_1${\rm .}  Semiabelic pairs
  with $(P,\Theta)$ with linearized $L=\cO(\Theta)$ are classified by
  a similar quotient of an open subset in $H^1(\delta/Y_1, \hbM^*)$ of
  positive\/{\rm -}\/definite classes{\rm .}
\endproclaim

\demo{Proof}
  Write $X=X_1\oplus X'$. As an easy computation shows, a $1$-cocycle
  in this case is a function $\tau(y,d,x_1,x')$ which by linearity can 
  be written as $b(y,x_1)a(x_1)^d c(y,x')$ and functions $a,b,c$
  satisfy:
 \vglue4pt
\hangindent=36pt\hangafter =1
 1.  $b:Y_1\times X_1\to \bC^*$ is bilinear and its restriction to
    $Y_1\times Y_1$ is symmetric and positive definite,
 
 \vglue4pt
\hangindent=18pt\hangafter =1 2.  $a(y_1+y_2)=a(y_1)a(y_2)b(y_1,y_2)$,

 \vglue4pt
\hangindent=18pt\hangafter =1  3.  $c:Y_1\times X' \to \bC^*$ is bilinear.
  \vglue4pt
  The first two pieces of the data are equivalent to giving an abelian
  variety $A$ and an ample sheaf $\cM$ on it. The third piece is a
  homomorphism $X'\to \hY_1=(\bC^*)^g$. The automorphism group of
  $\delta/Y_1$ is a finite extension of the group of linear
  transformations of $X$ given by arbitrary homomorphisms $\phi:X'\to
  X_1$. This means that we have to divide the set of $c$'s by all
  homomorphisms obtained from $b$ via some $\phi$. Hence, the revised
  third piece of the data is a homomorphism $X'\to \hY_1/j(X_1)$. But
  the above is precisely $\Pic^0 A$, and a homomorphism $X'\to \Pic^0
  A$ is equivalent to the extension $1\to T'\to G\to A\to 0$. Hence, a
  class in $H^1(\delta/Y_1, \hbM^*)$ modulo automorphisms is the same as
  giving the semiabelian variety $G$ and ample sheaf $\cM$ on $A$. The
  rest follows easily.
\enddemo

We can now take up the general case. As we already know, an arbitrary
polarized stable semiabelic variety (resp.\ pair) defines a complex
$\Delta'$ of lattice polytopes referenced by $\oXR'/Y'$ (resp.\ a  pointed 
complex) and a polarization of type $(d_1,\dots d_a)$ on the abelian
part $A$. Again, write $G=\hX/Y_1=(\bC^*)^g/Y_1$ and consider
\begin{itemize}
\item[1.] the complex $\Deltad$ which is the pullback of $\Delta'$ under
  $X\to X'$,
\item[2.] the complex $\Delta=\Deltad/Y_1$ with compact cells. It is
  referenced by $\oXR/Y$, where $0\to Y_1\to Y\to Y'\to 0$ is a
  natural split extension.
\end{itemize}
As before, we have a natural embedding $\Sym^2 Y_1\to
H_1(\Delta,\ubX)$, so that we can talk about positive definite cohomology
classes.

\specialnumber{5.5.6}\proclaim{Theorem}\label{thm:SSAVs_and_pairs_over_C}
  Over $\bC${\rm ,} polarized stable semiabelic varieties of type $\Delta$
  {\rm (}\/resp.\ pairs of type $(\Delta,C)${\rm )} with a fixed $T$\/{\rm -}\/action are
  classified by an open subset of positive definite classes in
  $H^1(\Delta,\uhbX)$ {\rm (}\/resp.\ $H^1((\Delta,C), \hbM^*)${\rm ),} modulo the
  group of symmetries of $\Delta${\rm .}
\endproclaim

\demo{Proof}
  Indeed, $Z^1$ is precisely the gluing data for the elementary blocks 
  that we have just considered, and the action of $C^0$ corresponds to 
  changing the origins in these blocks.
\enddemo

5.5.7.
  Finally, we want to define a moment map for a stable semiabelic pair
  $(P,\Theta)$ corresponding to a complex $\Delta$ referenced by
  $\oXR/Y$. Let $\theta\in H^0(P,L)$ be an equation of $\Theta$. By
  choosing a representation $A=(\bC*)^g/Y_1$ as above we find covers
  $P^{\dag}$ and $\wP$ with a not necessarily algebraic toric action
  such that $(P,L)=(P^{\dag},L^{\dag})/Y = (\wP,\wL)/X$. We will work
  with the second cover. As above, the sheaf $\wL$ comes with a
  section $\widetilde\theta=\sum_X \xi_x$ which is a (converging
  everywhere) power series. Just as in the finite case, we define$$ 
    \mom: \wP\to |\wDelta|; \qquad
    p\mapsto \frac{\sum |\xi_x(p)|^2 x}{\sum |\xi_x(p)|^2}.
  $$ 
  We leave it to the reader to check that:
  \begin{itemize}
  \item[1.] The above power series is everywhere converging as well;
    this reduces to the fact that the power series $\sum n \rho^{n^2}
    z^n$ has a radius of convergence $+\infty$ if $|\rho|<1$.
  \item[2.] $\mom$ commutes with the $X$-actions on $\wP$ and $\wDelta$,
    and hence descends to a moment map $\mom: P\to |\Delta|$.
  \item[3.] A fiber over a point in the interior of a $k$-dimensional cell
    $\delta$ is isomorphic to $(S^1)^k$.
  \end{itemize}

\demo{{\rm 5.6.} Mumford{\rm -}Faltings\/{\rm -}\/Chai\/{\rm '}\/s uniformization of abelian
  varieties}\hfill\break\noindent 
 Mumford-Faltings-Chai construction provides a uniformization of an
  abelian variety over a quotient field $\cK$ of a normal
  excellent ring $\cR$ complete with respect to an ideal $I$. We would
  now like to show that the data describing this uniformization are formally the same data as those for a stable
semiabelic variety of
  type $[c,c^t,\cM]$. These data contain an element of the group
  $Z^1(\Delta,\uhbX)$ as defined in Section~5.3. The complex $\Delta$ consists of one cell
  $\XR/Y$, where $X=X'$ is the group of characters for the toric part
  of the central fiber. The situation, therefore, is very similar to
  the case of field $\bC$.

\proclaimtitle{Mumford-Faltings-Chai}
\specialnumber{5.6.1}\proclaim{Theorem}
  There is a $1$\/{\rm -}\/to\/{\rm -}\/$1$ correspondence between the following\/{\rm :}
  \begin{itemize}
  \ritem{1.} A semiabelian scheme over $\cR$ together with an ample sheaf
    such that 
    \begin{itemize}
    \ritem{(a)} the generic fiber over $\cK$ is abelian{\rm ,}
    \ritem{(b)} the special fiber over $\cR/I$ is semiabelian with a split
      torus part{\rm .}
    \end{itemize}
  \ritem{2.} The degeneration data which contain an element $$\tau\in Z^1(\XR/Y,
    \uhbX[c,c^t,A,\cM])$$ as   defined in
    Section~{\rm 5.3,} with $c,c^t,A,\cM$
    over $\cR${\rm ,} and with the positivity condition\/{\rm :} the projection of
    $\tau$ to the functions on $\Sym^2Y$ has to be positive definite
    with respect to the powers of $I${\rm .}
  \end{itemize}

\endproclaim

\vglue-8pt
{\it Proof}.
  This is the reformulation of the content of Chapters II and III of
  \cite{FC}; see especially Theorems II.4.1 and
  II.5.1. The functions $a(y)$ and $\psi(y)$, there, in our notation are 
  the functions on $[y](1,0)$, and the functions $b(y,x)$ and
  $\tau(y,x)$ are the functions on $[y](0,x)$.
\hfill\qed\vglue4pt

{\it Remark} 5.6.2.
  Note that one has an element of $Z^1$ and not of $H^1$ because we
  look for an abelian variety over $\cK$ and not a torsor.

\demo{{\rm 5.7.} One\/{\rm -}\/parameter families of stable semiabelic pairs}
  In the next three sections we shall follow very closely the
  linearized toric case considered in 2.8, 2.9, 2.10 and do most things by analogy. We will
  also use the same or very similar notation. Hence, in this section
  $\cR$ denotes a DVR with the quotient field $\cK$
  etc. $(P_{\eta},\Theta_{\eta})$ is a stable semiabelic pair with the 
  action of a semiabelian group scheme $G_{\eta}$ over $\cK$.
\enddemo

{\it Reduction Step} 1.
  First, we recall that by the stable reduction theorem, after a
  finite base change $G_{\eta}$ can be extended to a semiabelian group 
  scheme $G$ over $\cR$. After an \'etale base change we can assume
  that the toric part of the central fiber is split with a character
  group $X\simeq\bZ^r$.
\vglue12pt

  One important difference from the purely toric case is this: if the
  toric parts of $G_{\eta}$ and $G_0$ have different dimensions,
  i.e.\ $G_{\eta}$ degenerates, then we assume $\cR$ to be
  {\it complete}. This does not hurt our chances of proving properness 
  of the moduli space: in the valuative criterion of properness the
  valuation rings considered can be assumed to be complete DVRs, in
  fact, even with an algebraically closed residue field.

\specialnumber{5.7.1}\proclaim{Theorem}
  After a finite{\rm ,} possibly ramified{\rm ,} base change the pair
  $(P_{\eta},\Theta_{\eta})$ 
  can be extended to a family $(P,\Theta)$
  of stable semiabelic pairs over $\cR${\rm .} Moreover{\rm ,} such an extension
  is unique{\rm ,} up to an isomorphism{\rm .}
\endproclaim

{\it Proof}.
   We will first prove this theorem in two special cases:
  \begin{itemize}
  \ritem{1.} when $G_{\eta}$ does not degenerate,
  \ritem{2.} when $G_{\eta}$ is abelian (the case considered in the
    Mumford-Faltings-Chai uniformization).
  \end{itemize}
  After this is done, the general case is obtained as a
  straightforward combination of the two.

\vglue6pt {\it Case} 1.
  In this case $G$ is a global over $\cR$ extension of a split torus
  $T$ by an abelian scheme $A$.
  
  As in Section~2.8, consider an even more special case first: when
  $(P_{\eta},\Theta_{\eta})\break \otimes\cK$ corresponds to a single
  polytope $Q\subset \XR$. After an \'etale base change we can put
  $(P_{\eta},\Theta_{\eta})$ in the standard form:
  $P_{\eta}[\delta,\cM_{\eta}]$. By properness of $A^t$ and
  $\Pic_{A/S}^{\lambda}$ over $S$ the sheaf $\cM_{\eta}$ can be
  extended to an ample sheaf $\cM$ over $A$. We have a $\cK$-algebra$$ 
    R_{\eta}=\oplus H^0(A,\cM_{\chi}) {\displaystyle\mathop{\otimes}_{\cR}}\cK.
  $$ 
  For each $x\in X$ we have an invertible $\cO_A$-module $\cM_{\chi}$
  over $A$. $\Gamma(A,\cM_{\chi})$ is a free $\cR$-module of rank $d$.
  Choose a basis $(e_1,\dots,e_d)$ in this module. The
  $(1,x)$-homogeneous component of $\theta$ can be written in this
  basis as$$ 
    \xi_x= \sum c_i(x) e_i^{(x)}.
  $$ 
  Now define $\psi(x)=\min\val\, c_i(x)$. In the same way as in  
  Section~2.8,
  the function $\psi$ defines
  \vglue2pt
   1.  a pointed cell decomposition $(\Delta_0,C_0)$ of $Q$,
  \vglue2pt 2.  an $\cR$-subalgebra $R$ in $R_{\eta}$.
  \vglue2pt\noindent 
  As before, after the necessary reductions, take $(P,L)=(\Proj_A R,
  \cO(1))$. The central fiber is a polarized variety
  $(P,L)[\Delta_0,\tau_0=1,c_0,c_0^t=1,\cM_0]$ and so is a stable
  semiabelic variety. The restriction of $\theta_0$ to each minimal
  orbit (i.e.\ corresponding to a vertex of one of the polytopes
  $\delta_i$) is not zero, so our condition on orbits is satisfied.
  Hence, we have obtained the required extension. The uniqueness
  follows by the next argument. We have a morphism to $A$ defined
  over $\cR$. Choose a section $\Spec\cR\to A$. Then over this section
  we have a family of stable toric varieties $P'$ and the type of the
  central fiber $P'_0$ is the same as the type of $P_0$. Hence, the
  uniqueness of the type of the central fiber, the corresponding
  subalgebra, and of the family follow by the toric case. This special
  subcase is completed.
  
  To complete the case 1, after an \'etale base change we represent
  $(P_{\eta},L_{\eta}),\Theta_{\eta}$ in the form
  $(P,L,\Theta)[c_{\eta},c^t_{\eta},\cM_{\eta},\tau_{\eta},\hf_{\eta}]$
  with $\tau_{\eta}$ being a collection of functions on
  $C_1/B_1(\ubX)$ with values in appropriate $\bG_m$-torsors. 
  In other words, $R_{\eta}$ is  a union of the algebras as in the
  special subcase, appropriately glued.
  
  Since $A/S$ is proper, $c^t_{\eta}:\pi_1\Delta\to A_{\eta}$ can be
  extended to $c^t:\pi_1\Delta\to A$ over the whole $\cR$. After this,
  the $\cR$ algebra $R$ is constructed formally by exactly the same
  formulas as in Section~2.8.
 
\vglue6pt {\it Case} 2.
  In this case, we have the Mumford-Faltings-Chai uniformization of
  $G_{\eta}= A_{\eta}$.  Again, we use formally the same formulas.
  However, the complex $\Delta$ now is $\XR/Y$ and the homologies are
  the group homologies with respect to the $Y$-action, as computed in
  Section~5.6. The integral-valued function $\psi(x)$ becomes a
  function defined on finitely many points $(1,x_i)$ --
  representatives of $X/Y$ in $X$.  Now, $\psi(x)$ is  propagated to an
  integral-valued function on the whole $X$ by the cocycle $\tau$
  according to the rule $\psi(x_i+y)=\psi(x_i)+\val\tau^{[y](1,x_i)}$.
  Recall that the function $\val\tau^{[y](1,x_i)}$ is a nonhomogeneous
  quadratic function of $y$ whose homogeneous part is positive
  definite. The convex hull $\ConvHull_{\psi}$ is an ``infinite
  polytope'' in the space $\XR\oplus \bR$. Because of the
  positive-definiteness condition, its lower envelope has finite
  faces. The projection of this lower envelope gives a $Y$-periodic
  decomposition of $\XR$ into lattice polytopes, i.e., a cell
  decomposition of $\XR/Y$.
  
  We can now formally define the algebra $\wR$  by the same formulas
  as before where $\wR$ is an $\cR$-subring of$$ 
    R_{\eta}=\oplus_{\chi=(d,x),x\in X,d>0} 
    H^0(A,\cM_{\chi})\otimes \cK.
  $$ 
  Although $R$ is not finitely generated over the  $\cO_A$-algebra itself,
  $\Proj_A R$ is covered by $\Spec_A $'s of finitely generated
  algebras. Hence, $(\wP,\wL)= (\Proj_A\wR,\cO(1))$ is a polarized
  scheme which is locally of finite type over $A$.
  
  After the reductions, the rest proceeds as in   Mumford's
  construction; see \cite{FC}, \cite{AN}. For the
  central fiber we have the proper variety $(P_0,L_0)=(\wP_0,\wL_0)/Y$
  and it is a stable semiabelic variety. It comes with a section
  $\theta_0$ which gives a divisor satisfying our condition of not
  containing orbits. For any infinitesimal neighborhood $\cR/I^{n+1}$
  we similarly have the pair
  $(\wP_n,\wL_n)=(\wP,\wL){\displaystyle\mathop{\otimes}_{\cR}} \cR/I^{n+1}$
  and the quotient $(P_n,L_n)=(\wP_n,\wL_n)/Y$ is proper over\break $\Spec\, \cR/I^{n+1}$. The formal scheme defined
by the collection
  $\{(P_n,L_n)\}$ is uniquely algebraizable by a fundamental theorem
  of Grothendieck to a polarized scheme $(P,L)$ over $\cR$. All the
  infinitesimal thickenings $\{(P_n,L_n)\}$ come with compatible
  sections $\theta_n$ which algebraize to a section $\theta$ of $L$.
  Finally, we recover $(P_{\eta},L_{\eta},\theta_{\eta})$ as the
  generic fiber of $(P,L,\theta)$. This proves the existence of the
  extension. 
  
  By   Grothendieck's formal algebraization theorem the family
  $(P,L,\Theta)$ is equivalent to a system $\{(P_n,L_n,\Theta_n)\}$
  over $\cR/I^{n+1}$. So, to prove the uniqueness, it is sufficient to
  prove the uniqueness of this system.  By
  Section~\ref{sec:Linearization of torus action} it defines a system
  $\{(\wP_n,\wL_n,\tilde\Theta_n)\}$ with a $Y$-action. By the finite
  case, this system uniquely defines a decomposition $\wDelta$ and the
  graded $\cO_A$-algebra $\wR=\oplus \tilde f_*\wL^d$ together with
  the $Y$-action. This gives the uniqueness.
\hfill\qed\vglue4pt

{\it Definition} 5.7.2.
  The cell decomposition $\Delta_0$ of $\XR/Y$ in the second case was
  defined, similarly to the finite case, by a lifting
  property. Accordingly, we will call such decompositions
  {\it regular}. 

\vglue4pt 5.7.3.
  We note one special case, when $Y=X$. A cell decomposition of
  $\XR/X$ in this case is the same as an $X$-periodic decomposition of
  $\XR$ with vertices in $X$. A {\it regular} decomposition in this
  case is the one defined by the following lifting property: it is the
  projection of the lower envelope of the convex hull of the points
  $\{(x,a(x)) \,|\, x\in X\}$, where $a(x)$ is a nonhomogeneous
  quadratic function with a positive-definite homogeneous part. Such
  decompositions have been known for at least a hundred years. They
  are usually called the {\it Delaunay} decompositions.
\vglue4pt

{\it Remark} 5.7.4.
  Note that if $\Delta_1$ and $\Delta_2$ are the complexes for the
  special, resp.\ generic, fibers then $|\Delta_1|$ and $|\Delta_2|$ can
  be identified, providing we are working over $\bZ[1/d]$ and the
  characteristic of the residue field $\cR/m$ does not divide $d$.

\vglue4pt

 5.8. {\it The modified complex} ${\bM}_*$.
  Let $(\Delta,C)$ be a pointed complex of lattice polytopes
  referenced by $\oXR/Y$. As usual, we represent it as the quotient
  $(\Deltad,\Cd)/Y$ for some $Y\subset X$.  We also choose
  representatives $\{\delta_i\}$ in $\Deltad$ of the maximal cells of
  $\Delta$; there are finitely many of them.  We will only construct
  the complex $S\bM_*$ and construct the moduli space in the next
  section in the case when $\rho:|\Deltad|\to\XR$ is injective and the
  image is convex.  An argument  similar to 2.10.1
  shows that this is the property which is constant in connected
  families. In this special case one has $\pi_1\Delta=Y$.
\enddemo

{\it Definition} 5.8.1.
  We will define the complex $SC_*(\ubX)$ in degrees $0,1,2$ in the
  same way as in the finite case but with the following modifications:
  \begin{itemize}
  \item[1.] The chain group $SC_1(\ubX)$ will as generators have symbols
    $(\bar\chi, \chi_{i_0})$, where $\bar\chi=\sum y_s.\chi_{i_s}$
    such that $\bar\chi$ and $\chi_{i_0}$ have the same image in
    $\Cone\, \Deltad\subset \XR$. We will also add generators
    $(y.\bar\chi,y.\chi_{i_0})=(\bar\chi, \chi_{i_0})$.
  \item[2.] The groups $\uFun_{\ge0,i}$ will be replaced as follows.  For
    each $x\in C_i$ in a polytope $\delta_i$ instead of just having
    the group $\bZ$ we will take $\bZ\oplus\bN^{d-1}$ ($d$ will be
    $h^0(\cM)$ for the ample sheaf on the abelian part). For each
    $x\in X\cap\delta_i\setminus C_i$ we will take $\bN^d$.
  \end{itemize}

\enddemo

  We should explain the form of the modified groups
  $\uFun_{\ge0,i}$. In each space $H^0(A,\cM_{\chi})$ we will choose a 
  basis. With such a basis, an element of this space is the same as a
  semigroup homomorphism from $\bN^d$ to the semigroup $k$ for
  multiplication. For $x\in C_i$ the condition  is that
  the section  not be $0$. The corresponding parameter space is
  $\bA^d\setminus 0$ and is not affine. As usual, we will cover it by
  $d$ affine sets and these will correspond to the semigroups
  $\bZ\oplus\bN^{d-1}$.

\specialnumber{5.8.2}\proclaim{Lemma}
  With these modifications{\rm ,} the semigroups $SC_1/SB_1(\ubX)$ and
  $S\bM_1/ B_1(S\bM_*)$ are finitely generated and finitely
  presented{\rm .} 
\endproclaim
 
\demo{Proof}
  Let us start with $SC_1/SB_1(\ubX)$. It follows from the relations
  that the group is generated by the symbols
  $(y_1.\chi_{i_1}+y_2.\chi_{i_2}, \chi_{i_0})$ with
  $\rho(y_1.\chi_{i_1})+\rho(y_2.\chi_{i_2})= \rho(\chi_{i_0})$ and
  with $\chi_{i_0},\chi_{i_0},\chi_{i_0}$ from a fixed finite set of
  generators for the corresponding $\bX\cap\Cone\,\delta$.  Choose an
  arbitrary norm on $\XR$. We claim that if some $y$  in such a
  generator is ``large'' with respect to this norm then it can be
  represented as a sum of generators with ``shorter'' $y$'s. (This
  statement is easy to quantify in terms of the sizes of $\delta_i$'s
  but we will not do it.) For the fixed $i_0,i_1,i_2,y_1,y_2$ we get
  semigroups of integral points in finitely generated rational cones,
  hence finitely generated. Therefore the whole group is finitely
  generated. As we had mentioned earlier, it is automatic that it is
  also finitely presented.
  
  To prove our claim, first of all take a very large $\chi'_{i_0}$ in
  $\Cone\,\delta_{i_0}$. Then
  \begin{eqnarray*}
     (y_1.\chi_{i_1}+y_2.\chi_{i_2}, \chi_{i_0})&=&   
    (y_1.\chi_{i_1}+y_2.\chi_{i_2}+2\chi'_{i_0}, 
    \chi_{i_0}+2\chi'_{i_0})\\
    & = & (y_1.\chi_{i_1}+\chi'_{i_0}, y'_1\chi'_{i_1})+
    (y_2.\chi_{i_1}+\chi'_{i_0}, y'_2\chi'_{i_2})\\
&&+\
    (y'_1.\chi'_{i_1}+y'_2.\chi'_{i_2}  , y_3\chi_{i_3}).
  \end{eqnarray*}
  I now claim that with the restraints on the sizes of $\chi$'s (they
  are from a finite set of generators) it is truly obvious that the
  sizes of $y_1,y_2$, if one of them is big, are roughly equal, and
  the sizes of the new $y'_1,y'_2,y_3$ are roughly half of that (just
  make a picture). This proves the statement for $SC_1/SB_1(\ubX)$.
  For $S\bM_1/ B_1(S\bM_*)$ note that with the relations that we have
  adding $\uFun$ adds several copies of $\bZ$ and $\bN$, not more than
  $d$ times the rank of $H_0(\uFun)$, a finite number.
\enddemo

\specialnumber{5.8.3}\proclaim{Lemma}\label{lem:main_action_finite2}
  The action of the torus $\Spec\, \bZ[\bSM_0]=\Spec\, \bZ[C_0(\ubX)]$ on
  the scheme $\Spec\, \bZ[\bSM_1/B_1(\bSM_*)]$ is proper with finite
  stabilizers{\rm .}
\endproclaim

\demo{Proof}
This is exactly the same as that of \ref{cor:main_action_finite}.
\enddemo

{\it Definition} 5.8.4.
  For each $\Delta'\ge\Delta$ and $(\Delta',C')\ge (\Delta,C)$ we have
  ideals $I[\Delta']$ and $I[\Delta',C']$ in $\bZ[\bSM_1/B_1(\bSM_*)]$
  defined as in the finite case.   

\demo{{\rm 5.9.} Test families of stable semiabelic pairs}
For each decomposition $\Delta'$ of $\oXR'/Y'$ with $\dim X'\le g$ we
are going to construct a complete Artin stack with a universal family
of $g$-dimensional stable semiabelic pairs over it. In the end,
this stack will be a completion of an open substack in our moduli
stack along the stratum corresponding to~$\Delta'$.
\enddemo

\hglue22pt A. {\it The totally degenerate case\/{\rm :} $\dim\Delta'=g$}.
Over $$U=\Spec\, \bZ[S\bM_1/B_1(S\bM_*)]$$ there exists a natural family
$(P^{\dag},\Theta^{\dag})_U$ of pairs with the $T$-action. The schemes
$P^{\dag}_U$ here are locally of finite type. These are would-be
$Y$-covers of the real family $(P,\Theta)_U$. The family
$(P^{\dag},\Theta^{\dag})_U$ is constructed in a natural way since the
semigroup $S\bM_1/B_1(S\bM_*)$ is precisely the semigroup of relations
one has to have in such a family. For each cell $\delta$ we consider
the semigroup algebras $\bZ[\Cone\, \delta]$ and $\bZ[\Cone\, y.\delta]$
for all its $Y$-translates. These cover the algebra graded by the
lattice $\bX$ with a representative in each degree $(d,x)$ with $d>0$.
A ring homomorphism $\bZ[S\bM_1]\to D$ is precise and allows us to glue the
rings $D[\Cone\, y.\delta]$ together and to give an element $\thetadag$
in degree 1. The condition for us to  actually get a ring is that this
homomorphism factors through $\bZ[S\bM_1/B_1(S\bM_*)]$.

Hence we have a ``universal'' $\bX$-graded algebra $\Rdag_U$ and we
define $(\Pdag,\Ldag)_U$ as $(\Proj_U \Rdag,\cO(1))$ which  comes with a
relative Cartier divisor $\Thetadag_U$ defined by the section
$\thetadag_U$ of $\Ldag_U$. By its very definition,
$(\Pdag,\Ldag,\thetadag)_U$ also comes with the $Y$-action, and our
next aim is to define the quotient. The map $\Pdag_U\to U$ is only
locally of finite type and many of the fibers are simply isomorphic to
$T$, so that this is just a bit delicate.

An ideal $I[\Delta]\subset \bZ[S\bM_1/B_1(S\bM_*)]$ defines
$\Stratum[\Delta]_U$ over which all fibers have type $\Delta$. Over
this stratum the $Y$-action is properly discontinuous in Zariski
topology, and we can divide the family $(\Pdag,\Thetadag)_U$ directly.
Now, if we are working over $\bC$, we can look at a small  neighborhood $V$ open in
classical topology  over which the $Y$-action on the
family is properly discontinuous. Over $V$ we can take the quotient
$(P,L,\theta)_U$ and it will be a proper family with an ample sheaf
and a section.

In general, the method is the same as the essential trick of  
Mumford's construction \cite{Mum1}. We look at the
thickenings$$ 
  U_n= \Spec\, \bZ[S\bM_1/B_1(S\bM_*)]/I[\Delta]^{n+1}  .
$$ 
Over each of these thickenings the $Y$-action on the family
$(\Pdag,\Ldag,\thetadag)_{U_n}$ is properly discontinuous in Zariski
topology, and we can take the quotient $(P,L,\theta)_{U_n}$. Each
sheaf $L_{U_n}$ is relatively ample, hence we have a compatible system
of projective polarized varieties. By   Grothendieck's formal
existence theorem \cite[5.4.5]{EGA} this family is equivalent to a
projective polarized family $(\hP,\hL)_U$ over the completion $\hU$.
The sections $\theta_{U_n}$ glue to a section $\widehat\theta$ of
$\hL$. Also, the family $(\Pdag,\Ldag,\thetadag)_{U_n}$ comes with a
choice of origins for each of the finitely many cells $\delta_{i_0}$.
Choose one of these sections, say $e$, and look at an open subscheme
of $\Pdag$ given by the condition that on every fiber $\Pdag_{s}$ it
is the union of components of the smooth locus containing $y.e(s)$.
This open subscheme is the $T$-orbit of $\{y.e(s)\}$ and it is a group
scheme over $U$. After completion and division by $Y$, it    gives a
semiabelian group scheme $\hG$ over $\hU$ acting on $\hP_U$.

Our last step is to divide $U$ and $\hU$ by the equivalence relation
$R$ given by choosing different origins and by a finite group of
symmetries of $\Delta$. On every $U_n$ this is given by a torus action
of $\Spec\, \bZ[S\bM_0] \times U_n$. This equivalence relation has
finite stabilizer since the equivalence relation $R$ on $U$ has finite
stabilizers by Lemma~\ref{lem:main_action_finite2}.  Finally, we
divide the family by the action of a finite group $\Sym\, \Delta$.  The
system $[U_n/R_n]$ gives a complete Artin stack which we will denote
$\cM[\Delta]$. It comes with a universal family $({\cal P},\Theta)$.

\medbreak  \hglue22pt B.  {\it The general case}.
We merely repeat the previous construction but over an appropriate
stack classifying the moduli of the abelian part.  We start with the
stack parametrizing semiabelian varieties with abelian parts of
constant rank $r$ and toric parts of constant rank $r$.  We take the
stack $A_{a,d}$ for the abelian part $(A,\lambda)$. Over it we take
$\bHom_{\bZ}'(Y'\times X',A\times A^t)$ which means the closed
subscheme of $\bHom_{\bZ}(Y',A)\times \bHom_{\bZ}(X',A^t)$ that
corresponds to our basic commuting condition $c\circ i= \lambda \circ
c$. This is the parameter space for the pairs $(c,c^t)$. Then we add a
direct summand $\bPic^{\lambda}$ parametrizing the rigidified sheaf
$\cM$.  Over this family we consider the functions on
$\bM_1/B_1(\bM_*)$ with the values in the appropriate torsors. These
are exactly what universally describes the graded algebra $\Rdag$ over
the abelian part together with a section $\thetadag$ of degree 1.  At
this point, we have a universal family of what would be ``universal
covers'' $(\Pdag,\Thetadag)$ together with the $Y'$-action. We finish
by applying the same completion / dividing by $Y'$ / algebraization
construction as before.
 
Finally, we divide the family by three equivalence relations:
\begin{itemize}
\item[1.] by the action of the torus $\Spec\,\bZ[\bM_0]$,
\item[2.] by the action of $A$ corresponding to changing the origin in
  the abelian part,
\item[3.] finally, by the symmetries of $\oXR'/Y'$ and the pairs
  $(A,\lambda)$.
\end{itemize}
The first of these equivalence relations is given by a proper action
with finite stabilizers.  This follows as before by Lemma~\ref{lem:main_action_finite2}. The second action is also
proper and with finite stabilizers. The third equivalence relation is finite.
The outcome is a stack  $\cM[\Delta']$ with a ``universal'' family of
stable semiabelic pairs with an arbitrary abelian part.
  \enddemo

5.10. {\it Moduli of stable semiabelic pairs}.

\specialnumber{5.10.1}\proclaim{Theorem}\label{thm:moduli_infperiodic_case}
{\rm 1.} The component $\overline{\cAP}_{g,d}$ of the moduli stack of
    stable semiabelic pairs containing $\cAPgd$ and pairs of the same
    numerical type is a proper Artin stack with finite stabilizers\/{\rm ;}
 \begin{itemize} \ritem{2.} It has a coarse moduli space $\overline{\AP}_{g,d}$ as a
    proper algebraic space.
  \end{itemize}  

\endproclaim

{\it Proof}.
  First, we claim that we can describe all infinitesimal deformations
  of a pair $(G_0,P_0,\Theta_0)/\Spec\, k$. Indeed, let $(G,P,\Theta)/S$
  be any infinitesimal family which has this pair as the central fiber, with
  local artinian $S$. Then the ranks of abelian and toric parts are
  constant, and $G/S$ is a global extension of $A/S$ by $T/S$. By
  taking an \'etale cover we can make the toric part to be split. By
  the results of Section \ref{sec:Linearization of torus action} we
  have, in a canonical way, a family $(G,\wP,\wTheta)$ with the sheaf
  $\wL=\cO(\wTheta)$ $T$-linearized. Let $(G,\Pdag,\Thetadag)$ be a
  connected component.  By choosing a section we have a projection $p$
  to the abelian part.  Any finite union of irreducible components
  $(G,P,\Theta)_f$ of $(G,\Pdag,\Thetadag)$ is a stable semi\-abelic
  pair. By \ref{thm:cohomologies_polarized_STVs} the higher direct
  images $R^ip_{f*}L_f$ vanish. We conclude that describing any finite
  part of $(G,\Pdag,\Thetadag)$ is equivalent to describing the graded
  $\cO_A$-algebra $R_f=\oplus_{d\ge0} L_f^d$ with an element
  $\theta_f$ of degree~1. Hence, describing the family
  $(G,\Pdag,\Thetadag)$ itself is equivalent to describing a graded
  $\cO_A$-algebra $\Rdag$ with an element $\thetadag$ of degree 1.
  Because of the universal character of our relations in
  $S\bM_1/B_1(S\bM_*)$ we see that the family $(P,\Theta)/S$ is the
  pullback of the universal family under a unique morphism from $S$ to
  the stack $\cM[\Delta]$.
  
  At this point, we are ready to apply Theorem 5.3 of
  [Art1]; cf.\ the case of surfaces of general type
  \cite[5.5]{Art2} there. Our groupoid of families of
  stable semiabelic pairs has a deformation theory and obstruction
  theory (even better than that) and for any point we have a
  pro-representable family of projective objects. Each pair
  $(P,\Theta)$ has a finite automorphism group. The openness of
  versality follows from the fact that for $\Delta'\le\Delta$ the
  semigroup model $[\Spec\,\bZ[S\bM_1]/B_1(S\bM_*)/ \Spec\, \bZ[\bM_0]]$
  for $\Delta'$ is an open substack of the semigroup model for
  $\Delta$.  Hence, all conditions of Artin's method are satisfied and
  we have an algebraic stack in   Artin's sense.  Our stack is
  complete by Section~5.7. Also, the formal models we have constructed are
  separated.  Since every 1-parameter family can be completed and the
  type of the central fiber in a 1-parameter family is uniquely
  defined, our stack is also proper. It has a coarse moduli space as a 
  proper algebraic space by [KM].\hfill\qed

\vglue12pt  5.11. {\it Approximation of the moduli space}.
Similarly to Section~2.11, we are now going to describe our moduli space in more detail
in the case of the complex $\oXR/X$, i.e.\ the component containing
principally polarized abelian torsors. We shall construct a
simplification of our moduli space and identify the main irreducible
component with a particular toroidal compactification of $A_g$. This we will do   in the totally degenerate case, over
the 0-cusp only. The general case at this point would not add anything but bookkeeping
complications.

In the lattice $C_0(\oXR/X,\uFun)/B_0(\oXR/X,\ubL)$, which by
Lemma~\ref{lem:homologies_of_L} is a lattice of dimension
$g(g+1)/2+(g+1)$, consider an infinite ``polytope'' $\sigma$ with  
vertices $1_x$, $x\in X$. This is an analog of a simplex from  the finite case in Section~2.11. There is a natural
projection
$\phi:\sigma\to (1,X)$ and the groups we are about to define will be constructed from the fibers of
this ``polytopal'' map.
\vglue12pt 
5.11.1.
  For every $q\in X_{\bR}$ we define the convex figure $\sigma_q$ as
  the fiber of $\sigma\to X$ over $(1,q)$. This is an infinite
  ``polytope'' of dimension $g(g+1)/2$ lying in the space parallel to
  $H^0(\oXR/X,\ubL)$ and, similar to the finite case, it has an
  explicit description as follows. Let $\delta$ be a $k$-dimensional
  ($0\le k\le g$) simplex fitting in $X_{\bR}/X$, i.e.\ whose
  $X$-translates intersect $\delta$ only at faces which contains $q$
  in its interior.  Say that $\delta$ has vertices $x_0,\dots x_k$ and
  $(1,q)=\sum c_i (1,x_i)$. Then there is a vertex of $\sigma_q$
  corresponding to $\delta$ which is an element of
  $C_0(\uFun)/B_0(\ubL)$ of the form$$ 
    \frac{1}{k+1} \sum c_i 1_{x_i}.
  $$ 
  Vertices of $\sigma_q$ correspond in a 1-to-1 way to simplices in
  $\oXR/X$ whose translation class contains $q$ in the interior, and
  higher-dimensional faces correspond to nonsimplicial cells with
  this property. The dual fan $\Fan(\sigma_q)$ divides $\Sym^2
  X^*_{\bR}$, i.e.\ the space of quadratic functions on $X_{\bR}$, into
  cones according to whether their Delaunay decompositions contain the
  same cell with $q$ in the interior. 
 
\demo{Definition {\rm 5.11.2}}
  For a quadratic function $\psi$ on $X$ consider the lower envelope
  of the points $(x,\psi(x))\in X_{\bR}\times \bR$ (this only makes
  sense if $\psi$ is semipositive definite). The projection of this
  lower envelope is called the Delaunay decomposition of $\psi$.
\enddemo

5.11.3.
  Now consider the Minkowski sum $\Sec(\oXR/X)$ of several
  $\sigma_{\alpha}$. The only condition is that we should have enough
  points $\alpha$ so that they can distinguish any two Delaunay
  triangulations according to which cells contain them (it is enough
  to consider the points appearing as centers of simplices, and these
  are among the finitely many, modulo $X$, rational points with
  denominator $(r+1)!$). Going through the same arguments with
  Minkowski averages as in Section~2.12, we see that the faces of $\Sec(\oXR/X)$ are in  
  1-to-1 correspondence with Delaunay decompositions, and that the
  dual fan $\Fan(\Sec(\oXR/X))$ divides the cone of semipositive
  definite symmetric quadratic forms into cones according to whether
  they define the same Delaunay decomposition. This is precisely the
  classical definition of the second Voronoi decomposition.  
 
\demo{Definition {\rm 5.11.4}}
  Now, let $\Delta$ be a triangulation of $\oXR/X$ into lattice
  polytopes. For every $q$ we can define a saturated semigroup
  $K_{\delta,q}$ in $H^0(\ubL)$ defined by differences $\chi_1-\chi_2$
  with $\chi_1\in\Cone\, \sigma$, $\chi_2\in \Cone\, \delta$,
  $\phi(\chi_1)= \phi(\chi_2) = N(1,q)$ for some cell $\delta$ of
  $\Delta$ with $q\in \delta^0$. This gives an affine $\Spec\, K_{\delta,q}$.
\enddemo

The group $K_{\delta,q}$ is one of the groups associated with the
infinite ``polytope'' $\sigma_q$. The advantage is that the former is
a semigroup of integral vectors in a finite rational polyhedral cone.
Note that for $x\in X$ the groups $K_{q,\delta}$ and
$K_{q+x,\delta+x}$ are naturally identified because$$ 
  (y.\chi_1-y.\chi_2)-(\chi_1-\chi_2) \in B_0(\ubL).
$$ 
As in the finite case, once we have chosen a decomposition $\Delta$,
there is a finite decomposition of $\oXR/X$ into polytopal regions
according to the type of a fiber $K_{\delta,q}$. Hence, we can define$$ 
  U[\Delta]= {\displaystyle\lim_{{\llar}}} \Spec\, \bZ[K_{\delta,\alpha}] =
  \Spec\, \bZ[{\displaystyle\lim_{\longrightarrow}}\,  K_{\delta,\alpha}].
$$ 

\specialnumber{5.11.5}\proclaim{Lemma}
  There is a finite morphism from the coarse moduli space $M[\Delta]$
  of the stack $\cM[\Delta]$ to the completion of $U[\Delta]$ along
  $\Stratum[\Delta]${\rm .}
\endproclaim

{\it Proof}.
  It is sufficient to check this on the semigroup models before the
  completion, for the morphism$$ 
   \Spec\bZ[\ker(S\bM_1/B_1(S\bM_*) \to \bM_0 ]
   \to U[\Delta].
  $$ 
  Here is how this map is constructed. Over
  $\Spec\bZ[S\bM_1/B_1(S\bM_*)]$ we have a family
  $(\Pdag,\Ldag,\thetadag)$. Write $\thetadag=\sum_X \xi_x$ as the
  infinite sum of its homogeneous components. For every finite
  collection $(n_x)\in \Fun(X, \bZ_{\ge0})$ we have an element $\prod
  \xi_x^{n_x}$, and for two such products in the same homogeneous
  space, we can compare them. The products with $(n_x)\in \Cone\,\delta$,
  for $\delta\in \Delta$, generate the corresponding 1-dimensional
  eigenspace $R_{\chi}$. Hence, we can write$$ 
    \prod \xi_x^{m_x} = c\left( (m_x),(n_x) \right) \prod \xi_x^{m_x}.
  $$ 
  Because $\Ldag$ and $\thetadag$ are invariant under the
  $Y=X$-action, one has$$ 
    c\left( y.(m_x),y.(n_x) \right) = c\left(
    (m_x),(n_x) \right)    .
  $$ 
  Hence, we have a coefficient $c$ defined for every pair of integral
  vectors in the same fiber of $\Cone(\sigma)\to \bX$.  As in the
  finite case, these coefficients $c$ define a homomorphism
  $K_{\delta,q}\to \bZ[S\bM_1/B_1(S\bM_*)]$ and by taking the limit we
  have the required homomorphism $\Spec\bZ[S\bM_1/B_1(S\bM_*)]\to
  U[\Delta]$.  Finally, note that this construction does not depend on
  the choice of origins, and so we obtain the required homomorphism.

  As in the finite case, both spaces have natural stratification by
  decompositions $\Delta'\ge \Delta$ and the strata, which are 
  respectively $\Spec\, \bZ[H^1(\Delta',\bM_*)]$ and $\Spec\, \bZ[
  H^0(\Delta', \ubL) ]$, are isogenic.
\hfill\qed\vglue8pt

As $\Delta$ goes over all triangulations of $\oXR/X$, the completions
of $U[\Delta]$ cover the simplification $M_{\simp}$ of our moduli
$\overline{\AP_g}$, and we obtain a finite morphism  $\overline{\AP_g} 
\to M_{\simp}$.

\specialnumber{5.11.6}\proclaim{Theorem}\label{saynum:main_comp_2ndVoronoi}
  The toroidal compactification of $A_g$ for the second Voronoi
  decomposition is isomorphic to the main irreducible component of
  $\overline{\AP}_g${\rm ,} the one containing
  ${\rm A}_{g}=\AP_{g}${\rm .}
\endproclaim

{\it Proof}.
  Again, we check this on the level of semigroup models, before the
  completion. The main component of $\overline{\AP}_g$ maps to
  $M_{\simp}$. On the other hand, for every Delaunay triangulation
  $\Delta$ there is a cone $\Cone[\Delta]$ in  
  $\Fan(\Sec(\oXR/X))$, i.e.\ the second Voronoi decomposition with a
  natural {\it surjective} map ${\displaystyle\lim_{\longrightarrow}}\,  K_{\delta,q}\to
  \Cone[\Delta]$.\pagebreak
  
  This means that the main component of $M_{\simp}$ is already normal
  and is isomorphic to the toroidal compactification of $A_g$ for the
  second Voronoi compactification. We have a finite map
  $\overline{\AP}_g\to M_{\simp}$ which is an isomorphism on the main
  stratum, the one corresponding to $\oXR/X$. Hence the main
  irreducible component of the space $\overline{\AP}_g$ is the same.
\hfill\qed\vglue8pt

5.12. {\it Generalized secondary polytopes for periodic
  decompositions}.
 For a polytopal decomposition of $\oXR/Y$ (where $Y\subset X$ has
finite index) the definition of the secondary polytope $\Sec(\Delta)$
is a straightforward generalization of that for the finite case which
we considered in Section~2.12:

\demo{Definition {\rm 5.12.1}}
  $\Sec(\Delta)$ is defined as the image of $\prod \Sec(\delta_i)
  \subset C_0(\Delta,\uFun_{\bR})$ in
  $C_0(\Delta,\uFun_{\bR})/B_0(\Delta,\ubL_{\bR})$, where $\delta_i$
  go over maximal-dimensional (i.e.\ minimal by the reverse of
  inclusion order) cells. 
\enddemo

{\it Example} 5.12.2.
  Let $E_n$ be the decomposition of $\bR^n$ into standard Euclidean
  cubes. There is only one maximal-dimensional cell modulo
  translation. Then $\Sec(E_n)$ is simply the image of the secondary
  polytope of the cube in a lower-dimensional space, of dimension
  $n(n-1)/2$ instead of $2^n-n-1$.
\vglue12pt

Similar to Section~2.12,
faces of this secondary polytope are in   1-to-1 correspondence with
regular subdivisions of $\Delta$. 

For an arbitrary decomposition of $\oXR/Y$ we can define
$\Sec(\Delta)$ as an appropriate Minkowski sum of ``polytopes'' some
of which are infinite, just as we have done   in the previous section
for $\oXR/X$. This gives an ``unbalanced'' version of $\Sec(\Delta)$.
Below we construct a balanced, $\GL(X)$-invariant version
combinatorially equivalent to the above, and giving the same dual
fan. 

\vglue12pt 5.12.3.
  For every $q\in X_{\bR}$ we define the balanced convex figure
  $\delta_q\bal$ as follows. Let $\delta$ be a $k$-dimensional ($0\le
  k\le r$) simplex fitting in $X_{\bR}/X$, i.e.\ whose $X$-translates
  intersect $\delta$ only at faces which contain $q$ in their interiors.
  Consider all the $k+1$ translates $\delta^0,\dots ,\delta^k$ of it
  that contain the origin $0$, and let $q_0,\dots, q_k$ be the
  corresponding translates of $q$ in their interiors. Say, $\delta^j$
  has vertices $x^j_0,\dots, x^j_k$ and $(1,q_i)=\sum c_i (1,x_i^j)$.
  Then define the vertex of $\sigma_q$ corresponding to $\delta$ as
  the image in $C_0(\uFun)/B_0(\ubL)$ of$$ 
    \frac{1}{k+1} \sum c_i 1_{x_i^j}.
  $$ 
\vglue12pt

  Taking the Minkowski sum of these balanced ``polytopes'' gives the
  explicit ``balanced'' coordinate definition for the secondary
  polytope obtained by   averaging the coordinates of the vertices
  above, similar to the equation~(2) in Definition 2.12.4.

\demo{Definition {\rm 5.12.4}}
  For each triangulation $\Tr$ of $\oX_{\bR}/X$ let
  $\phi_{\Tr}$ be an element of $C_0(\uFun)/B_0(\ubL)$ defined by the
  formula
  \begin{equation}\label{eqn:vertices_of_secpol_XR}
    \phi_{\Tr}= 
      \sum_{\delta_i: x\in\delta_i,0\in\delta_i} 
      \Vol\, \delta_i \cdot (1_x - 1_0). \speqnu{4}
  \end{equation}
  Then, for any $x\ne0$,
  \begin{equation}
    \phi_{\Tr}(x)= 
      \sum_{\delta_i: x\in\delta_i,0\in\delta_i} ,
      \Vol\, \delta_i,  \speqnu{5}
  \end{equation}
  and for the origin we have $\phi_{\Tr}(0)= -r(r+1)!$, which comes
  from $r=$ number of vertices in a simplex minus one, and $(r+1)!=$
  volume of $\Star(0)$ in the lattice measure of volume.

  Define $\Sec(\oX_{\bR}/X)$ as the convex hull of the points
  $\phi_{\Tr}$ as $\Tr$ goes over all the triangulations of
  $\oX_{\bR}/X$. This is a ``polytope'' which may (and indeed does)
  have infinitely many  vertices. 
\enddemo

\specialnumber{5.12.5}\proclaim{Lemma}
  $\Sec(\oX_{\bR}/X)$ is invariant with respect to the action of
  $\GL(X)${\rm .}
\endproclaim

\demo{Proof}
  Indeed, for any $g\in\GL(X)$ and any triangulation $\Tr$ one has 
  $g.\phi_{\Tr}=\phi_{g.\Tr}$.
\enddemo

\specialnumber{5.12.6}\proclaim{Lemma}
  $\Sec(\oX_{\bR}/X)$ lies in the subspace
  $H_0(\ubL)=C_0(\ubL)/B_0(\ubL)${\rm .}
\endproclaim

\demo{Proof}
  We need to check that the vertices $\phi_{\Tr}$ in $C_0(\uFun)$ have
  the same image when projected to $C_0(\ubX)=\bX$.  Consider one  
  such vertex corresponding to a triangulation $\Tr$. Let $\sigma_1$
  be one of the maximal-dimensional simplices in this triangulation.
  In the translation equivalence class of $\sigma$ there are precisely
  $r+1$ simplices that contain $0$: they are obtained by shifting one
  of the $r+1$ vertices to the origin. For every edge $\{0,x\}$ in one 
  of these $r+1$ simplices there is precisely one edge $\{0,-x\}$ in
  its copy shifted by $(-x)$. Therefore, all terms in the definition
  of $\phi_{\Tr}$ split into groups of the form $(1_x+1_{-x}-2\cdot
  1_0)$, and each of these groups has image $0$ in $\bX$.
\enddemo

\proclaimtitle{of the proof}
\specialnumber{5.12.7}\proclaim{{C}orollary}
  When considered as an element of $\Sym^2 X${\rm ,}
  \begin{equation}
    \phi_{\Tr}= \frac{1}{2}
      \sum_{\delta_i: x\in\delta_i,0\in\delta_i} 
      \Vol\,\delta_i \cdot x\otimes x. \speqnu{6}
  \end{equation}
\endproclaim

\demo{Proof}
  Indeed, under the isomorphism $H_0(\ubL)\isoto \Sym^2 X$, the
  element $(1_x + 1_{-x}-2\cdot 1_0)$ goes to $x\otimes x$.
\enddemo

\specialnumber{5.12.8}\proclaim{{C}orollary}\label{cor:2ndVor_projective}
  The toroidal compactification over $\bZ$ of $A_{g}$ corresponding to 
  the second Voronoi decomposition is projective{\rm .}
\endproclaim

\demo{Proof}
  This follows by   Tai's criterion in \cite{AMRT} which requires a
  $\GL(X)$-invariant polarization function on the fan. In our case it
  is provided by $\Sec(\XR/X)$. Even though \cite{AMRT} works over
  $\bC$, the result applies in every case when the Satake-Baily-Borel
  compactification is constructed as a projective scheme, which is
  done over $\bZ$ in \cite{FC}. See also \cite[IV]{Cha} 
  where   Tai's criterion is explained over $\bZ[1/2]$.
\enddemo

5.13. {\it Periodic non\/{\rm -}\/Delaunay decompositions}.
  It is very easy to construct nonregular periodic decompositions
  when the degree of polarization is big: simply take any finite
  nonregular decomposition and repeat it periodically filling the
  holes in an arbitrary periodic pattern.  The lifting property will
  not be satisfied even locally. The principally polarized case is
  significantly harder. It can be checked by hand that in dimensions
  $g\le3$ all periodic decompositions are regular, i.e.\ Delaunay.
  There is an example in dimension 4, however.

\demo{Example {\rm 5.13.1}}
  Erdahl and Ryshkov in \cite[3.4]{ER} describe a ``red
  triangular'' lattice of dimension 4 which has the following Delaunay
  decomposition $\Delta$: two cyclic polytopes $C(4,6)$ (of dimension
  4, with 6 vertices) and 18 simplices $\alpha_4$. The two $C(4,6)$'s
  are mirror-symmetric with respect to  a vertex.
  
  We shall  construct a subdecomposition of $\Delta$ which is
  not Delaunay itself.  
  If each cell is simplicial, so that  the proper faces cannot be subdivided
  any further, the secondary polytope $\Sec\Delta$ is
  $\prod\Sec(\delta_i)$, the direct product of the secondary polytopes
  of the cells.
  
  That is the case here, because a cyclic polytope is simplicial.
Now,  $\Sec(\alpha_4)$ is a point and $\Sec(C(4,6))$ has dimension
  $6-4-1=1$ and is a closed interval.  The end points correspond to
  two possible splittings of $C(4,6)$ into three simplices.
  Therefore, $\Sec(\Delta)$ is a square and the corresponding toric
  scheme is $P_{\Sec(\Delta)}=\bP^1\times\bP^1$.
  
  There are exactly six nonregular subdecompositions of $\Delta$. They
  correspond to splittings of the cyclic polytopes which are not
  mirror symmetric. They are not Delaunay because every Delaunay
  decomposition has the symmetry $x\mapsto -x$.
\enddemo

  The moduli interpretation of this example is that the connected
  component of $\AP_4$ containing $A_4$ (the moduli of PPAVs of
  dimension 4) also contains another irreducible component. It is
  2-dimensional.

\demo{Example {\rm 5.13.2}}
  There is a very elegant description of Delaunay decompositions of
  root lattices $A_n,D_n,E_{6,7,8}$ and their duals (weight lattices)
  in terms of their Coxeter-Dynkin diagrams: it was given in
  \cite{CS}, \cite{MP}; see also
  \cite{Erd}. In particular, the Delaunay decomposition of the
  $E_8$ lattice consists of 135 crosspolytopes $\beta_8$, each
  generating an index-2 sublattice, and 1920 triple-volume simplices
  $\alpha_8$.
  
  Since crosspolytopes (see \cite[7.2]{Cox}) are
  simplicial, the regular sub\-decom\-positions of $E_8$ are described by
  the faces of a polytope $\Sec(E_8)$ which is again the direct
  product of the secondary polytopes of the cells.
  
  It can be easily checked that the secondary polytope of an
  $n$-dimensional crosspolytope is an $(n-1)$-dimensional simplex:
  $\Sec(\beta_n)=\alpha_{n-1}$, and that the corresponding secondary
  variety is the projective space $\bP^{n-1}$. Therefore,
  $\Sec(E_8)$ is the direct product of 135 copies of $\alpha_7$, and
  $P_{\Sec(E_8)}$ is the direct product of 135 copies of $\bP^7$.
\enddemo

  The moduli interpretation of this example is that the connected
  component of $\AP_8$ containing $A_8$ (the moduli of PPAVs of
  dimension 8) also contains another irreducible component of
  dimension 945. Recall that $A_8$ itself (moduli space $A_8$, not the
  lattice $A_8$) has dimension $36$.

 \section{Further questions}

We list just some of the many questions that arise in connection to
our moduli spaces:
\begin{itemize}
\item[1.] How should one define level structures for semiabelic pairs that 
  behave well at infinity?
\item[2.] In the case of the  nonprincipal polarization  dimension of the main component of $\APgd$, is 
  $(d-1)$ larger than $\dim A_{g,d}$? Is there a simple way to cut   this
  dimension back to $\dim \Agd$ by considering not all pairs but
  perhaps only some symmetric ones, possibly on a finite cover
  corresponding to a level structure that behaves well at infinity?
  The canonical level structure seems to be a natural choice here.
  What is the right definition of the canonical level structure for
  the pairs?
\item[3.] Is there a way to isolate $\avor$ in $\overline{\APg}$ using a
  log structure? Is there a good definition for the moduli of pairs
  with a log structure?
\item[4.] According to Mumford and Namikawa \cite{Nam1}
  the Torelli map $M_g\to A_g$ extends to a morphism from the
  Mumford-Deligne compactification $\overline{M}_g$ to $\avor \subset
  \overline{\APg}$. One would expect that with our moduli
  compactification there must be a functorial map $\overline{M}_g \to
  \overline{\APg}$. We show that this is indeed the case in
  \cite{Ale1}.
\item[5.] Are the moduli spaces we have constructed always reduced? I
  suppose that the answer is no.
\item[6.] Are they connected, once we fix the numerical type? Again, I
  think the answer in general may be no. (A recent example of
  F.~Santos \cite{San} shows that the answer is in general no for
  $M_Q$, where $Q$ is a lattice polytope.)
\item[7.] The last question is directly related to the generalized Baues
  problem in geometric combinatorics. We hope that the geometric
  constructions of the present paper, in particular, the configuration
  of the generalized secondary polytopes that we have constructed, may
  provide new clues for the resolution of this problem.
\item[8.] Explicit description of various secondary polytopes describing
  parts of our moduli spaces is an interesting combinatorial problem
  in itself.
\item[9.] We only considered the semiabelian group action in this
  paper. What about other algebraic groups?
\end{itemize}

\AuthorRefNames [SGA71]

\bye